\documentclass[10pt,a4paper,reqno]{amsart}
\usepackage{CJKutf8}
\usepackage{tabls} 
\usepackage{url}
\usepackage[linktocpage = true]{hyperref}
\hypersetup{
 colorlinks = true,
 citecolor = blue}
\usepackage{color}
\definecolor{red}{rgb}{1,0,0}
\definecolor{green}{rgb}{0,1,0}
\definecolor{blue}{rgb}{0,0,1}
\definecolor{refkey}{gray}{.625}
\definecolor{labelkey}{gray}{.625}
\usepackage[lite]{amsrefs}
\usepackage[a4paper]{geometry}
\usepackage{amsfonts}
\usepackage{amssymb}
\usepackage{mathrsfs}
\usepackage{amsmath}
\usepackage{amsthm}
\usepackage{amscd}
\usepackage{stmaryrd,pigpen}
\usepackage{enumerate}
\usepackage{paralist}
\usepackage{xypic}
\usepackage{graphicx}
\usepackage{mathtools}
\usepackage[all,cmtip]{xy}
\usepackage{kantlipsum}
\usepackage{tikz}
\usepackage{tcolorbox}
\allowdisplaybreaks

\newcommand{\At}{\operatorname{At}}

\newcommand{\id}{\operatorname{id}}
\newcommand{\R}{\mathbb{R}}

\newcommand{\g}{\mathfrak{g}}
\newcommand{\h}{\mathfrak{h}}
\renewcommand{\k}{\mathfrak{k}}

\newcommand{\K}{\mathbb{K}}

\newcommand{\bsection}[1]{\Gamma(#1)}

\makeatletter
 \def\title@font{\normalsize\bfseries}
 \let\ltx@maketitle\@maketitle
 \def\@maketitle{\bgroup%
 \let\ltx@title\@title%
 \def\@\title{\resizebox{\textwidth}{!}{%
  \mbox{\title@font\ltx@title}%
 }}%
 \ltx@maketitle%
 \egroup}
\makeatother

\theoremstyle{plain}
\numberwithin{equation}{section}

\newtheorem*{theorem*}{Theorem}

\newtheorem{Def}[equation]{Definition}

\newtheorem{Ex}[equation]{Example}

\newtheorem{Convention}[equation]{Convention}
\newtheorem*{Que*}{Objective}

\newtheorem{defn}[equation]{Definition}
\newtheorem{lem}[equation]{Lemma}
\newtheorem{prop}[equation]{Proposition}
\newtheorem{thm}[equation]{Theorem}
\newtheorem{cor}[equation]{Corollary}
\newtheorem{rmk}[equation]{Remark}

\begin{document}
\def\ArtA{\mathscr{A}}
\def\ArtB{\mathscr{B}}
\def\ArtC{\mathscr{C}}
\def\bott{\mathrm{Bott}}
\def\bp{\begin{proof}}
\def\ce{\mathrm{CE}}
\def\C{\mathbb{C}}
\def\CE{\mathrm{CE}}
\def\D{\mathcal{D}}
\def\ev{\mathrm{ev}}
\def\F{\mathscr{F}}
\def\g{\mathfrak{g}}
\def\H{\textbf{H}}
\def\j{{J}}
\def\L{\mathscr{L}}
\def\E{\mathscr{E}}

\def\M{\mathcal{M}}
\def\m{\mathfrak{m}}
\def\t{\mathfrak{t}}
\def\O{\mathcal{O}}
\def\P{\mathcal{P}}
\def\r{\mathcal{R}}
\def\U{\mathcal{U}}
\def\w{\mathfrak{W}}
\def\L{\mathcal{L}}
\def\X{\mathbb{X}}
\def\Y{\mathbb{Y}}
\def\spec{\text{spec}}
\def\Im{\text{Im}}
\def\coker{\operatorname{coker}}
\def\Ext{\operatorname{Ext}}
\def\End{\operatorname{End}}
\def\pr{\operatorname{pr}}
\def\id{\operatorname{id}}
\def\Der{\operatorname{Der}}
\def\Hom{\operatorname{Hom}}
\def\Jet{\operatorname{Jet}}
\def\Map{\operatorname{Map}}
\def\Mod{\operatorname{Mod}}
\def\sgn{\operatorname{sgn}}
\def\sh{\operatorname{sh}}
\def\Ad{\operatorname{Ad}}

\newcommand{\Atiyahanalogue}{S^\nabla}
\newcommand{\zeroaction}{\chi}

\newcommand{\Linfty}{ L_\infty  }
\newcommand{\curvature}{\mathcal{R}}
\newcommand{\piepartial}{\eth}

\newcommand{\inputvariable}{\cdot}

\newcommand{\Artinring}{\mathfrak{a}}

 \newcommand{\MCgA}{\textrm{MC}_\g(\v)}
\newcommand{\MCLPgA}{\textrm{MC}_{\Omega^\bullet_A(B)}(\v)}
\newcommand{\MCfunctor}{{\rm MC}_{\Omega^\bullet_A(B)}}

\newcommand{\algDefFctLA}{{\rm algDef}_{{(L,A)}}}
\newcommand{\weakDefFctLA}{{\rm wkDef}_{{(L,A)}}}

\newcommand{\firstDef}{{\rm Def}_\g}
\newcommand{\MCfunctorg}{{\rm MC}_\g}
\newcommand{\diffE}{\mathfrak{D}(E)}
\newcommand{\LdiffE}{\mathfrak{D}_L(E)}
\newcommand{\LAdiffE}{\mathfrak{D}_{L/A}(E)}

\newcommand{\pairing}[2]{\langle #1,#2\rangle}
\newcommand{\rank}{\mathrm{rank}}
\newcommand{\dA}{d_A}
\newcommand{\dAE}{d^{\End(E)\oplus B}_{A}}
\newcommand{\dEndE}{d^{\mathend(E)}_{A}}

\newcommand{\Linftythree }{L_{\leqslant 3}}
\newcommand{\infDefFctLA}{{\rm infDef}_{{(L,A)}}}
\newcommand{{\newboundary}}{\kappa}
\newcommand{\Diff}{\mathrm{Diff}}
\newcommand{\CinfMK}{C^\infty(M,\K)}

\newcommand{\projection}{\mathrm{pr}}
\newcommand{\projectionA}{\mathrm{pr}_A}
\newcommand{\projectionB}{\mathrm{pr}_B}

\def\iA{i_{_A}}
\def\iB{i_{_B}}
\def\prA{\pr_{_A}}
\def\prB{\pr_{_B}}
\def\sigA{\sigma_{_A}}
\def\sigB{\sigma_{_B}}
\def\sigL{\sigma_{_L}}
\def\ioA{{\iota}_{_A}}
\def\ioB{\iota_{_B}}
\def\ioL{{\iota}_{_L}}
\def\rhoA{{\rho}_{_A}}
\def\rhoB{\rho_{_B}}
\def\rhoL{{\rho}_{_L}}

\def\fD{\mathfrak{D}}
\def\fk{\mathfrak{k}}
\def\sA{\mathscr{A}}
\def\sJ{{J}}
\def\rd{\mathrm{d}}

\newcommand{\longto}{\longrightarrow}
\newcommand{\longhookrightarrow}{\lhook\joinrel\longrightarrow}
\newcommand{\longhookleftarrow}{\longleftarrow\joinrel\rhook}

\newcommand{\llparen}{(\!(}
\newcommand{\rrparen}{)\!)}

\newcommand{\parbar}{\bar{\partial}}
\newcommand{\nabplus}{\scalebox{0.8}{$\stackrel{\nabla}{\oplus}$}}
\newcommand{\splus}{\scalebox{0.8}{$\stackrel{s}{\oplus}$}}
\newcommand{\tilnab}{\widetilde{\nabla}}
\newcommand{\tilnabplus}{\scalebox{0.8}{$\stackrel{\tilnab}{\oplus}$}}
\newcommand{\barnab}{\overline{\nabla}}
\newcommand{\barnabplus}{\scalebox{0.8}{$\stackrel{\barnab}{\oplus}$}}

\newcommand{\AdiffE}{\mathfrak{D}_{A}(E)}
\newcommand{\BdiffE}{\mathfrak{D}_{B}(E)}
\newcommand{\BdiffI}{\fD_{i_{_B}}(E)}
\newcommand{\AtExtE}{\End(E) \nabplus B}
\newcommand{\AtAlgE}{\sA_L(E,\nabla)}

\newcommand{\deltap}{\delta}
\newcommand{\deltas}{s}
\newcommand{\KScorrespondence }{\mathrm{KSc}}
\newcommand{\KSmap}{\mathrm{KS}}
\newcommand{\ArtCat}{\textbf{Art}}
\newcommand{\PhiA}{\Phi_A}
\newcommand{\PhiB}{\Phi_B}
\newcommand{\adLb}{{\mathrm{ad}_b  }}
\newcommand{\adLa}{{\mathrm{ad}_a  }}
\newcommand{\maximalidealofartin}{\m_\v}
 \newcommand{\bigiota}{\mathcal{I}}

\newcommand{\SpecofA}{\mathbb{P}}

\newcommand{\FINISH}{\textcolor{red}{FINISH}}

\newcommand{\unifAtExt}{\mathfrak{D}_{L/A}(E)}
\newcommand{\mathend}{\mathrm{End}}

\title[Atiyah constructions for Lie algebroid connections on fiber bundles]{Atiyah constructions for Lie algebroid connections on fiber bundles}


\author{Chen He}
\address{North China Electric Power University, Beijing, China}
\email{\href{che@ncepu.edu.cn}{che@ncepu.edu.cn}}

\author{Dadi Ni$^\diamond$}
\address{School of Mathematics and Statistics, Henan University} 
\email{\href{mailto:nidd@henu.edu.cn}{nidd@henu.edu.cn}}

\author{Zhuo Chen}
\address{Department of Mathematics, Tsinghua University} 
\email{\href{mailto:chenzhuo@mail.tsinghua.edu.cn}{chenzhuo@mail.tsinghua.edu.cn}}

\thanks{C.~H. is supported by the Natural Science Foundation of
	Beijing (No.~1244049) and the Fundamental Research Funds
	for the Central Universities of China (No.~2025MS082). D.~N. is supported by the Natural Science Foundation of Henan Province (No.~252300421766).}

\thanks{$^\diamond$ Corresponding author.}

\begin{abstract} 
	To address the need for a unified framework that incorporates Lie algebroid connections on both vector and principal bundles,
	this paper investigates a generalized Atiyah algebroid structure and its short exact sequence.
	Building on this generalization, we describe Atiyah-type extensions and sequences that represent Atiyah classes through three explicit constructions designed to encode Lie algebroid connections compatible with specified sub-structures.
	As illustrative examples, we work out an enriched Atiyah algebroid construct,
	providing a systematic tool to characterize certain key properties of holomorphic connections and invariant connections.
\end{abstract}

\maketitle

\begin{itemize}
	\item 
	\item 
\end{itemize}

\setcounter{tocdepth}{2}
\tableofcontents
\parskip = 0em

\section{Introduction}

\vskip 10pt

\subsection{The motivation}

In differential geometry, the theory of connections provides a fundamental tool for tensor calculus on smooth manifolds. 
Ehresmann introduced the modern definition of connections on a fiber bundle in terms of horizontal distributions, significantly influencing Kobayashi and Nomizu's comprehensive treatise~\cite{KN} on the foundations of differential geometry.
Atiyah~\cite{Atiyah57} provided an alternative definition of connections on a principal $G$-bundle $P$ over a manifold $M$. This approach utilizes vector-bundle splittings of the \textbf{Atiyah sequence}, which afterwards bears his name: 
\begin{equation*} 
	\xymatrix@C=2pc@R=2pc{
		0\ar[r] & \Ad(P) \ar[r] & {\At(P)} \ar@<0.5ex>[r] & TM \ar[r]\ar@<0.5ex>@{-->}[l] & 0.
	}
\end{equation*}
Here, $\At(P):=(TP)/G$ represents a vector bundle that encodes the connections on $P$. 

Subsequently, Libermann~\cite{Libermann} noted that Atiyah's construction, $\At(P)$, possesses a Lie-theoretic structure, now recognized as a Lie algebroid (refer to Section~\ref{Subsec: LieAlgoid}).
This structure gives $\At(P)$ the trademark of the \textbf{Atiyah algebroid} associated with the principal bundle $P$ (refer to Section~\ref{subsubSec:expansionLiealgebroid}).
Notably, Atiyah algebroid $\At(P)$ served as a crucial example in both Mackenzie's~\cite{Mac87}  and Kubarski's~\cite{Kubarski} investigations of transitive Lie algebroids.
Recently, the studies of Atiyah algebroid and sequence have various aspects from global analysis on manifolds~\cite{GraKP11}, non-commutative geometry~\cite{LazM12}, and Lie groupoid and stack~\cite{BisN14}, to name a few.

Atiyah algebroids also arise naturally in theoretical mechanics. Lagrangian functions on Lie algebroids are important objects in both mathematics and theoretical physics.
Specifically, when a Lagrangian defined on a principal $G$-bundle $P$ (a configuration space possessing symmetries under the action of the structural Lie group $G$) remains invariant under this action, it is inherently defined on the Atiyah algebroid $A=(TP)/G$.
This foundational observation, as reviewed by Weinstein~\cite{Wei96a} in the context of unifying internal and external symmetry, and further explored in the groupoid approach to mechanics~\cite{Wei96b}, has spurred the development of multiple Lie algebroid generalizations of Lagrangian and Hamiltonian formalisms~\cites{Mar01, LMM05, Mes05, GGU06}.
 
Note that the classic theory of connections has been mostly set up for connection operators with respect to tangent fields.
Nowadays, compared with tangent bundles, Lie algebroids sometimes emerge as more accurate representations of the underlying (potentially singular) geometric structures.
The concept of connections with respect to Lie algebroids, to the best of our knowledge, has been first explored in two primary contexts.
First, Xu~\cite{Xu1999} considered Lie algebroid connections on vector bundles, defining them in terms of covariant differential operators (see Definition~\ref{def:L-Conn}), which generalize Lie algebroid representations on vector bundles as described by Mackenzie~\cite{Mac87}.
Second, Fernandes~\cite{Fernandes02} studied Lie algebroid connections on principal bundles, characterizing them through horizontal distributions (see Definition~\ref{def:L-ConnP}).
 
It is now natural to investigate the Atiyah-type constructions for Lie algebroid connections, which is the main goal of this paper.


\vskip 10pt

\subsection{Objectives}
Building on the preceding findings and motivation, we outline three key considerations, states our objectives, and briefly describes our proposed solutions. 

\vskip 5pt
\subsubsection{Generalized Atiyah algebroids} 

Mackenzie~\cite{Mac87} addressed connections for a transitive Lie algebroid $L$ with surjective anchor map onto $TM$ by considering a splitting of the short exact sequence: 
\begin{equation*} \xymatrix@C=2pc@R=2pc{ 0\ar[r]& \ker \rho \ar[r]& {L} \ar@<0.5ex>[r]^(0.45){\rho}& TM \ar[r]\ar@<0.5ex>@{-->}[l]& 0. } \end{equation*} 
Kubarski~\cite{Kubarski} also employed this approach in his theory of characteristic classes for regular Lie algebroids.

Inspired by Mackenzie and Kubarski's perspective, we aim to construct an Atiyah-type sequence, such that a Lie algebroid connection (without the transitivity restriction) can still be understood as a splitting of this sequence.

\textbf{Objective 1} -- \textit{To develop a generalized version of Atiyah algebroids/sequences capable of representing Lie algebroid connections on both vector bundles and principal bundles.}

The classic Atiyah algebroids associated with a vector bundle $E$ or a principal $G$-bundle $P$ – constructions that we will revisit in Section~\ref{subsubSec:expansionLiealgebroid} – are Lie algebroid extensions of the tangent bundle $TM$. Expanding upon this foundation, modern research has revealed numerous applications for connections of Lie algebroids more general than $TM$, particularly at the intersection of Lie theory and Poisson geometry. In this paper, we generalize the notion of Atiyah algebroids and Atiyah sequences by replacing the role of $TM$ with a general Lie algebroid $L$. Specifically, given a Lie algebroid $L$ and a vector bundle $E$ over a common base manifold, we define $\LdiffE$, the \textbf{$L$-Atiyah algebroid} of $E$, characterized by having its symbols residing in $L$. Analogously, given a principal $G$-bundle $P$, we introduce $\At_L(P)$, the $L$-Atiyah algebroid of $P$ (see Definitions~\ref{def:LdiffE} and~\ref{def:LdiffP}). Furthermore, we establish a bijective correspondence (Theorems~\ref{thm:conn=splitting} and~\ref{thm:Pconn=splitting}) between Lie algebroid connections and vector-bundle splittings of these generalized Atiyah sequences: 
\begin{equation*}
	\begin{split}
		\xymatrix@C=2.5pc@R=1pc{
			0\ar[r]&  \mathend(E) \ar[r]&
			{\LdiffE} \ar@<0.5ex>[r]&
			L\ar[r]\ar@<0.5ex>@{-->}[l]& 0,\\
			0\ar[r]&  \Ad(P) \ar[r]&
			{\At_L(P)} \ar@<0.5ex>[r]&
			L\ar[r]\ar@<0.5ex>@{-->}[l]& 0.
		}
	\end{split}
\end{equation*}

\vskip 5pt
\subsubsection{The Atiyah extensions} 

Atiyah's seminal work~\cite{Atiyah57} aimed to establish criteria for the existence of holomorphic connections on a holomorphic vector bundle $\E$ over a complex manifold $M$. A key element of this characterization is a cohomology class, now known as the \textbf{Atiyah class}, denoted $\alpha(\E)\in H^{1,1}\big(M,\End(\E)\big)$, which serves as an obstruction to the existence of such connections. To briefly review its construction, consider the vector bundle ${J}^1 \E$ of first jets of holomorphic sections of $\E$. This bundle fits into the following canonical short exact sequence of holomorphic vector bundles over $M$, often referred to as the first jet sequence:
\begin{equation*}
	\xymatrix@C=2pc@R=2pc{
		0\ar[r]& \mathscr{T}^\vee M\otimes \E  \ar[r]&
		{J}^1 \E \ar[r]&
		\E \ar[r]& 0,
	}
\end{equation*}
where $\mathscr{T}^\vee M$ represents the holomorphic cotangent bundle of $M$. The Atiyah class of $\E$ is then defined as the extension class representing this short exact sequence:
$$\alpha(\E)\in\Ext^1_M\big(\E,\mathscr{T}^\vee M \otimes \E\big)\cong H^1\big(M,\mathscr{T}^\vee M \otimes{\rm End} \,\E\big)\cong H^{1,1}\big(M,\End(\E)\big),$$
as detailed in~\cites{Atiyah57}. 

The notion of the Atiyah class is generalized from a holomorphic vector bundle to a \textbf{Lie triad} $(L,A,E)$, where $L$ is a Lie algebroid, $A$ is a Lie sub-algebroid of $L$, and $E$ is a vector bundle endowed with a flat $A$-connection~\cite{CMP16}. 
Central to this generalization is the construction of a jet-bundle-like object, denoted by ${J}^1_{L/A}(E)$,
which participates in a canonical short exact sequence of vector bundles equipped with flat $A$-connections:
\begin{equation*} 
	\xymatrix@C=2pc@R=2pc{ 0\ar[r]& (L/A)^\vee\otimes E \ar[r]& {J}^1_{L/A}(E) \ar[r]& E \ar[r]& 0. }
\end{equation*}
The extension class of this sequence defines the generalized Atiyah class, denoted by $\alpha_{_{L,A}}(E)$, which lives in the first extension group $\Ext^1_A\big(E,(L/A)^\vee\otimes E\big)$.
This group is isomorphic to the first Lie algebroid cohomology group $H^1\big(A, (L/A)^\vee\otimes \End (E)\big)$. Crucially, the Atiyah class $\alpha_{_{L,A}}(E)$ serves as an obstruction to the existence of Lie algebroid connections of $L$ on $E$ that are compatible with the given flat $A$-connection.

While ${J}^1_{L/A}(E)$ serves as a generalization of ${J}^1\E$, its structure is considerably complex. A more algebraically tractable alternative for understanding the Atiyah class $\alpha(\E)$ of a holomorphic vector bundle $\E\to M$ involves considering a dual version of the holomorphic first jet sequence:

\begin{equation}\label{seq:AtiyahEdualversion}
	\xymatrix@C=2pc@R=2pc{
		0 \ar[r]&\End(\E)\ar[r]&\mathfrak{D}(\E)\ar[r]
		&\mathscr{T}M\ar[r]&0,
	}
\end{equation}

The extension class of this sequence is precisely the Atiyah class $\alpha(\E)$, see~\cite{Chan-Suen}. Indeed, the Atiyah bundle $\mathfrak{D}(\E)$, comprised of holomorphic covariant differential operators on $\E$, is a holomorphic Lie algebroid (as detailed in Section~\ref{SubSec:holomAtiyahalgebroid}). This contrasts with ${J}^1 \E$, which lacks comparable geometric or algebraic structures. Therefore, motivated by the relative simplicity of $\mathfrak{D}(\E)$ compared to ${J}^1\E$, we are led to the following:

\textbf{Objective 2} -- \textit{Instead of relying on $ {J}^1_{L/A}(E)$, we aim to construct Atiyah-type objects that encode compatible Lie algebroid connections.}

The central construction in our approach is a vector bundle, denoted $\mathfrak{D}_{L/A}(E)$, which we term the \textbf{Atiyah extension} of the Lie triad $(L,A,E)$. This extension allows us to define an analogous sequence to~\eqref{seq:AtiyahEdualversion} for the Lie triad $(L,A,E)$, specifically:
\begin{equation}\label{seq:AtiyahDBE0}
	\xymatrix@C=2pc@R=2pc{
		0 \ar[r]&\End(E)\ar[r]^{ }&\mathfrak{D}_{L/A}(E)\ar[r]^{ }
		&L/A\ar[r]&0.
	}
\end{equation}
(See Theorem~\ref{thm:DBE}.) In this context, $B:=L/A$ (equipped with an $A$-representation) replaces the role of $T^{1,0}X$, the underlying smooth vector bundle of $\mathscr{T}X$ with a canonical $T^{0,1}X$-representation, as found in Sequence~\eqref{seq:AtiyahEdualversion}. Critically, the extension class of Sequence~\eqref{seq:AtiyahDBE0} represents the Atiyah class $\alpha$.

We propose that constructing the Atiyah class via this method offers advantages over approaches based on the jet object ${J}^1_{L/A}(E)$. Specifically, $\mathfrak{D}_{L/A}(E)$ can be realized as a quotient of a natural Lie algebroid $\mathfrak{D}_{L}(E)$ by its Lie sub-algebroid $A$ (as detailed in Section~\ref{subsec:relativeAtiyah}). In contrast, ${J}^1_{L/A}(E)$ lacks a readily apparent geometric or algebraic structure. 
In fact, we present three equivalent characterizations of vector-bundle extensions equipped with flat $A$-connections (see Theorems~\ref{thm:DBE},~\ref{thm:DBEI}, and~\ref{thm:BdiffI=AtExtE}):
\begin{equation*}
	\xymatrix@C=2pc@R=1pc{
		0 \ar[r]&\End(E)\ar[r]
		&\mathfrak{D}_{L/A}(E)\ar[r]&L/A\ar[r]&0,\\
		0\ar[r]&\mathend(E)\ar[r]
		& \mathfrak{D}_{i_{_{L/A}}}(E)\ar[r]&L/A\ar[r]&0,\\
		0\ar[r]&\mathend(E)\ar[r]
		& \End(E) \nabplus L/A \ar[r]&L/A\ar[r]&0,
	}
\end{equation*}
where the third sequence provides a direct derivation of the Atiyah cocycle.
Based on these, we  establish a bijection between compatible Lie algebroid connections and compatible splittings of these Atiyah-type sequences,
and demonstrate that the Atiyah class serves as an obstruction to the compatibility conditions, (see Proposition~\ref{prop:AtiyahCocyle=0} and Theorem~\ref{thm:AtClass} ).

Finally, we assemble  the functorial relations among the objects  $\LdiffE$, $\fD_A(E)$ and $\BdiffE$ in terms of a centered hexagon~\eqref{Dia:Hexagon} in Appendix~\hyperlink{App:A}{A}:
\begin{equation*}
	\begin{split}
		\xymatrix@C=0.5pc@R=4pc{
			&L\ar[rr]&&B&\\
			A\ar[ru]\ar[rr]\ar@<-0.5ex>[dr]
			&&\fD_{L}(E)\ar[lu]
			\ar[ru]\ar[rr]
			&&\BdiffE\ar[lu]\\
			&\fD_A(E)\ar[ru]\ar@<-0.5ex>[lu]
			\ar@<-0.5ex>[rr]
			&&\End(E).\ar[lu]\ar[ru]
			\ar@<-0.5ex>[ll]&
		}
	\end{split}
\end{equation*}

 
 
 \vskip 5pt
\subsubsection{Matched connections, holomorphic and equivariant Atiyah algebroids}

As examples of these constructions, we 
study in Section~\ref{Sec:matchedpairHolomAtiyahalgebroid} the Atiyah algebroids and Atiyah extensions 
of a distinctive type of Lie triads $(L, A, E)$ where 
the quotient bundle $B = L/A$ is also a Lie subalgebroid of $L$, leading to $L$ being a matched-sum Lie algebroid of the matched pair $(A,B)$, denoted as $L \cong A \bowtie B$ (as detailed in Subsection~\ref{subsec:matchedPair}).

For instance, Atiyah's original construction, intended for finding holomorphic connections, was developed within the holomorphic category. In this context, the holomorphic tangent bundle $\mathscr{T} M$ and the Atiyah bundle $\mathfrak{D}(\E)$ are holomorphic Lie algebroids. According to~\cite{LSX08}, holomorphic Lie algebroids over a complex manifold $M$ are equivalent to matched pairs that result in a matched-sum  Lie algebroid $T^{0,1}M \bowtie L$, where~$L$ represents a smooth complex Lie algebroid.

The other case of matched pairs that we are interested in arises from the $\g$-equivariant setting: given a Lie algebroid $L$ equipped with an action of a Lie algebra $\g$ (Definition~\ref{dfn:g-algebroid}),
one obtains a matched-sum  Lie algebroid $(\g\ltimes M)\bowtie L$. In this scenario, we can study the existence of equivariant connections on an equivariant vector bundle.


More generally, when dealing with a matched-type Lie triad $(L=A\bowtie B,A,E)$, one can consider Lie algebroid connections that are compatible with the matched structure. We term such connections \textbf{matched connections} (Definition~\ref{Def:matchedconnection}). 
We wish to enrich the notion of Atiyah algebroids so that they encode matched Lie algebroid connections on both vector bundles and principal bundles. Consequently, we have:

\textbf{Objective 3} -- \textit{To develop  an enriched Atiyah algebroid framework capable of capturing the properties of  matched connections. }
  
The Atiyah extension $\BdiffE$ associated with the matched-type Lie triad $(L=A\bowtie B, A, E)$ provides the desired Lie algebroid structure. Specifically, $A$ and $\BdiffE$ combine to form the Lie algebroid isomorphism $A\bowtie \BdiffE \cong \mathfrak{D}_{A\bowtie B}(E)$. Moreover, this Atiyah extension fits into a commutative diagram of Lie algebroid sequences, as presented in Theorem~\ref{thm:matchedAtExt}:
 
 \begin{equation*}
 	\xymatrix@C=2pc@R=2pc{
 		0 \ar[r]&\End(E)\ar[r] \ar@{=}[d]
 		&\mathfrak{D}_{B}(E)\ar[r]\ar@{^{(}->}[d]&B\ar[r]\ar@{^{(}->}[d]&0\,\\
 		0 \ar[r]&\End(E)\ar[r]
 		&\mathfrak{D}_{A\bowtie B}(E)\ar[r]&A\bowtie B\ar[r]&0.\\
 	}
 \end{equation*}

This general template now nicely incorporates the connection theory for holomorphic or equivariant Lie algebroids by means of holomorphic or equivariant Atiyah algebroids.

\vskip 10pt

\subsection{Conventions and notations.}

Throughout the paper, $\K$ stands for the base field $\R$ or $\C$. All smooth manifolds are assumed to be Hausdorff, paracompact and smooth without boundary.
Some commonly used
symbols are listed below:

\begin{enumerate}
	\item $(L, A)$: a Lie pair;   $L$ is a Lie algebroid and $A$ is a Lie sub-algebroid in $L$.
	\item $(L, A, E)$: a Lie triad; $(L, A)$ is a Lie pair and $E$ a vector bundle equipped with a flat $A$-connection.
	\item $\diffE$: the Atiyah algebroid  of a  vector bundle $E$.
	\item $\LdiffE$: the $L$-Atiyah algebroid  of a  vector bundle $E$.
	\item ${\End(E)\nabplus L}$:   the virtually split $L$-Atiyah algebroid of $E$ with respect to an $L$-connection $\nabla$.
	\item $\BdiffE$: the $B$-Atiyah extension of $(L, A, E)$.
	\item $\BdiffI$: the embedded $B$-Atiyah extension of $(L, A, E)$.
	\item $\AtExtE$: the virtually split $B$-Atiyah extension associated with $(L,A,E)$ and $\nabla$.
\end{enumerate}

\vskip 25pt

\section{Preliminaries}\label{Sec:Pre}

\vskip 10pt

We give a short account of some basic knowledge of Lie algebroid, Lie pair and Lie triad.

\vskip 10pt
\subsection{Lie algebroids,   connections and representations}\label{Subsec: LieAlgoid}
\begin{Def}
	A \textbf{Lie $\K$-algebroid} $(L,[~\cdot,\cdot~]_{_L},\rho_{_L})$ over a smooth manifold $M$ consists of 
	\begin{enumerate}[(i)]
		\item a $\K$-vector bundle $L$ over $M$;
		\item a Lie bracket $[~\cdot,\cdot~]_{_L}$ on the $\K$-vector space of sections $\Gamma(L,\K)$;
		\item a vector bundle map $\rho_{_L}\colon L\to TM\otimes_{\R}\K$, called the anchor, 
	\end{enumerate}
	such that $\rho_{_L}$ induces a Lie algebra homomorphism from $\Gamma(L,\K)$ to $ \bsection{TM}\otimes_{\R}{\K}$ satisfying the Leibniz rule
	$$[l_1,fl_2]_{_L}=\big(\rho_{_L}(l_1)f\big)l_2+f[l_1,l_2]_{_L}, 
	\qquad \forall \, f \in C^{\infty} (M,{\K}),\, l_1,l_2 \in \Gamma(L,\K). $$
\end{Def}

\begin{Convention}
	We often omit to write $\K$ for simplicity if no confusion occurs, and we assume the geometric objects all have a fixed base manifold $M$ unless otherwise mentioned.
\end{Convention}

\begin{Def}
	A base-preserving vector bundle map $\Psi:L\to L'$ between Lie algebroids on $M$
	is called a \textbf{Lie algebroid morphism} if $\rho_{_L}=\rho_{_{L'}}\circ\Psi$ and $\Psi([~\cdot,\cdot~]_{_L})=[\Psi(\cdot),\Psi(\cdot)]_{_{L'}}$.
\end{Def}

\begin{Def}\label{def:L-Conn}
A \textbf{connection} of a Lie algebroid $L$ on a vector bundle $E$, both over $M$,  is a $\K$-bilinear map:
\[
\nabla\colon\bsection{L}\times\bsection{E}\longrightarrow\bsection{E},\quad (l,e)\longmapsto \nabla_{_l} e,
\] 
satisfying $\nabla_{_{fl}}e = f\nabla_{_l} e$ and $\nabla_{_l}(fe) = f\nabla_{_l} e + \big(\rho_{_L}(l)f\big)e$
for all $l\in\bsection{L}$, $e\in\bsection{E}$, and $f\in C^{\infty} (M)$.\\
\end{Def}

\begin{Def}\label{def:L-Curv}
The \textbf{curvature} of $\nabla$ is the $C^\infty(M)$-linear map $R^{\nabla}\colon \Lambda^2 \big(\Gamma(L)\big)\to \Gamma\big(\mathend(E)\big)$ defined by
\[R^{\nabla}(l_1,l_2)=\nabla_{l_1}\nabla_{l_2}-\nabla_{l_2}\nabla_{l_1}-\nabla_{[l_1,l_2]},\qquad \forall\, l_1,l_2\in\bsection{L}. \]
\end{Def}


View $\nabla$ as a linear map $\Gamma(E)\to\bsection{ L^\vee \otimes E}$, and set  $\Omega^k_L(E):=\bsection{ \wedge^k L^\vee \otimes E}$, 
then we can extend $\nabla$ to be a degree $(+1)$ derivation
$${d^\nabla_L}\colon\Omega^k_L(E)\longto \Omega^{k+1}_L(E),$$
using the Leibniz rule: ${d^\nabla_L}(\beta\otimes e)=d_L\beta\otimes e+(-1)^n \beta\wedge \nabla e$, for all $\beta\in \bsection{\Lambda^n  L^\vee }, e\in\bsection{E}$.
More explicitly, for $\omega\in\Omega^k_L(E)=\Gamma\big(\Hom(\wedge^k L, E)\big)$ and $l_i\in\bsection{L}$, we have that:
\[
(d^\nabla_L \omega)(l_0,l_1,\ldots,l_k) = \sum_{i=0}^k (-1)^i\nabla_{l_i}\big(\omega(l_0,\ldots,\hat{l_i},\ldots,l_k)\big)+\sum_{0\leqslant i<j\leqslant k}(-1)^{i+j}\omega\big([l_i,l_j],l_0,\ldots,\hat{l_i},\ldots,\hat{l_i},\ldots,l_k\big).
\]


\begin{Def}\label{dfn:L-repr}
A vector bundle $E$ endowed with a \textbf{flat} $L$-connection $\nabla$, i.e., 
$R^{\nabla}=0$, is called an \mbox{\textbf{$L$-representation}}.
Equivalently, the derivation ${d^\nabla_L}\colon\Omega^k_L(E)\to \Omega^{k+1}_L(E)$ satisfies $(d^\nabla_L)^2=0$, and will be called the \textbf{Chevalley-Eilenberg differential}.
The corresponding cohomology groups are denoted by $H^\bullet_{\mathrm{CE}} (L,E)$.	 
\end{Def}

\begin{Def}\label{dfn:morphConn}
Let $\nabla$ be an $L$-connection on $E$, and $\nabla'$ an $L'$-connection on $E'$.
A \textbf{morphism between two connections $\nabla,\,\nabla'$} consists of two base-preserving bundle maps:
a Lie algebroid morphism $\Psi:L\to L'$ and a vector bundle map $\Phi:E\to E'$
satisfying $\Phi(\nabla_{_l} e)=\nabla'_{_{\Psi(l)}} \big(\Phi(e)\big)$  for all $l\in \Gamma(L),\,e\in \Gamma(E)$.
\end{Def}

\vskip 10pt
\subsection{Lie pairs}\label{subsubSec:Liepair}
\begin{Def}
	A \textbf{Lie algebroid pair} $(L,A)$, simply called a \textbf{Lie pair}, consists of a Lie algebroid $L$ and a Lie subalgebroid $A\overset{i_{_A}}{\lhook\joinrel\longrightarrow} L$ 
	over the same base manifold $M$.
\end{Def}

\begin{Ex}\label{ex:antiholom}
	Let $(M,J)$ be an almost complex manifold.
	By the Newlander--Nirenberg theorem, the almost complex structure $J$ is integrable if and only if the $(0,1)$-tangent bundle $T^{0,1}M\subseteq T_\C M=TM\otimes \C$ is involutive,
	i.e., $T^{0,1}M$ is a Lie subalgebroid of $T_\C M$ and $(T_\C M,T^{0,1}M)$ forms a Lie pair.
\end{Ex}

\begin{Ex}\label{ex:distribution}
	Let $M$ be a real manifold. A Lie subalgebroid of the tangent Lie algebroid $TM$ is an involutive distribution $F\subseteq TM$ of constant rank.
	We then have a Lie pair $(TM,F)$.
\end{Ex}

\begin{Ex}\label{ex:generalizedCplx}
	A generalized complex structure $\mathcal{J}$ on $M$ is defined by Hitchin~\cite{Hitchin03} to be an endomorphism $\mathcal{J}\colon TM\oplus T^*M\to TM\oplus T^*M$ satisfying $\mathcal{J}^2=-\id$ and that the $(\pm 1)$-eigenbundles $L_\pm$ are involutive under the Courant bracket.
	If $\mathcal{J}$ is regular in the sense that $\rho(L_+)\cap \rho(L_-)$ is of constant rank for the complexified projection $\rho\colon T_\C M \oplus T^*_\C M \to T_\C M$, then the involutive distribution $\rho(L_-)\subseteq T_\C M$ is also of constant rank.
	Such a $\mathcal{J}$ gives a Lie pair $\big(T_\C M,\rho(L_-)\big)$.
\end{Ex}

The quotient bundle $B:=L/A$
and the projection ${\rm pr}_{_B}\colon L\rightarrow B$ gives a short exact sequence of vector bundles,
which we call the \textbf{quotient sequence of $(L,A)$}:
\begin{equation}\label{eq:LAsequence}
	\xymatrix@C=2pc@R=2pc{
		0\ar[r]&  A \ar[r]^{i_{_A}}& L \ar[r]^{{\rm pr}_{_B}}&
		B\ar[r]& 0.
	}
\end{equation}
There is a natural flat connection of $A$ on $B$ called the \textbf{(flat) Bott connection} or \textbf{Bott representation}:
\begin{equation}\label{eq:Bott}
	D: \Gamma(A)\times\Gamma(B)\longrightarrow \Gamma(B), \quad (a,b)\longmapsto D_a \,b:={\rm pr}_{_B}\big([a,\tilde{b}]_{_L}\big), \qquad \forall\, a\in \Gamma(A),\,b\in \Gamma(B),
\end{equation}
where $a\in \Gamma(A)$ is identified with its embedded image $\iA(a)\in \Gamma(L)$ and $\tilde{b}$ is any element in ${\rm pr}_{_B}^{-1}(b) \subset \Gamma(L)$.
The flatness of $D$ comes from the Jacobi identity for $[~\cdot,\cdot~]_{_L}$. 


\begin{Def}
	A \textbf{Lie pair morphism} between Lie pairs $(L,A)$ and $(L',A')$
	consists of two Lie algebroid morphisms 
	$\Psi_L:L\to L'$ and $\Psi_A:A\to A'$ such that
	$\Psi_L\circ i_{_A}=i_{_{A'}} \circ\Psi_A$, see the left commutative square in diagram~\eqref{eq:morphismOfLAsequence} below. 
\end{Def}

The following is easy to check.

\begin{prop}\label{prop:morphismOfBottConn}
	A Lie pair morphism $(\Psi_L,\,\Psi_A)$ between two Lie pairs $(L,A)$ and $(L',A')$ induces a natural morphism between their quotient sequences:
	\begin{equation}\label{eq:morphismOfLAsequence}
		\begin{split}
			\xymatrix@C=2.5pc@R=2pc{
			0\ar[r]&  A \ar[r]^{i_{_A}} \ar[d]_{\Psi_A}
			\ar@{}[dr]|(0.5){\scalebox{1.2}{$\circlearrowleft$}}
			& L \ar[r]^{{\rm pr}_{_B}} \ar[d]^{\Psi_L}
			\ar@{}[dr]|(0.5){\scalebox{1.2}{$\circlearrowleft$}}
			& B \ar[r] \ar[d]^{\Psi_B} & 0\,\\
			0\ar[r]&  A' \ar[r]^{i_{_{A'}}}& L' \ar[r]^{{\rm pr}_{_{B'}}}&
			B'\ar[r]& 0,
		}
		\end{split}
	\end{equation}
	where the induced vector-bundle map $\Psi_B$ is evaluated on sections as
	\begin{align}\label{eq:PsiB}
		\Psi_B\colon \Gamma(B)\longrightarrow \Gamma(B'),\quad \prB(l)\longmapsto \mathrm{pr}_{_{B'}}\big(\Psi_L(l)\big),\qquad\forall\, l\in \Gamma(L),
	\end{align}
	and is a morphism between the Bott connection $D$ of $A$ on $B$ and the Bott connection $D'$ of $A'$ on $B'$. 
\end{prop}
%

There are always non-canonical vector-bundle decompositions $L\cong A\oplus B$
by specifying right splittings of quotient sequence~\eqref{eq:LAsequence}
via vector-bundle embeddings $\iB\colon B\to L$ such that ${\rm pr}_{_B} \circ i_{_B}=\id_B$,
which also correspond to left splittings via vector-bundle projections
$\prA\colon L\to A$ such that $i_{_A}\circ {\rm pr}_{_A}=\id_A$.
The mutual determination between $\iB$ and $\prA$ is that
$i_{_A}\circ {\rm pr}_{_A} +i_{_B}\circ {\rm pr}_{_B}=\id_L$.
Consequently, quotient sequence~\eqref{eq:LAsequence} is supplemented with a reversed short exact sequence of vector bundles:
\begin{equation}\label{Seq:ABdecomposition}
	\xymatrix@C=2.5pc@R=2pc{
		0\ar[r]&  A \ar@<0.5ex>[r]^{i_{_A}} &
		L \ar@<0.5ex>[r]^{{\rm pr}_{_B}} \ar@<0.5ex>@{-->}[l]^{\,\,\,{\rm pr}_{_A}}&
		B\ar[r] \ar@<0.5ex>@{-->}[l]^{{i}_{_B}} & 0.
	}
\end{equation}

If we have two different right splittings $\iB,\,i'_{_B}:B\to L$,
their images must differ by $i'_{_B}(b)-\iB(b)$ in $\iA\big(\Gamma(A)\big)\subseteq \Gamma(L),\,\forall\, b\in \Gamma(B)$.
Hence, there is a vector bundle map measuring the difference:
\begin{align}\label{eq:iB-i'B}
	I\colon B\longto A,\quad\mbox{s.t.}\quad i'_{_B}(b)-\iB(b) = \iA \big(I (b)\big), \qquad \forall\, b\in \Gamma(B).
\end{align}

Now, fixing such a decomposition of $L\cong A\oplus B$ via $\iB$ and $\prA$,
we get a (non-canonical) $B$-operation $\eth$ on $A$, although $B$ is not an honest Lie algebroid:
\begin{align}\label{eq:eth}
	\eth:\Gamma(B)\otimes\Gamma(A)\longrightarrow \Gamma(A),\quad
	(b,a)\longmapsto \eth_{_b} a:={\rm pr}_{_A}\big[\iB(b),a\big]_{_L}, \qquad \forall \, a\in \Gamma(A),\,b\in \Gamma(B).
\end{align}
After identifying $b\in \Gamma(B)$ with its embedded image $\iB(b)\in \Gamma(L)$, equations~\eqref{eq:Bott}\,\eqref{eq:eth} give:
\begin{align}\label{eq:Lbracket}
	[a,b]_{_L}=D_a b-\eth_{_b} a, \qquad \forall \, a\in \Gamma(A),\,b\in \Gamma(B).
\end{align}

\vskip 10pt
\subsection{Lie triads and  Atiyah classes}
Here, we introduce a notion of Lie triad.
The properties of its Atiyah class is predominantly sourced from~\cite{CMP16}. 
\begin{Def}\label{defn:triad}
	A \textbf{Lie triad} is a triple $(L,A,E)$ such that $(L,A)$ forms a Lie pair and
	$E$ is a vector bundle on $M$ equipped with a flat $A$-connection $\barnab$.
\end{Def}

\begin{rmk}\label{rmk:LAP}
	The vector bundle $E$ can also be replaced by a principal bundle $P$ equipped with a flat $A$-connection in the sense of later Definition~\ref{def:Pcurv}, resulting in the principal bundle version of Lie triad~$(L,A,P)$.
\end{rmk}

\begin{Ex}[Continuing Example~\ref{ex:antiholom}] \label{Ex:triadOfHoloBundle}
	Notice that a holomorphic vector bundle $\E$ on a complex manifold $M$
	is equivalent to a Dolbeault operator $\parbar$ satisfying $\parbar^{\,2}=0$ 
	on its underlying smooth $\C$-vector bundle $E$.
	That is, $\parbar$ is a flat connection of $T^{0,1} M$ on $E$.
	It gives a complex-valued Lie triad $(T_\C M ,\,T^{0,1} M,\, E)$.
\end{Ex}

\begin{Ex} [Continuing Example~\ref{ex:distribution}]\label{ex:foliatedVB}
	By the Frobenius theorem, let $\mathcal{F}$ be the regular foliation that integrates the involutive sub-bundle $F\subseteq TM$ such that $T\mathcal{F}=F$.
	A vector bundle $E$ over $M$ is $\mathcal{F}$-foliated if $E$ is equipped with a flat $F$-connection $\barnab$, as introduced by Molino~\cite{Mo71} and Kamber--Tondeur~\cite{KT74}.
	Thus, we get a Lie triad $(TM,\,T\mathcal{F},\,E)$.
\end{Ex}

\begin{Ex}[Continuing Example~\ref{ex:generalizedCplx}]\label{ex:generalizedHolVB}
	A generalized holomorphic vector bundle on a generalized complex manifold $(M,\mathcal{J})$ is defined by Gualtieri~\cite{Gualtieri} to be a complex-valued bundle $E$ with a flat $L_-$-connection~$\parbar$.
	If $\mathcal{J}$ is regular, then $\rho(L_-)\subseteq T_\C M$ is a Lie subalgebroid and inherits the flat connection $\parbar$.
	Hence, there is a Lie triad $\big(T_\C M, \rho(L_-), E\big)$.
\end{Ex}


\begin{defn}\label{def:extendNabla}
	An \textbf{extending $L$-connection} of the Lie triad $(L,A,E)$
	is an $L$-connection $\nabla$ on $E$ such that $\nabla_a=\barnab_a$ for all $a\in \Gamma(A)$.
\end{defn}

\begin{rmk}
	The existence of an extending connection $\nabla$ is guaranteed by using an argument of partition of unity, 
	or simply by defining $\nabla_{a+\iB(b)}:=\barnab_{a}+\nabla'_{{\iB(b)}},\,\forall a\in \Gamma(A),b\in \Gamma(B)$
	where $\nabla'$ is any $L$-connection on $E$ and $\iB:B\to L$ a right splitting.
\end{rmk}

\begin{Ex}
	For the Lie triad $(L=T_\C M ,\,A=T^{0,1} M,\,E)$ in Example~\ref{Ex:triadOfHoloBundle} with the usual Dolbeault operator $\parbar$,
	a $T_\C M$-connection that extends $\parbar$ is a $(1,0)$-connection. 
	For example, the Chern connection arising from a Hermitian metric on $E$   is determined by requiring the compatibility with the Hermitian metric and extending $\parbar$.
\end{Ex}

Take any extending connection $\nabla$ and pick a right splitting $\iB\colon B\to L$,
we can restrict the curvature $R^{\nabla}$ on $A\otimes B$ to be the $C^\infty(M)$-linear map $R^{\nabla}_{A\otimes B}\triangleq R^{\nabla}|_{A\otimes B}\colon A\otimes B\longto \mathend(E)$ given by
\begin{equation}\label{eq:AtiyahCocycle}
	R^{\nabla}_{A\otimes B}(a\otimes b)=\nabla_a\nabla_{\iB(b)}-\nabla_{\iB(b)}\nabla_a-\nabla_{[a,\iB(b)]},\quad \forall a\in\bsection{A}, b\in\bsection{B},
\end{equation}
which in fact is independent from the choice of $\iB$ due to the flatness of the $A$-connection $\barnab$.


\begin{Def}\label{def:compatibleConn}
	An  extending connection $\nabla$   is said to be \textbf{$A$-compatible}  if
	it  
	satisfies
	\[\nabla_a\nabla_l-\nabla_l\nabla_a=\nabla_{[a,l]},\quad \forall a\in\bsection{A}, l\in\bsection{L}. \quad \mbox{Or equivalently, } R^{\nabla}_{A\otimes L}\triangleq R^{\nabla}|_{A\otimes L} \mbox{ is zero}.\]
\end{Def}

The vanishing of $R^{\nabla}_{A\otimes L}$ equates with that of $R^{\nabla}_{A\otimes B}$.
In~\cite{CMP16}, Sti\'enon, Xu, and one of the  authors proved the following properties.

\begin{thm}[{\cite[Thm.~2.5]{CMP16}}]\label{thm:CSX}
	Note that $R^{\nabla}_{A\otimes B} \in\Gamma\big({A^\vee}\otimes {B^\vee}\otimes \mathend(E)\big)=\Omega_A^1\big({B^\vee}\otimes \mathend(E)\big)$, the following hold:
	\begin{enumerate}[(a)]
		\item $R^{\nabla}_{A\otimes B}$ is a $1$-cocycle with respect to the Chevalley-Eilenberg differential for $\Omega_A^\bullet\big({B^\vee}\otimes \mathend(E)\big)$;
		\item\label{item:AtiyahIndepOfConn} $R^{\nabla}_{A\otimes B}$ defines a cohomology class $\alpha(L,A,E)\in H^1_{\mathrm{CE}}\big(A,\,B^\vee\otimes\mathend(E)\big)$ actually independent from the choice of an extending connection $\nabla$;
		\item\label{item:AtiyahVanish} the cohomology class $\alpha(L,A,E)$ vanishes if and only if there exists an $A$-compatible $L$-connection on $E$.
	\end{enumerate}
\end{thm}

\begin{defn}
	The $1$-cocycle $R^{\nabla}_{A\otimes B}$ is called the \textbf{Atiyah cocycle of $E$ relative to $(L,A,\nabla)$}.
	Its cohomology class $\alpha(L,A,E)\in H^1_{\mathrm{CE}}\big(A,\,B^\vee\otimes\mathend(E)\big)$ is called the \textbf{Atiyah class of $E$ relative to $(L,A)$}, which in our current terminology
	is also called the \textbf{Atiyah class of the Lie triad $(L,A,E)$}.
\end{defn}

\begin{Ex}[Continuing Example~\ref{Ex:triadOfHoloBundle}]
	It was Atiyah's original idea~\cite{Atiyah57} that the existence of holomorphic connections on the holomorphic vector bundle $\E$  is precisely obstructed by the Atiyah class $\alpha(T_\C M,T^{0,1}M,E)$.
\end{Ex}

\begin{Ex}[Continuing Example~\ref{ex:foliatedVB}]
	The $\mathcal{F}$-foliated vector bundle $E$ admits certain compatible connections, called projectable~\cite{Mo71} or basic~\cite{KT74}, if and only if the Atiyah class $\alpha(TM,\,T\mathcal{F},\,E)$ vanishes.
\end{Ex}

\begin{Ex}[Continuing Example~\ref{ex:generalizedHolVB}]
	The generalized holomorphic vector bundle $E$ on a regular generalized complex manifold $(M,\mathcal{J})$ admits a generalized holomorphic connection if and only if the Atiyah class $\alpha\big(T_\C M,\, \rho(L_-),\, E\big)$ vanishes, see~\cite{LangLiu23}.
\end{Ex}

\vskip 25pt

\section{Generalized Atiyah algebroids and Lie algebroid connections} \label{Sec:generalAtiyahalgebroid}

\vskip 10pt

In this part, we begin with a quick review of the classic Atiyah algebroids that stem from the usual connections on principal bundles and vector bundles,
then extend them to the setting of general Lie algebroid connections respectively on vector bundles and principal bundles,
and also obtain a structural description of the Atiyah algebroid for a flat Lie algebroid connection.

\vskip 10pt
 
\subsection{Classic Atiyah algebroids}\label{subsubSec:expansionLiealgebroid}
Let us recall the notions of the Atiyah algebroid $ (TP)  /G$ of a principal $G$-bundle $P$
and the Atiyah algebroid $\diffE$ of a  vector bundle $E$ (see~\cite{Mac05}).

For a principal $G$-bundle $P$, the tangent bundle $TP$ is naturally endowed with a free right $G$-action.
The Lie algebroid structure of $TP$ over $P$ then descends to form a Lie algebroid  $(TP)/G=:\At(P)$ over $M$, called  the  \textbf{Atiyah algebroid of $P$}, which fits in an exact sequence of Lie algebroids
\begin{equation*}\label{seq:At(P)}
	\xymatrix@C=2pc@R=2pc{
		0\ar[r]&  P\times_{_G} {\g} \ar[r]& (TP) /G \ar[r]^{\quad\widetilde{T\pi}}&
		TM \ar[r]& 0,
	}
\end{equation*}
called the \textbf{Atiyah sequence of $P$}.
Here $P\times_{_G} {\g}=:\Ad(P)$ is the associated vector bundle of $P$ for the adjoint $G$-action on Lie algebra $\g$ and has a Lie algebra bundle structure,
and $\widetilde{T\pi}$ is induced from the tangent map of the projection $\pi\colon P\to M$. Atiyah sequence of $P$ is also often expressed as
\begin{equation*}
	\xymatrix@C=2pc@R=2pc{
		0\ar[r]&  \Ad(P) \ar[r]&
		{\At(P)} \ar[r]&
		TM \ar[r]& 0.
			}
\end{equation*}

Given a vector bundle $E$, one has the frame bundle $P=\mathrm{Fr(E)}$ with structure group $G=GL\big(\rank(E),\K\big)$.
This $(TP)/G$ is the \textbf{Atiyah algebroid} $\diffE$ whose fiber at each $m\in M$ consists of
linear covariant derivative operators $\delta_m:\Gamma(E)\to E_m$ 
equipped with a tangent vector symbol $X_m\in T_m M $ such that
$$\delta_m(fe)=f(m)\delta_m(e) + X_m(f)e(m),\qquad \forall\,f\in C^\infty(M),\,e \in\Gamma(E).$$
The section space $\Gamma\big(\diffE\big)$ 
is then the set of global linear covariant differential operators
${\delta\colon \Gamma(E)\to \Gamma(E)}$ with symbols $X\in \Gamma(TM) $ which satisfy
$$\delta(fe)=f\delta(e) + X(f)e,\qquad \forall\, f\in C^\infty(M),\, e \in\Gamma(E).$$
As a Lie algebroid, the anchor of $\diffE$  is the symbol map $\sigma:\delta\mapsto X$
and   the Lie bracket on $\Gamma\big(\diffE\big)$ is just the usual commutator of two linear operators on $E$. Moreover, $\diffE$ fits in the exact sequence of Lie algebroids
\begin{equation*}\label{Seq:diffE}
	\xymatrix@C=2pc@R=2pc{
		0\ar[r]&  \mathend(E) \ar[r]& {\diffE} \ar[r]^{\sigma}&
		TM \ar[r]& 0,
	}
\end{equation*}
called the \textbf{Atiyah sequence of $E$}.
Here $\End(E)$ has zero anchor and its Lie bracket on $\Gamma\big(\End(E)\big)$ is the usual commutator of linear endomorphisms on $E$.

\vskip 10pt
\subsection{$L$-Atiyah algebroids of vector bundles} \label{subsec:DLE}

\vskip 5pt
\subsubsection{Construction of an  $L$-Atiyah algebroid}

For a Lie algebroid $(L,[~\cdot,\cdot~]_{_L},\rho_{_L})$ over $M$,
its anchor map  $\rho_{_L}\colon L\to TM $ together with
the anchor map of Atiyah algebroid $\sigma\colon \diffE\to TM  $
gives a pullback construction: 
\begin{equation}\label{eq:DLE}
	\LdiffE \triangleq \Big\{(\delta_m,l_m) \in \bigcup_{m\in M}
	\big(\diffE_m\times L_m\big) \bigm|\sigma(\delta_m)=\rho_{_L}(l_m) \in T_m M  \Big\}.
\end{equation}
Due to the surjectivity of $\sigma$, we see that $\LdiffE$ is a vector bundle over $M$
and it projects onto $L$ with the kernel $\mathend(E)$ same as $\ker \sigma$.
Diagrammatically, $\LdiffE$ fits in a morphism between exact sequences of vector bundles:
\begin{equation}\label{Dia:LdiffE-pullback}
	\begin{split}
		\xymatrix@C=2pc@R=2pc{
			0\ar[r]& \mathend(E)
			\ar[r]^{{\iota}}
			\ar@{=}[d]& 
			\LdiffE
			\ar[r]^{\sigma_{_L}}
			\ar[d]_{f_{_L}}
			\ar@{}[dr]|(0.25){\text{\pigpenfont J}}&
			L
			\ar[r]
			\ar[d]^{\rho_{_L}}&
			0\,
			\\
			0\ar[r]& \mathend(E)
			\ar[r]
			& 
			\diffE
			\ar[r]^{\sigma}
			&
			TM 
			\ar[r]
			&
			0,}
	\end{split}
\end{equation}
where the right square is the pullback diagram for $\LdiffE$ and ${\iota},\,\sigma_{_L},\,f_{_L}$ are explicitly given as
\begin{equation*}\label{eq:symbol-DLE}
	{\iota}(\varphi_m)=(\varphi_m,0),\qquad f_{_L}(\delta_m,l_m)=\delta_m,\qquad \sigma_{_L}(\delta_m,l_m)=l_m,
\end{equation*}
for all $m\in M$, $\varphi_m\in \End(E_m)$, and $(\delta_m,l_m)\in\LdiffE_m$.

By construction~\eqref{eq:DLE}, the section space of $\LdiffE$ is
\begin{align}\label{eq:ΓDLE}
	\Gamma\big(\LdiffE\big)=\big\{(\delta,l)\in \Gamma\big(\diffE\big)\times \Gamma(L) \bigm | \sigma(\delta)=\rho_{_L}(l)\in \Gamma(TM )\big\}.
\end{align}
Such a section $(\delta,l)$ has the underlying covariant differential operator $\delta$ whose symbol is $\sigma(\delta)=\rho_{_L}(l)$.
Hence, we call $(\delta,l)\in\Gamma\big(\LdiffE\big)$ a \textbf{covariant $L$-differential operator on $E$} and call $l=\sigL(\delta,l)$ the \textbf{$L$-symbol} of $(\delta,l)$.
In this perspective, $\LdiffE$ becomes the bundle of covariant $L$-differential operators on $E$.

The vector bundle $\LdiffE$ is naturally a Lie algebroid with its Lie bracket and anchor defined for $(\delta,l),\,({\delta}_i,l_i)$ in $\Gamma\big(\LdiffE\big)$ to be
\begin{equation}\label{eq:LdiffE}
\begin{aligned} 
	\big[(\delta_1,l_1),({\delta}_2,l_2)\big]_{\LdiffE}&\triangleq\big([\delta_1,\delta_2]_{\diffE},\,[l_1,l_2]_{_L}\big),\\
	\rho_{_{\LdiffE}}(\delta,l)&\triangleq\sigma(\delta)=\rho_{_L}(l), \qquad \text{i.e., } \rho_{_{\LdiffE}}=\rho_{_L}\circ\sigma_{_L}=\sigma\circ f_{_L}.
\end{aligned}
\end{equation}
The sequence of vector bundles
\begin{equation}\label{seq:LdiffE}
	\xymatrix@C=2pc@R=2pc{
		0\ar[r]&  \mathend(E) \ar[r]^{{\iota}}& {\LdiffE} \ar[r]^{\quad \sigma_{_L}}&
		L\ar[r]& 0,
	}
\end{equation}
then becomes an exact sequence of Lie algebroids over $M$.

\begin{Def}\label{def:LdiffE}
	We call $\LdiffE$ the \textbf{$L$-Atiyah algebroid of $E$},
	and   Sequence~\eqref{seq:LdiffE} the \textbf{$L$-Atiyah sequence of $E$}.
\end{Def}

The functoriality for Atiyah algebroids $\LdiffE$ with respect to $L$ can be stated as follows.

\begin{prop}\label{prop:functoriality}
	Let $\Psi:L\to L'$ be a Lie algebroid morphism over $M$,
	then it naturally induces a Lie algebroid morphism $f_{_\Psi}$
	between Atiyah algebroids that fits in a morphism between Atiyah sequences:
	\begin{equation*}\label{diag:functoriality}
		\begin{split}
			\xymatrix@C=2pc@R=2pc{
				0 \ar[r]
				&\mathend(E) \ar[r]^{{\iota}}\ar@{=}[d]
				&\mathfrak{D}_L(E) \ar[r]^{\quad\sigma_{_L}}\ar[d]_{f_{_\Psi}}
				\ar@{}[dr]|(0.25){\text{\pigpenfont J}}
				&L \ar[r]\ar[d]^{\Psi}
				&0\,
				\\
				0 \ar[r]
				&\mathend(E) \ar[r]^{{\iota}}
				& \mathfrak{D}_{L'}(E) \ar[r]^{\quad\sigma_{_{L'}}}
				&L'\ar[r] 
				&0,
			}
		\end{split}
	\end{equation*}
	where the right square is a pullback of Lie algebroids on $M$, and 
	the morphism $f_{_\Psi}$  is
	\begin{align*}\label{eq:functoriality}
		f_{_\Psi}\colon \Gamma\big(\mathfrak{D}_L(E)\big)\longto \Gamma\big(\mathfrak{D}_{L'}(E)\big),
		\quad (\delta,l)\longmapsto \big(\delta,\Psi(l)\big).
	\end{align*} 
\end{prop} 

\vskip 5pt
\subsubsection{Splitting of the $L$-Atiyah algebroid by an $L$-connection}\label{subsubsec:splitAtiyah}
Just like the interpretation~\cite{Atiyah57} of a classic connection as a splitting of the classic Atiyah sequence, similar theorem holds for $L$-Atiyah sequences.

\begin{thm}\label{thm:conn=splitting}
	The following data are equivalent:
	\begin{enumerate}[(1)]
		\item an $L$-connection $\nabla$ on $E$ (Definition~\ref{def:L-Conn});
		\item a vector-bundle right splitting
		$s\colon L\to \LdiffE$ of $L$-Atiyah sequence~\eqref{seq:LdiffE},
		or equivalently, a vector-bundle left splitting $\theta\colon\LdiffE \to \mathend(E)$ of~\eqref{seq:LdiffE} with the mutual determination 
		${\iota}\circ\theta+s\circ\sigma_{_L}=id_{\LdiffE}$:
		\begin{equation*}\label{seq:LdiffE2}
			\xymatrix@C=2.5pc@R=2pc{
				0\ar[r]&  \mathend(E) \ar@<0.5ex>[r]^{\,\,{\iota}}&
				{\LdiffE} \ar@<0.5ex>[r]^{\quad \sigma_{_L}} \ar@<0.5ex>@{-->}[l]^{\,\,\theta}&
				L\ar[r]\ar@<0.5ex>@{-->}[l]^{\quad s}& 0;
			}
		\end{equation*}
		\item a vector-bundle decomposition $\LdiffE\cong \mathend(E)\oplus L$ such that $\iota(\mathend(E))=\mathend(E)\oplus 0$ and ${\sigL(0\oplus L)=L}$;
		\item a vector-bundle map $\bar{s} \colon L\to \diffE$ that makes the following lower triangle commute:
		\begin{equation*}
			\begin{split}
				\xymatrix@C=2.5pc@R=2pc{
					\LdiffE
					\ar@<0.5ex>[r]^{\sigL}
					\ar[d]_{f_{_L}}
					\ar@{}[dr]|(0.25){\text{\pigpenfont J}}|(0.7){\scalebox{1.2}{$\circlearrowleft$}}
					& L\, \ar[d]^{\rhoL}
					\ar@{-->}@<0.5ex>[l]^{s} \ar@{-->}[ld]|(0.5){\scalebox{1}{$\bar{s}$}}
					\\
					\diffE 
					\ar[r]^{\sigma}
					&TM.}
			\end{split}
		\end{equation*}
		(The correspondence between $s$ and $\bar{s}$ is that $\bar{s}=f_{_L}\circ s$.)
	\end{enumerate}
\end{thm}

\begin{proof}
	$(1)\Leftrightarrow(2)\colon$ Given an $L$-connection $\nabla$,
	we have a $C^\infty(M)$-linear map ${s\colon \Gamma(L)\to \Gamma\big(\LdiffE\big),\,l\mapsto (\nabla_l,l),}$
	which satisfies $\sigL\circ s=\id_L$, and hence splits $L$-Atiyah sequence~\eqref{seq:LdiffE} on the right.
	Conversely, given a vector-bundle right splitting $s\colon L\to \LdiffE$ satisfying $\sigL\circ s=\id_L$,
	we set a connection
	${\nabla}_l=s(l)\in \Gamma\big(\LdiffE\big)$
	whose $L$-symbol is exactly $l$.
	
	$(2)\colon$ The equivalence between left and right splittings is clear.
	
	$(2)\Leftrightarrow(3)\colon$ The vector-bundle splittings $s$ or $\theta$, together with their corresponding $L$-connection $\nabla$,
	give an explicit decomposition $\LdiffE\cong \mathend(E)\oplus L$:
	\begin{equation}\label{eq:LdiffE=EndE+L}
		\begin{aligned}
			\Gamma\big(\LdiffE\big) &\longrightarrow \Gamma\big(\mathend(E)\oplus L\big):
			&(\delta,l)&\longmapsto \theta(\delta,l)\oplus \sigL(\delta,l) =(\delta-\nabla_{_l})\oplus l,\\
			\Gamma\big(\mathend(E)\oplus L\big) &\longrightarrow \Gamma\big(\LdiffE\big):
			& \varphi \oplus l &\longmapsto {\iota}(\varphi)+s(l)=(\varphi+\nabla_{_l},l),
		\end{aligned} 
	\end{equation}
	where elements of $\Gamma\big(\mathend(E)\oplus L\big)$ are written as
	$\varphi \oplus l$ for $\varphi\in\Gamma\big(\mathend(E)\big),\,l\in\Gamma(L)$.
	The converse is evident.
	
	$(2)\Leftrightarrow(4)\colon$ This is due to the universal property of $\LdiffE$ as the pull-back of $L$ and $\diffE$.
\end{proof}

\begin{cor} \label{cor:conn=splitting}
	With notations as in Theorem~\ref{thm:conn=splitting}, we have the following conclusions:  
	\begin{enumerate}[(a)]
		\item $L$-connections $\nabla$ on $E$ always exist, and they form an affine space modelled on $\Gamma\big(\Hom\big(L,\End(E)\big)\big)$;
		\item an $L$-connection $\nabla$ is flat if and only if
		the corresponding splitting $s$ is a Lie algebroid morphism making $L$ a Lie subalgebroid of $\LdiffE$,
		equivalently, if and only if $\bar{s}$ is a Lie algebroid morphism.
	\end{enumerate}
\end{cor}
\begin{proof}
	Clause (a) follows from the fact that vector-bundle right splittings $s\colon L\to \LdiffE$ of~\eqref{seq:LdiffE} always exist
	and any two of such differ by an element in $\Gamma\big(\Hom(L,\End(E))\big)$.
	For (b), the map $\bar{s}$ or the splitting $s$ respects Lie algebroid structure 
	iff $[s(l_1),s(l_2)]-s([l_1,l_2])=\big([\nabla_{l_1},\nabla_{l_1}]-\nabla_{[l_1,l_2]},0\big)=\big(R^{^\nabla}(l_1,l_2),0\big)$
	vanishes for any $l_1,\,l_2\in \Gamma(L)$.
\end{proof}

With the above identification~\eqref{eq:LdiffE=EndE+L} of $\LdiffE$ and $ \mathend(E)\oplus L$, the following fact is immediate.

\begin{thm}\label{thm:EndE+L}
	Given an $L$-connection $\nabla$ on $E$, the identifications in~\eqref{eq:LdiffE=EndE+L} transfer Lie algebroid structure~\eqref{eq:LdiffE} of $\LdiffE$ to Lie algebroid structure of $\mathend(E)\oplus L$ as follows:
	\begin{enumerate}[i)]
		\item the anchor $\rho^\nabla\colon \mathend(E)\oplus L\to TM$ of $\varphi \oplus l\in \Gamma\big(\mathend(E)\oplus L\big)$ is
		\begin{equation*}\label{eq:rho_EndE+L}
			\rho^\nabla(\varphi \oplus l)\triangleq\rho_{_L}(l);
		\end{equation*}
		\item \label{item:[]_EndE+L}the Lie bracket $[~\cdot,\cdot~]^\nabla\colon \Gamma\big(\mathend(E)\oplus L\big)\times\Gamma\big(\mathend(E)\oplus L\big)\to\Gamma\big(\mathend(E)\oplus L\big)$ for $\varphi_i \oplus l_i\in \Gamma\big(\mathend(E)\oplus L\big)$ reads
		\begin{equation*}
			[\varphi_1 \oplus l_1,\varphi_2\oplus l_2]^\nabla
			\triangleq\Big([\varphi_1,\varphi_2]+\nabla^{\mathend(E)}_{l_1}\varphi_2-\nabla^{\mathend(E)}_{l_2}\varphi_1+R^{\nabla}(l_1,l_2)\Big) \oplus [l_1,l_2],
		\end{equation*}
		where $\nabla^{\mathend(E)}$ is the induced $L$-connection on $\mathend(E)$ evaluated on $\varphi\in \Gamma\big(\End(E)\big),\,e\in \Gamma(E)$ as:
		\begin{equation*}\label{eq:nablaEndE}
			\Big(\nabla^{\mathend(E)}_l\varphi\Big)(e)\triangleq\nabla_l\big(\varphi(e)\big)-\varphi\big(\nabla_l(e)\big), \quad\text{i.e.}, \nabla^{\mathend(E)}_l\varphi=[\nabla_l,\varphi].
		\end{equation*}
	\end{enumerate}
	After the transfer, the identifications in~\eqref{eq:LdiffE=EndE+L} automatically become Lie algebroid isomorphisms.
\end{thm}

\begin{proof}
	The transferred anchor $\rho^\nabla(\varphi \oplus l)$ is expected to be $\rho_{_{\LdiffE}}\big(\varphi+\nabla_{_l},l\big)=\rho_{_L} (l)$.
	The transferred Lie bracket $[\varphi_1 \oplus l_1,\varphi_2\oplus l_2]^\nabla$ is expected to be related with the corresponding Lie bracket in $\Gamma\big(\LdiffE\big)$:
	\begin{align*}
			\big[\big(\varphi_1+\nabla_{_{l_1}},l_1\big),\big(\varphi_2+\nabla_{_{l_2}},l_2\big)\big]
			&=\Big([\varphi_1,\varphi_2]+[\nabla_{_{l_1}},\varphi_2]+[\varphi_1,\nabla_{_{l_2}}]+[\nabla_{_{l_1}},\nabla_{_{l_2}}],[l_1,l_2]\Big)
		\end{align*}
	which can be identified via~\eqref{eq:LdiffE=EndE+L} back to be an element in $\Gamma\big(\End(E)\oplus L\big)$:
	\begin{align*}
			&\Big([\varphi_1,\varphi_2]+[\nabla_{_{l_1}},\varphi_2]-[\nabla_{_{l_2}},\varphi_1]+[\nabla_{_{l_1}},\nabla_{_{l_2}}]-\nabla_{_{[l_1,l_2]}}\Big) \oplus [l_1,l_2]\\
			=&\Big([\varphi_1,\varphi_2]+\nabla^{\mathend(E)}_{l_1}\varphi_2-\nabla^{\mathend(E)}_{l_2}\varphi_1+R^{\nabla}(l_1,l_2)\Big) \oplus [l_1,l_2].
		\end{align*}
	The Jacobi identity and the Leibniz rule are also spontaneously transferred to hold for $\rho^\nabla,\,[\cdot,\cdot]^\nabla$.
\end{proof}

\begin{defn}\label{def:EndE+L}
	We denote the above Lie algebroid $\big(\mathend(E)\oplus L,[~\cdot,\cdot~]^\nabla,\rho^\nabla\big)$
	by  ${\End(E)\nabplus L}$  for short,
	and call it the \textbf{virtually split $L$-Atiyah algebroid of $E$ with respect to $\nabla$}.
\end{defn}

$L$-Atiyah sequence~\eqref{seq:LdiffE} is then transferred to the below format:
\begin{equation}\label{seq:LdiffE3}
	\xymatrix@C=2pc@R=2pc{
		0\ar[r]&  \mathend(E) \ar[r]^{{\iota}\quad}& {\End(E)\nabplus L} \ar@<0.5ex>[r]^{\,\,\qquad \sigma_{_L}} 
		&L\ar[r]\ar@<0.5ex>@{-->}[l]^{\,\qquad s}& 0,
	}
\end{equation}
where the Lie algebroid morphisms $\iota,\,\sigL$ and the vector-bundle right splitting $s$ take very simple forms: $$\iota(\varphi)=\varphi\oplus 0 ,\quad\sigL(\varphi\oplus l)=l,\quad s(l)=0\oplus l.$$
We call~\eqref{seq:LdiffE3} the \textbf{virtually split $L$-Atiyah sequence of $E$ with respect to $\nabla$}.

\begin{rmk}\label{rmk:indepNabla}
	Though the notation ${\End(E)\nabplus L}$ formally emphasizes its dependence on an $L$-connection $\nabla$,
	it is essentially unique up to Lie algebroid isomorphism  if we have some other $L$-connection $\nabla'$ on $E$,
	$${\End(E)\nabplus L} \cong \LdiffE\cong \End(E)\scalebox{0.8}{$\stackrel{\nabla'}{\oplus}$} L,$$
	due to the Lie algebroid identification $\mathfrak{D}_L(E)\cong {\End(E)\nabplus L} $ in Theorem~\ref{thm:EndE+L}.
\end{rmk}

\begin{rmk}\label{rmk:EndE+L}
	A direct computation of the Jacobi $3$-form of $[~\cdot,\cdot~]^\nabla$ gives
	\begin{align*}
		\big[[\varphi_1 \oplus l_1,\varphi_2\oplus l_2]^\nabla,\varphi_3\oplus l_3\big]^\nabla +c.p.
		= \Big(\big(\big[[\varphi_1 ,\varphi_2],\varphi_3\big]+c.p.\big) - \big(d_L^\nabla R^{\nabla}\big)(l_1,l_2,l_3)\Big) \oplus \big(\big[[l_1,l_2],l_3\big]+c.p.\big)
	\end{align*}
	where $c.p.$ stands for cyclic permutations, and $d_L^\nabla R^{\nabla} \in \Gamma\big(\wedge^3 L^\vee \otimes \mathend(E)\big)$ is the Bianchi $3$-form.
	Thus, we see that the Bianchi identity $d_L^\nabla R^{\nabla}=0$ is equivalent the Jacobi identity for ${\End(E)\nabplus L}$,
	the latter of which is true via the transfer from $\LdiffE$.
\end{rmk}

\begin{rmk} The Lie bracket on a Lie algebra extension, as shown in Theorem~\ref{thm:EndE+L}.(\ref{item:[]_EndE+L}), was explicitly presented by Alekseevsky et al.~\cite{Alek04}. 
They noted that this formula has been known for a long time, at least implicitly by Hochschild and Mackenzie.
Also, the construction of virtually split $L$-Atiyah algebroid/sequence
with respect to a choice of $\nabla$ has appeared as a special case of Abad and Crainc's~\cite[Proposition~3.9]{Abad-Crainic12}.
\end{rmk}



\vskip 10pt
\subsection{$L$-Atiyah algebroids of   principal $G$-bundles}\label{subsec:AtLP}
The classic notion of a connection on a principal $G$-bundle $P$ is defined as a $G$-invariant horizontal distribution in $TP$,
which naturally works with respect to a Lie algebroid $L$ to be a $G$-invariant ``horizontal $L$-distribution'' in $TP$.

\begin{defn}[\cite{Fernandes02}]\label{def:L-ConnP}
	A \textbf{connection} of a Lie algebroid $L$ on a principal $G$-bundle $\pi\colon P\to M$, is a vector bundle map $h\colon \pi^* L \to TP$ with identity map on the base $P$
	such that it is $G$-equivariant and fits in the commutative diagram of vector bundles:
	\begin{equation*}
		\begin{split}
			\xymatrix@C=2.5pc@R=2pc{
				\pi^* L \ar[r]^{h}\ar[d]_{\pi_*} \ar@{}[dr]|(0.5){\scalebox{1.2}{$\circlearrowleft$}}
				& TP \ar[d]^{\pi_*}\,\\
				L\ar[r]_{\rho_{_L}} &	TM.}
		\end{split}
	\end{equation*}
\end{defn}

We will show that the construction of $L$-Atiyah algebroid for a vector bundle and its interpretation in terms of $L$-connections from Subsections~\ref{subsec:DLE},
also hold for a principal $G$-bundle $\pi\colon P\to M$.

\vskip 5pt
\subsubsection{Construction of the $L$-Atiyah algebroid $\At_L(P)$ }

Firstly, we consider the pullback between the tangent map ${T\pi}$ and the pulled anchor $\pi^*(\rho_{_{L}})$, both of which are vector-bundle maps over $P$:
\begin{equation}\label{Dia:LdiffE-pullback1}
	\begin{split}
		\xymatrix@C=2pc@R=2pc{
			0\ar[r] & P\times {\g} \ar[r]\ar@{=}[d]& 
			T^L P
			\ar[r]
			\ar[d]
			\ar@{}[dr]|(0.3){\text{\pigpenfont J}}&
			\pi^*L
			\ar[r]
			\ar[d]^{\pi^*\rho_{_L}}&
			0\,\\
			0\ar[r]&  P\times {\g} \ar[r] &
			TP  \ar[r]^{T\pi\,\,}&
			\pi^*TM \ar[r]& 0.}
	\end{split}
\end{equation}
The pullback results in a vector bundle on $P$ which we denote by $T^L P$.

Next, we take the quotients of these $G$-equivariant bundles by the free $G$-actions to obtain a diagram of vector-bundle maps over $M$:
\begin{equation*}\label{Dia:LdiffE-pullback2}
	\begin{split}
		\xymatrix@C=2pc@R=2pc{
			0\ar[r] & P\times_{_G} {\g} \ar[r]^{\iota}\ar@{=}[d]& 
			(T^L P)/G
			\ar[r]^(0.55){\sigL}
			\ar[d]_{f_{_L}}
			\ar@{}[dr]|(0.3){\text{\pigpenfont J}}&
			L
			\ar[r]
			\ar[d]^{\rho_{_L}}&
			0\,\\
			0\ar[r]&  P\times_{_G} {\g} \ar[r] &
			(TP)/G \ar[r]^(0.55){\widetilde{T\pi}}&
			TM \ar[r]& 0,}
	\end{split}
\end{equation*}
where the right square is a pullback diagram of vector-bundle maps.

Specifically, we have
\begin{equation*}\label{eq:AtLP}
	(T^L P)/G = \Big\{(\xi_m,l_m) \in \bigcup_{m\in M}
	\Big(\big((TP)/G\big)_m\times L_m\Big) \bigm| \widetilde{T\pi}(\xi_m)=\rho_{_L}(l_m) \in T_m M  \Big\},
\end{equation*}
with the section space $\Gamma\big((T^L P)/G\big)=\big\{(\xi,l)\in \Gamma\big((TP)/G\big)\times \Gamma(L) \bigm| \widetilde{T\pi}(\xi)=\rho_{_L}(l)\in \Gamma(TM )\big\}$.
The maps $\iota,\,\sigma_{_L},\,f_{_L}$ are then given by
\begin{equation*}\label{eq:symbol-AtLP}
	\iota(\xi')=(\xi',0),\qquad f_{_L}(\xi,l)=\xi,\qquad \sigma_{_L}(\xi,l)=l,
\end{equation*}
for all $\xi'\in \Gamma\big(P\times_{_G} {\g}\big),\,(\xi,l)\in \Gamma\big((T^L P)/G\big)$.

Note that the maps $\widetilde{T\pi}$ and $\rho_{_{L}}$ are anchors of Lie algebroids.
Hence, the pullback vector-bundle $(T^L P)/G$ is naturally a Lie algebroid on $M$
whose anchor and Lie bracket can be expressed by those of $L$ and the Atiyah algebroid $\At(P)=(TP)/G$, similar to 
\eqref{eq:LdiffE} for $\LdiffE$:
\begin{equation}\label{eq:AtLP}
	\begin{aligned} 
		\big[(\xi_1,l_1),({\xi}_2,l_2)\big]_{(T^L P)/G}&\triangleq\big([\xi_1,{\xi}_2]_{(TP)/G},\,[l_1,l_2]_{_L}\big),\\
		\rho_{_{(T^L P)/G}}(\xi,l)&\triangleq\widetilde{T\pi}(\xi)=\rho_{_L}(l).
	\end{aligned}
\end{equation}
The sequence of vector bundles
\begin{equation}\label{seq:LdiffP}
	\xymatrix@C=2pc@R=2pc{
		0\ar[r]&  P\times_{_G} {\g} \ar[r]^(0.48){\iota}& (T^L P)/G \ar[r]^(0.65){\sigL}&
		L\ar[r]& 0,
	}
\end{equation}
then becomes an exact sequence of Lie algebroids over $M$.

\begin{Def}\label{def:LdiffP}
	We call $\At_L(P):=(T^L P)/G$ the \textbf{$L$-Atiyah algebroid} of the principal $G$-bundle $P$,
	and call sequence~\eqref{seq:LdiffP} the \textbf{$L$-Atiyah sequence} of $P$, which we also frequently write as:
	\begin{equation*}
		\xymatrix@C=2pc@R=2pc{
			0\ar[r]&  \Ad(P) \ar[r]^{\iota}& \At_L(P) \ar[r]^(0.64){\sigL}&
			L\ar[r]& 0.
		}
	\end{equation*}
	
\end{Def}

Let $\Psi\colon L\to L'$ be a Lie algebroid morphism over $M$.
It naturally induces a Lie algebroid morphism
between Atiyah algebroids $\At_L(P)$ and $\At_{L'}(P)$ similar to Proposition~\ref{prop:functoriality}.

\vskip 5pt
\subsubsection{Splitting of the $L$-Atiyah algebroid by an $L$-connection}\label{subsubsec:splitAtiyahP}

As in Theorem~\ref{thm:conn=splitting}, we provide the equivalences of $L$-connections.

\begin{thm} \label{thm:Pconn=splitting}
	The following data are equivalent:
	\begin{enumerate}[(1)]
		\item an $L$-connection on $P$ (Definition~\ref{def:L-ConnP});
		\item \label{item:Pconn=splitting} a vector-bundle right splitting
		$s\colon L\to \At_L(P)$ of $L$-Atiyah sequence~\eqref{seq:LdiffP},
		or equivalently, a vector-bundle left splitting $\theta\colon \At_L(P) \to \Ad(P)$ of~\eqref{seq:LdiffP}
		with the mutual determination ${\iota}\circ\theta+s\circ\sigma_{_L}=id_{(T^L P)/G}$:
		\begin{equation*}\label{seq:LdiffE4}
			\xymatrix@C=2pc@R=2pc{
				0\ar[r]&  \Ad(P) \ar@<0.5ex>[r]^{\iota}&
				{\At_L(P)} \ar@<0.5ex>[r]^(0.64){\sigma_{_L}} \ar@<0.5ex>@{-->}[l]^{\theta}&
				L\ar[r]\ar@<0.5ex>@{-->}[l]^(0.36){s}& 0;
			}
		\end{equation*}
		\item a vector-bundle decomposition $\At_L(P)\cong \Ad(P)\oplus L$ such that $\iota\big(\Ad(P)\big)=\Ad(P)\oplus 0$ and ${\sigL(0\oplus L)=L}$;
		\item \label{item:Pconn=sbar} a vector-bundle map $\bar{s} \colon L\to \At(P)$ that makes the following lower triangle commute:
		\begin{equation*}
			\begin{split}
				\xymatrix@C=2.5pc@R=2pc{
					\At_L(P)
					\ar@<0.5ex>[r]^{\sigL}
					\ar[d]_{f_{_L}}
					\ar@{}[dr]|(0.25){\text{\pigpenfont J}}|(0.7){\scalebox{1.2}{$\circlearrowleft$}}
					& L\, \ar[d]^{\rhoL}
					\ar@{-->}@<0.5ex>[l]^{s} \ar@{-->}[ld]|(0.5){\scalebox{1}{$\bar{s}$}}
					\\
					\At(P) 
					\ar[r]^{\sigma}
					&TM.}
			\end{split}
		\end{equation*}
	\end{enumerate}
\end{thm}
\begin{proof}
	The equivalences $(2)\Leftrightarrow(3)\Leftrightarrow(4)$ are explained the same way as in Theorem~\ref{thm:conn=splitting}.
	For instance, the left and right splittings $s,\,\theta$ give the explicit decomposition $\At_L(P)\cong \Ad(P)\oplus L$, cf.\eqref{eq:LdiffE=EndE+L} :
	\begin{equation}\label{eq:LAtP=AdP+L}
		\begin{aligned}
			\Gamma\big(\At_L(P)\big) &\longrightarrow \Gamma\big(\Ad(P)\oplus L\big):
			&(\xi,l)&\longmapsto \theta(\xi,l)\oplus \sigL(\xi,l) =\big(\xi-(f_{_L}\circ s)(l)\big)\oplus l,\\
			\Gamma\big(\Ad(P)\oplus L\big) &\longrightarrow \Gamma\big(\At_L(P)\big):
			& \xi' \oplus l &\longmapsto {\iota}(\xi')+s(l)=\big(\xi'+(f_{_L}\circ s)(l),\,l\big).
		\end{aligned} 
	\end{equation}
	
	To see $(1)\Leftrightarrow(4)$, after passing to the $G$-quotients,
	a $G$-equivariant vector-bundle map $h\colon \pi^* L \to TP$
	descends to a vector bundle map $\bar{s} \colon L\to (T P)/G$,
	and the commutating diagram in Definition~\ref{def:L-ConnP} amounts to the commuting lower triangle in (4).
\end{proof}

Consequently, instead of Definition~\ref{def:L-ConnP},
we adopt the following equivalent but more convenient definition of a connection using Theorem~\ref{thm:Pconn=splitting}.(\ref{item:Pconn=splitting}).

\begin{Def}\label{def:Pconn=splitting}
	An \textbf{$L$-connection} on the principal $G$-bundle $P$ is a vector-bundle right splitting of $L$-Atiyah sequence~\eqref{seq:LdiffP}:
	$${s\colon L\longto \At_L(P)}.$$
\end{Def}

The measurement of whether $s$ preserves Lie brackets:
\[R^s\colon \Gamma(L)\times \Gamma(L)\longto \Gamma\big(\At_L(P)\big),\quad
(l_1,l_2)\longmapsto [s(l_1),s(l_2)]-s([l_1,l_2]),
\]
is actually $C^\infty(M)$-linear, because $s$ respects anchors as shown by the diagram in Theorem~\ref{thm:Pconn=splitting}.(\ref{item:Pconn=sbar}).
In addition, applying the symbol map $\sigL$, 
we see the images of $R^s$ all have vanishing $L$-symbols meaning that $R^s$ takes values in $\Ad(P)$ by $L$-Atiyah sequence~\eqref{seq:LdiffP}.

\begin{Def}\label{def:Pcurv}
	The $C^\infty(M)$-linear map
	$$R^s\colon \Gamma(L)\times \Gamma(L)\longto \Gamma\big(\Ad(P)\big)$$
	as defined above is the \textbf{curvature} of the $L$-connection $s$ on $P$.
	If $R^s=0$, the connection $s$ is said to be \textbf{flat}.
\end{Def}

\begin{cor} With notations as in Theorem~\ref{thm:Pconn=splitting} and the above definitions, we have the following facts:
	\begin{enumerate}[(a)]
		\item $L$-connections $s$ on $E$ always exist, and they form an affine space modelled on $\Gamma\big(\Hom\big(L,\Ad(P)\big)\big)$;
		\item an $L$-connection $s\colon L \to \At_L(P)$ is flat if and only if 
		it is a Lie algebroid morphism making $L$ a Lie subalgebroid of $\At_L (P)$,
		equivalently, if and only if $\bar{s}\colon L\to \At(P)$ is a Lie algebroid morphism.
	\end{enumerate}
\end{cor}

The principal bundle version of Theorem~\ref{thm:EndE+L} is given below.

\begin{thm}\label{thm:AdP+L}
	Given an $L$-connection $s$ on $P$, the identifications in~\eqref{eq:LAtP=AdP+L} transfer Lie algebroid structure~\eqref{eq:AtLP} of $\At_L(P)$ to the Lie algebroid structure of $\Ad(P)\oplus L$ as follows:
	\begin{enumerate}[i)]
		\item the anchor $\rho^s\colon \Ad(P)\oplus L\to TM$ of $\xi'\oplus l \in\Gamma\big(\Ad(P)\oplus L\big)$ is
		\begin{equation*}\label{eq:rho_AdP+L}
			\rho^s(\xi' \oplus l)\triangleq \rho_{_L}(l);
		\end{equation*}
		\item the Lie bracket $[~\cdot,\cdot~]^s\colon \Gamma\big(\Ad(P)\oplus L\big)\times\Gamma\big(\Ad(P)\oplus L\big)\to\Gamma\big(\Ad(P)\oplus L\big)$ for $\xi'_i\oplus l_i \in\Gamma\big(\Ad(P)\oplus L\big)$ reads
		\begin{equation*}\label{eq:[]_AdP+L}
			[\xi'_1 \oplus l_1,\xi'_2\oplus l_2]^s
		\triangleq\Big([\xi'_1,\xi'_2]+\nabla^{\Ad (P)}_{l_1}\xi'_2-\nabla^{\Ad (P)}_{l_2}\xi'_1+R^s(l_1,l_2)\Big) \oplus [l_1,l_2],
		\end{equation*}
		where $\nabla^{\Ad (P)}$ is the induced $L$-connection on the associated adjoint vector bundle $\Ad (P)=P\times_{_G} {\g}$ evaluated on a section $\xi'\in \Gamma\big(\Ad(P)\big)$ as:
		\begin{equation*}\label{eq:nablaAdP}
			\nabla^{\Ad (P)}_l\xi'\triangleq [s(l),\xi']\in \Gamma\big(\Ad(P)\big)\subseteq \Gamma\big(\At_L(P)\big).
		\end{equation*}
		\end{enumerate}
		After this transfer, the identifications in~\eqref{eq:LAtP=AdP+L} automatically become Lie algebroid isomorphisms.
\end{thm}

\begin{defn}\label{def:AtP+L}
	We denote the above Lie algebroid $\big(\Ad(P)\oplus L,\,[~\cdot,\cdot~]^s,\,\rho^s\big)$
	by  ${\Ad(P)\splus L}$  for short,
	and call it the \textbf{virtually split $L$-Atiyah algebroid of $P$ with respect to $s$}.
\end{defn}

$L$-Atiyah sequence~\eqref{seq:LdiffP} of $P$ is then transferred to be:
\begin{equation}\label{seq:LdiffP1}
	\xymatrix@C=2pc@R=2pc{
		0\ar[r]&  \Ad(P) \ar[r]^{{\iota}\quad}& {\Ad(P)\splus L} \ar@<0.5ex>[r]^{\,\,\qquad \sigma_{_L}} 
		&L\ar[r]\ar@<0.5ex>@{-->}[l]^{\,\qquad s}& 0,
	}
\end{equation}
where the Lie algebroid morphisms $\iota,\,\sigL$ and the vector-bundle right splitting $s$ take very simple forms: $$\iota(\xi')=\xi'\oplus 0 ,\quad\sigL(\xi'\oplus l)=l,\quad s(l)=0\oplus l.$$
We call Sequence~\eqref{seq:LdiffP1} the \textbf{virtually split $L$-Atiyah sequence of $P$ with respect to $s$}.

In the remaining part of this paper, we mainly consider the Atiyah constructions of a vector bundle $E$, 
but most of our results are immediately applicable to the Atiyah constructions of a principal $G$-bundle $P$. 

\vskip 10pt

\subsection[Atiyah algebroids for Lie algebroid representations]{Atiyah algebroids for Lie algebroid representations, and a universal construction}\label{Sec:Atiyahalgebroid4reps}
If a vector bundle $E$ or a principal bundle $P$ has a flat $L$-connection, then the virtually split Atiyah algebroid actually splits as we shall see.

\vskip 5pt
\subsubsection{Atiyah algebroids with respect to flat connections}

Let $E$ be a representation of a Lie algebroid $A$, with a flat connection $\barnab$.
By Corollary~\ref{cor:conn=splitting},
the $A$-Atiyah algebroid $\fD_A(E)$ has a right splitting $s_{_A}:A\to\fD_A(E)$ 
which is a Lie algebroid inclusion and can be evaluated on sections as:
\begin{align}\label{eq:sA}
	s_{_A}\colon \Gamma(A)\longto \Gamma\big(\fD_A(E)\big),\quad a\longmapsto (\barnab_a,a),
\end{align}
where $\Gamma\big(\fD_A(E)\big)=\big\{(\delta,a)\in \Gamma\big(\fD(E)\big) \times \Gamma(A) \mid \sigma(\delta)=\rho_{_A}(a)\big\}$.

The split $A$-Atiyah algebroid ${\End(E) \scalebox{0.8}{$\stackrel{\barnab}{\oplus}$}  A}$
contains both $\End(E)$ and $A$ as Lie subalgebroids,
and by Theorem~\ref{thm:EndE+L}.(\ref{item:[]_EndE+L}) the Lie bracket between $a\in \Gamma(A)$ and $\varphi\in\Gamma\big(\End(E)\big)$ is 
\begin{align}\label{eq:BottEndE}[a,\varphi]=\barnab^{\End(E)}_{a} \varphi= \barnab_a\circ \varphi-\varphi\circ \barnab_a \in\Gamma\big(\End(E)\big),
\end{align}
where $\barnab^{\End(E)}$ is a flat $A$-connection on $\End(E)$, due to the flatness of $\barnab$ on $E$.
Thus, we have:

\begin{prop}\label{prop:semiDirect}
	${\End(E) \scalebox{0.8}{$\stackrel{\barnab}{\oplus}$} A}$ is a semi-direct product Lie algebroid $\End(E) \rtimes A$,
	with anchor $\rhoA$ and the below Lie bracket:
	\[
	[\varphi_1\oplus a_1,\varphi_2\oplus a_2]=\big([\varphi_1,\varphi_2]+\barnab^{\End(E)}_{a_1} \varphi_2-\barnab^{\End(E)}_{a_2} \varphi_1\big)\oplus[a_1,a_2],
	\]
	for all $a_1, a_2\in\Gamma(A),\,\varphi_1,\varphi_2\in\Gamma\big(\End(E)\big)$.
\end{prop}

\begin{Ex}\label{ex:holomEnd(E)}
	Given a holomorphic vector bundle $\E$ on a complex manifold $M$, there is naturally a flat $T^{0,1} M$-connection $\parbar$ on the underlying smooth vector bundle $E$.
	Hence, the $T^{0,1} M$-Atiyah algebroid of $E$ has a split form $\End(E)\rtimes T^{0,1} M$.
\end{Ex}

In the case of a principal $G$-bundle $P$ with a flat $A$-connection,
i.e., a Lie algebroid inclusion ${\bar{s}:A\to  \At_A(P)}$,
we can identify $\At_A(P)$ by a split vector bundle, similar to~\eqref{eq:LdiffE=EndE+L}:
\[
\At_A(P) \cong \Ad(P)\overset{\bar{s}}{\oplus}A: \quad \xi \longmapsto \big(\xi-\bar{s}\circ\sigA(\xi)\big)\oplus\sigA(\xi).
\]
As in Proposition~\ref{prop:semiDirect} above, we obtain

\begin{prop}
	$\Ad(P)\overset{\bar{s}}{\oplus}A$ is a semi-direct product Lie algebroid $\Ad(P)\rtimes A$,
	with anchor $\rhoA$ and the below Lie bracket for $\xi_i\oplus a_i \in \Gamma\big(\Ad(P)\rtimes A\big)$:
	\[
	[\xi_1\oplus a_1,\xi_2\oplus a_2]=\big([\xi_1,\xi_2]+\barnab^{\Ad(P)}_{a_1}\xi_2-\barnab^{\Ad(P)}_{a_2}\xi_1\big)\oplus [a_1,a_2],
	\]
	where $\barnab^{\Ad(P)}_{a}\xi=[\bar{s}(a),\xi]\in \Gamma\big(\Ad(P)\big)\subseteq \Gamma\big(\At_A(P)\big)$
	is the induced flat $A$-connection on the adjoint bundle~$\Ad(P)$.
\end{prop}

\vskip 5pt
\subsubsection{A universal construction}

For a Lie algebroid $L$ and a vector bundle $E$ over $M$, the Atiyah bundle $\LdiffE$ itself is a Lie algebroid, so we can consider the Atiyah algebroid $\fD_{\LdiffE}(E)$ whose section space is
\[
\Gamma\big(\fD_{\LdiffE}(E)\big)=\big\{(\delta',(\delta,l))\in \Gamma\big(\diffE\big)\times \Gamma\big(\LdiffE\big) \bigm | \sigma(\delta')=\sigma(\delta)=\rho_{_L}(l)\in \Gamma(TM )\big\}
\]
with the anchor map $(\delta',(\delta,l))\mapsto \rho_{_L}(l)$ and the Lie bracket using~\eqref{eq:LdiffE}:
\[
\big[(\delta'_1,(\delta_1,l_1)),(\delta'_2,(\delta_2,l_2))\big]=\big([\delta'_1,\delta'_2],([\delta_1,\delta_2],[l_1,l_2])\big).
\]
We will show that $\fD_{\LdiffE}(E)$ has a natural semi-direct product Lie algebroid structure and satisfies a universal property for $L$-connections on $E$.

By Proposition~\ref{prop:functoriality} on the functoriality of Atiyah algebroids, 
the anchor ${\rho\colon L\to TM}$ induces a Lie algebroid morphism $f_{_L}\colon \LdiffE\to \diffE$ .
Then applying Theorem~\ref{thm:conn=splitting}, we have the following commuting lower triangle of Lie algebroid morphisms
\begin{equation*}
	\xymatrix@C=2.5pc@R=2pc{
		\fD_{\LdiffE}(E)
		\ar@<0.5ex>[r]^{\sigma_{_{\LdiffE}}}
		\ar[d]_{f_{_{\LdiffE}}}
		\ar@{}[dr]|(0.25){\text{\pigpenfont J}}|(0.7){\scalebox{1.2}{$\circlearrowleft$}}
		& \LdiffE\, \ar[d]^{\sigL}
		\ar@{-->}@<0.5ex>[l]^{S} \ar@{-->}[ld]|(0.5){{f_{_L}}}
		\\
		\diffE 
		\ar[r]^{\sigma}
		&TM,} \qquad \quad
	\xymatrix@C=2pc@R=2pc{
		(\delta',(\delta,l))
		\ar[r]^{\sigma_{_{\LdiffE}}}
		\ar[d]_{f_{_{\LdiffE}}}
		& (\delta,l) \ar[d]^{\sigL}
		\\
		\delta' 
		\ar[r]^(0.28){\sigma}
		&\sigma(\delta')=\sigma(\delta)=\rhoL(l).}
\end{equation*}
The corresponding right splitting
$S\colon \LdiffE\to \fD_{\LdiffE}(E)$ respects the Lie algebroid structures, and hence provides a canonical flat $\LdiffE$-connection on $E$.

Specifically, the commutative diagram represented by solid arrows is shown in the above right square, whereas the dotted arrows $f_{_L}$ and $S$ are
\[
f_{_L}(\delta,l)=\delta,\qquad S(\delta,l)=(\delta,(\delta,l)),\qquad \forall\,(\delta,l)\in \Gamma\big(\LdiffE\big).
\]
The \textbf{canonical flat $\LdiffE$-connection} on $E$, denoted by $\breve{\nabla}$, can be given as
\[
\breve{\nabla}_{(\delta,l)}(e)=f_{_L}(\delta,l)(e)=\delta(e), \qquad \forall\, (\delta,l)\in \Gamma\big(\LdiffE\big),\,e\in \Gamma(E).
\]

Since $\breve{\nabla}$ is flat,
formula~\eqref{eq:LdiffE=EndE+L} and Proposition~\ref{prop:semiDirect} implies:

\begin{prop}
	$\fD_{\LdiffE}(E)$ is a semi-direct product Lie algebroid:
	\[
	\fD_{\LdiffE}(E)\cong \End(E)\scalebox{0.8}{$\stackrel{\breve{\nabla}}{\oplus}$} \LdiffE = \End(E)\rtimes \LdiffE:
	\quad (\delta',(\delta,l))\longmapsto (\delta'-\delta) \oplus (\delta,l).
	\]
	Here, $\End(E)\rtimes \LdiffE$ has the anchor $\rhoL$ and the Lie bracket for $\varphi_i\oplus (\delta_i,l_i) \in \Gamma\big(\End(E)\rtimes \LdiffE\big)$ to be:
	\[
	[\varphi_1\oplus (\delta_1,l_1),\varphi_2\oplus (\delta_2,l_2)]=
	\big([\varphi_1,\varphi_2]+\breve{\nabla}_{(\delta_1,l_1)}^{\End(E)}(\varphi_2)-\breve{\nabla}_{(\delta_2,l_2)}^{\End(E)}(\varphi_1)\big)\oplus \big([\delta_1,\delta_2],[l_1,l_2]\big)
	\]
	where $\breve{\nabla}_{(\delta,l)}^{\End(E)}(\varphi)=[\delta,\varphi]=\delta\circ\varphi-\varphi\circ\delta\in\Gamma\big(\End(E)\big)$ is a flat $\LdiffE$-connection on $\End(E)$.
\end{prop}

The universal property of $\breve{\nabla}$ is that, for any $L$-connection $\nabla$ on $E$ 
and its corresponding vector-bundle right splitting $s\colon L\to \LdiffE\colon l\mapsto(\nabla_l,l)$, we have the relation:
\[
\nabla_l=\breve{\nabla}_{s(l)}, \qquad \forall\,l\in\Gamma(L).
\]
In other words, we get a naturality property:

\begin{prop}
	The $L$-connection $\nabla$ is pulled back from $\breve{\nabla}$ via $s$.
\end{prop}

Likewise, $\At_{{\At_L(P)}}(P)$ also has a natural semi-direct product Lie algebroid structure and satisfies a universal property for $L$-connections on the principal bundle $P$.

\vskip 25pt

\section{Atiyah extensions, complementary connections, and the compatibility conditions} \label{Sec:Atiyahextension}

\vskip 10pt

In this section, we fix a Lie triad  $(L,A,E)$ over $M$
with a flat $A$-connection $\barnab$ on $E$,
and denote by $B$ the  quotient bundle $L/A$
which is a canonical $A$-representation via the Bott connection. 
We will introduce three $A$-representations that realize the Atiyah class $\alpha(L,A,E)$ as $A$-representation extensions.
They are:
\begin{enumerate}[(I)]
	\item the \textbf{Atiyah extension} $\BdiffE$,
	\item the \textbf{embedded Atiyah extension} $\BdiffI$,
	\item the \textbf{virtually split Atiyah extension} $\AtExtE$,
\end{enumerate}
among which $\BdiffE$ is independent of the choices of a splitting $\iB\colon B\to L$ and  an $L$-connection $\nabla$ extending~$\barnab$, 
while $\BdiffI$ and $\AtExtE$ depend respectively on the choice of  $\iB$ 
and the choice of $\nabla$. 

We shall show that these Atiyah extensions are isomorphic as $A$-representations,
point out their duality with the jet bundle~${J}^1_{B}(E)$,
and reveal how they directly represent the obstruction to the existence of compatible connections just as the Atiyah class does.

Furthermore, for a principal bundle version of Lie triad $(L,A,P)$ (see Remark~\ref{rmk:LAP}), the pertinent constructions and properties of Atiyah extensions, supposedly denoted by $\At_B(P),\,\At_{\iB}(P),\mbox{ and } \Ad(P)\nabplus B$, will be promptly attainable. So we skip the details of that version. 

\vskip 10pt
\subsection{First construction of Atiyah extensions as  quotient bundles}\label{subsec:relativeAtiyah}

For the Lie triad $(L,A,E)$, we already have two Atiyah algebroids, namely $\fD_A(E)$, and $\fD_L(E)$.
We shall also associate to the quotient bundle $B$ a similar construction, to be denoted by $\BdiffE$.

By Proposition~\ref{prop:functoriality} on the functoriality of Atiyah algebroids,
the Lie algebroid inclusion $i_{_A}\colon A\to L$ 
induces a morphism between the corresponding Atiyah sequences:

\begin{equation}\label{Dia:DALE-pullback}
	\begin{split}
		\xymatrix@C=2pc@R=2pc{
			0 \ar[r]
			&\mathend(E) \ar[r]^{{\iota}}\ar@{=}[d]
			&\mathfrak{D}_A(E) \ar@<0.5ex>[r]^{\quad\sigma_{_A}}\ar[d]_{f_{_{A,L}}}
			\ar@{}[dr]|(0.25){\text{\pigpenfont J}}
			&A \ar[r]\ar[d]^{i_{_A}} \ar@<0.5ex>[l]^{\quad s_{_A}}
			&0\,
			\\
			0 \ar[r]
			&\mathend(E) \ar[r]^{{\iota}}
			& \mathfrak{D}_L(E) \ar[r]^{\quad\sigma_{_L}}
			&L\ar[r]
			&0.
		}
	\end{split}
\end{equation}
The induced Lie algebroid morphism $f_{_{A,L}}$ is explicitly given as:
\begin{align}\label{eq:DAEtoDLA}
	f_{_{A,L}}\colon \Gamma\big(\fD_A(E)\big)\longrightarrow \Gamma\big(\fD_L(E)\big), \quad (\delta,a)\longmapsto (\delta,i_{_A}(a)).
\end{align}

\begin{prop}\label{prop:LiePairDLAE}
	The Lie algebroid morphism $f_{_{A,L}}\colon \fD_A(E)\rightarrow \fD_L(E)$ is injective,
	making $\big(\fD_L(E),\fD_A(E)\big)$ a Lie pair with the quotient bundle $\fD_L(E)/\fD_A(E)\cong B$.
	The quotient sequence is:
	\begin{equation}\label{seq:LiePairDLAE}
		\xymatrix@C=2.5pc@R=2pc{
			0 \ar[r]&\fD_A(E)\ar[r]^{ f_{_{A,L}}}&\mathfrak{D}_L(E)\ar[r]^{\quad{\rm pr}_{_B}\circ \sigma_{_L}}
			&B\ar[r]&0.
		}
	\end{equation}
\end{prop}

\begin{proof}
	Applying the Snake Lemma to the Atiyah sequences 
	with respect to $A$ and $L$ in Diagram~\eqref{Dia:DALE-pullback},
	we have an expanded commutative diagram with extra two vertical short exact sequences:
	\begin{equation}\label{dia:snakeOfAtiyah}
		\begin{split}
			\xymatrix@C=2.5pc@R=2pc{
				&&0\ar[d]&0\ar[d]&\\
				0 \ar[r]&\End(E)\ar[r]^{{\iota}}\ar@{=}[d]&
				\mathfrak{D}_A(E)\ar[r]^{\quad\sigma_{_A}}\ar[d]_{f_{_{A,L}}}&A\ar[r]\ar[d]^{i_{_A}}&0\\
				0 \ar[r]&\End(E)\ar[r]^{{\iota}}&
				\mathfrak{D}_L(E)\ar[r]^{\quad\sigma_{_L}}\ar[d]
				&L\ar[r]\ar[d]^{{\rm pr}_{_B}}&0\\
				&&B\ar@{=}[r]\ar[d]&B\ar[d]&\\
				&&0&0&
			}
		\end{split}
	\end{equation}
	where the kernel and the cokernel of $f_{_{A,L}}$ are inherited from $i_{_A}$.
\end{proof}

The Lie algebroid morphism $f_{_{A,L}}\circ s_{_A}\colon A\to \LdiffE$ composed from two inclusions~\eqref{eq:sA}~\eqref{eq:DAEtoDLA} is given by
\begin{align}\label{eq:A-DLE}
	s_{_{A,L}}\triangleq f_{_{A,L}}\circ s_{_A}\colon\Gamma(A)\overset{s_{_A}}{\longrightarrow} \Gamma\big(\fD_A(E)\big)\overset{f_{_{A,L}}}{\longrightarrow} \Gamma\big(\fD_L(E)\big),  \quad a\longmapsto (\barnab_a,a)\longmapsto\big(\barnab_a,i_{_A}(a)\big),
\end{align}
further making a Lie pair $\big(\fD_L(E),A\big)$ with a quotient map denoted by $P_{_B}$ 
and the quotient sequence:
\begin{equation}\label{seq:DLE-DBE}
	\xymatrix@C=2pc@R=2pc{
		0 \ar[r]&A\ar[r]^{s_{_{A,L}}\quad}&\mathfrak{D}_L(E)\ar[r]^{P_{_B}\,\,}
		&\fD_L(E)/A\ar[r]&0.
	}
\end{equation}
The quotient bundle $\fD_L(E)/A$ consists of equivalence classes from $\fD_L(E)$ which we write in double parentheses $\llparen\delta,l\rrparen\triangleq P_{_B}(\delta,l)$. So, its section space is:
\begin{align}\label{eq:ΓDBE}
	\begin{split}
	&\Gamma\big(\fD_L(E)/A\big)
	\cong \Gamma\big(\fD_L(E)\big)/\Gamma(A)\\
	=&\Big\{\llparen\delta,l\rrparen \bigm| (\delta,l)\in \Gamma\big(\fD_L(E)\big) \mbox{ s.t. } \llparen\delta,l\rrparen=\llparen\delta',l'\rrparen \mbox{ iff } (\delta,l)=(\delta',l')+(\barnab_a,a) \mbox{ for some } a\in \Gamma(A)\Big\}.
	\end{split}
\end{align}

According to~\eqref{eq:Bott}, we obtain that

\begin{prop}\label{prop:DBE}
	The flat Bott connection of $A$ on $\fD_L(E)/A$ is
	\begin{align}\label{eq:BottDBE}
		D_{a}\llparen\delta,l\rrparen=
		\big(\!\big([\barnab_a,\delta],[a,l]\big)\!\big),
		\qquad \forall\, a\in \Gamma(A),\,\llparen\delta,l\rrparen\in\Gamma\big(\fD_L(E)/A\big).
	\end{align}
\end{prop} 

\begin{thm}\label{thm:DBE}
	The Bott $A$-representation $\BdiffE \triangleq \mathfrak{D}_L(E)/A$ fits in an Atiyah-type short exact sequence of $A$-representations
	\begin{equation}\label{seq:AtiyahDBE}
		\xymatrix@C=2pc@R=2pc{
			0 \ar[r]&\End(E)\ar[r]^{\iota_{_{B}}}&\BdiffE\ar[r]^{\qquad\sigma_{_B}}
			&B\ar[r]&0,
		}
	\end{equation}
	where $\iota_{_{B}}(\varphi)=\llparen\varphi,0\rrparen,\,\forall\,\varphi\in \Gamma\big(\End(E)\big)$ and $\sigma_{_B}\llparen \delta,l\rrparen=\pr_{_B}(l),\,\forall\,\llparen \delta,l\rrparen\in \Gamma\big(\fD_L(E)/A\big)$.
\end{thm}

\begin{proof}
	The composition $\sigma_{_L} \circ s_{_{A,L}}\colon A{\to} \fD_L(E){\to} L \colon a\mapsto \big(\barnab_a,i_{_A}(a)\big)\mapsto i_{_A}(a)$ is exactly $i_{_A}$,
	which provides a morphism between Lie pairs $\big(\fD_L(E),A\big)$ and $(L,A)$.
	Thus, by Proposition~\ref{prop:morphismOfBottConn}, we get a morphism between their quotient sequences and a morphism $\sigma_{_B}$ between the Bott $A$-connections on $\fD_L(E)/A=\BdiffE$ and $L/A=B$. See the horizontal short exact sequences and their morphism in Diagram~\eqref{dia:snakeOfLiePair1}.
	
	Then, applying the Snake Lemma to these two horizontal short exact sequences, 
	we obtain an expanded commutative diagram with extra two vertical short exact sequences:
	\begin{equation}\label{dia:snakeOfLiePair1}
		\begin{split}
			\xymatrix@C=2.5pc@R=2pc{
				&&0\ar[d]&0\ar[d]&\\
				&&\End(E)\ar@{=}[r]\ar[d]^{{\iota}}&\End(E)\ar[d]^{\iota_{_B}}&\\
				0 \ar[r]&A\ar[r]^{s_{_{A,L}}\quad}\ar@{=}[d]
				&\mathfrak{D}_L(E)\ar[r]^{P_{_B}}\ar[d]^{\sigma_{_L}}
				&\BdiffE\ar[r]\ar[d]^{\sigma_{_B}}&0\\
				0 \ar[r]&A\ar[r]^{i_{_A}\quad}&L\ar[r]^{{\rm pr}_{_{B}}}\ar[d]&B\ar[r]\ar[d]&0\\
				&&0&0&
			}
		\end{split}
	\end{equation}
	where the kernel and the cokernel of $\sigB$ are inherited from $\sigL$.

	The expressions of $\iota_{_B}=P_{_B}\circ {\iota}$ and $\sigma_{_B}$ are obtained by diagram chasing.
	One easily verifies that the $A$-representation structures of $\End(E),\,\BdiffE$ and $B$ as given by ~\eqref{eq:BottEndE},~\eqref{eq:BottDBE} and~\eqref{eq:Bott} are preserved by $\iota_{_B}$ and $\sigma_{_B}$.
\end{proof}

\begin{defn}\label{def:DBE}
	We call Sequence~\eqref{seq:AtiyahDBE} the \textbf{$B$-Atiyah sequence of $E$},
	and call $\BdiffE$ the \textbf{$B$-Atiyah extension/bundle associated with $(L,A,E)$} .
\end{defn}

\begin{rmk}\label{rmk:DBE-Jet}
	To represent the Atiyah class $\alpha(L,A,E)$, Chen, Sti\'enon, and Xu ~\cite{CMP16} constructed the jet bundle ${J}^1_{B}(E)$ which is a natural $A$-representation that fits in the jet sequence:
	\begin{equation}\label{seq:JetB}
	\xymatrix@C=2pc@R=2pc{
		0 \ar[r]& {B^\vee}\otimes E \ar[r]&{J}^1_{B}(E)\ar[r]
		&E\ar[r]&0,
	}
	\end{equation}	  
	The relation between ${J}^1_{B}(E)$ and the Atiyah extension $\BdiffE$ can be seen from the perspective of $E$-duality, (cf.~\cite{Chen-Liu-Omni-Lie-algebroids}).
	
\end{rmk}

Similar to Proposition~\ref{prop:functoriality}, we state the functoriality of $B$-Atiyah sequences:

\begin{prop}\label{prop:functoriality1}
	Let $(\Psi_L,\Psi_A)$ be a morphism between Lie pairs $(L,A),\,(L',A')$ over $M$,
	and suppose that $\Psi_A:A\to A'$ intertwines the flat $A$-connection with the flat $A'$-connection  on $E$,
	i.e., ${\barnab_a={\barnab}_{\Psi_A(a)},\,\forall\, a\in\Gamma(A)}$.
	Denote $L/A=B,\,L'/A'=B'$.
	Then, there induces a morphism between relative Atiyah sequences:
	\begin{equation*}\label{diag:functoriality1}
		\begin{split}
			\xymatrix@C=2pc@R=2pc{
				0 \ar[r]
				&\mathend(E) \ar[r]^{\iota_{_B}}\ar@{=}[d]
				&\BdiffE \ar[r]^{\quad\sigma_{_B}}\ar[d]_{f_{_{\Psi_B}}}
				\ar@{}[dr]|(0.25){\text{\pigpenfont J}}
				&B \ar[r]\ar[d]^{\Psi_B}
				&0\,\,\\
				0 \ar[r]
				&\mathend(E) \ar[r]^{\iota_{_{B'}}}
				& \mathfrak{D}_{B'}(E) \ar[r]^{\quad\sigma_{_{B'}}}
				&B'\ar[r] 
				&0,
			}
		\end{split}
	\end{equation*}
	where the vertical arrows
	are morphisms between $A$-representations and $A'$-representations in the sense of Definition~\ref{dfn:morphConn}, and the middle one is explicitly given by
	$$f_{_{\Psi_B}}\colon \Gamma\big(\BdiffE\big)\longrightarrow \Gamma\big(\mathfrak{D}_{B'}(E)\big), 
	\quad\llparen\delta,l\rrparen\longmapsto \big(\!\big(\delta,\Psi_L(l)\big)\!\big).$$
	Moreover, the right square in the above diagram is a pullback of vector bundles on $M$.
\end{prop}

Further functorial relations among $\LdiffE$, $\fD_A(E)$ and $\BdiffE$ can be found in Appendix~\hyperlink{App:A}{A}.

\vskip 10pt
\subsection{Second construction of Atiyah extensions as sub-bundles}\label{subsec:embedAtiyah}

Fix a decomposition $L\cong A\oplus B$ as in~\eqref{Seq:ABdecomposition}
via a vector-bundle right splitting $i_{_B}\colon B\to L$
and its corresponding  left splitting $\pr_{_A}\colon L \to A$.
We will show that $i_{_B}$ leads to an embedding of the Atiyah extension $\BdiffE$ into the Atiyah algebroid~$\LdiffE$.

Recall pullback construction~\eqref{eq:DLE} of $\mathfrak{D}_L(E)$.
There is a similar pullback of $\diffE$ along $\rho_{_B}=\rho_{_L}\circ i_{_B}$ :
\begin{equation*}\label{eq:DBE}
	\begin{split}
		\fD_{i_{_B}}(E)\triangleq\Big\{(\delta_m,b_m) \in \bigcup_{x\in M}
		\big(\diffE_x\times B_x\big) \bigm| \sigma(\delta_m)=\rho_{_B}(b_m) \in T_m M\Big\},
	\end{split}
\end{equation*}
which is a vector bundle over $M$ and has the section space:
\begin{align}\label{eq:ΓDBEI}
	\Gamma\big(\BdiffI\big)=\big\{(\delta,b)\in \Gamma\big(\diffE\big)\times \Gamma(B) \mid \sigma(\delta)=\rho_{_B}(b)\big\}.
\end{align}

Again, we have an Atiyah-type sequence for $\BdiffI$ projecting onto $B$ with the kernel $\End(E)$.
Due to the functoriality of pullbacks, the inclusion $i_{_B}$ induces a morphism between Atiyah sequences:
\begin{equation}\label{Dia:BdiffE-pullback}
	\begin{split}
		\xymatrix@C=2pc@R=2pc{
			0\ar[r]& \mathend(E)
			\ar[r]^{\tilde{\iota}_{_B}}
			\ar@{=}[d]& 
			\BdiffI
			\ar[r]^{\quad\widetilde{\sigma}_{_B}}
			\ar[d]_{f_{_{\iB}}}
			\ar@{}[dr]|(0.3){\text{\pigpenfont J}}&
			B\ar[r]\ar[d]^{\iB}&
			0\,\,\\
			0 \ar[r]
			&\mathend(E) \ar[r]^{{\iota}}
			& \mathfrak{D}_L(E) \ar[r]^{\quad\sigma_{_L}}
			&L\ar[r] 
			&0,}
	\end{split}
\end{equation}
where explicitly we have for all $m\in M$ and $\varphi_m\in \End(E_m),\,(\delta_m,b_m)\in\BdiffI_m$ that
\begin{equation}\label{eq:symbol-DBE}
	\tilde{\iota}_{_B}(\varphi_m)=(\varphi_m,0),\qquad \widetilde{\sigma}_{_B}(\delta_m,b_m)=b_m,\qquad f_{_{\iB}}(\delta_m,b_m)=\big(\delta_m,\iB(b_m)\big).
\end{equation}

By the above last expression, the map $f_{_{\iB}}:\BdiffI\longrightarrow \LdiffE$
induced from the embedding $\iB:B\to L$ is also a vector-bundle embedding.
In fact, one checks the following sequence is split.

\begin{lem}\label{lem:DBEI}
	$\BdiffI$ fits in a quotient sequence of the Lie pair $\big(\LdiffE,A\big)$:
	\begin{equation}\label{seq:DLE-DBE2}
		\xymatrix@C=2.5pc@R=2pc{
			0 \ar[r]&A\ar[r]^{s_{_{A,L}}\quad}&\mathfrak{D}_L(E)\ar@<0.5ex>[r]^{\widetilde{P}_{_B}\,\,}
			&\BdiffI\ar[r]\ar@{-->}@<0.5ex>[l]^{f_{_{\iB}}}&0,
		}
	\end{equation}
	where $\widetilde{P}_{_B}(\delta,l)=\big(\delta-\parbar_{\prA(l)},\prB(l)\big), \,\forall\, (\delta,l)\in \Gamma\big(\LdiffE\big)$,
	and $f_{_{\iB}}$ is a vector-bundle right splitting.
\end{lem}

Exact sequences~\eqref{seq:DLE-DBE} and~\eqref{seq:DLE-DBE2}
realize both $\BdiffE$ and $\BdiffI$ as $\coker s_{_{A,L}}$.
There is an explicit identification $\BdiffI \cong \BdiffE$ together with its converse evaluated on sections as:
	\begin{equation}\label{eq:BdiffI=BdiffE}
		\begin{aligned}
			\Gamma\big(\BdiffI\big) &\overset{\cong}{\longto} \Gamma\big(\BdiffE\big):
			&(\delta,b) &\longmapsto \big(\!\big(\delta,\iB(b)\big)\!\big),\\ 
			\Gamma\big(\BdiffE\big) &\overset{\cong}{\longto} \Gamma\big(\BdiffI\big):
			&\llparen\delta,l\rrparen &\longmapsto \big(\delta-\barnab_{\prA(l)},\prB(l)\big).  
		\end{aligned}
	\end{equation}
Hence, we obtain the following transfer of Bott connections.

\begin{prop}\label{prop:DBEI}
	Fix an embedding $\iB\colon B\to L$,
	the identifications in~\eqref{eq:BdiffI=BdiffE} transfer flat Bott $A$-connection~\eqref{eq:BottDBE} on $\BdiffE$ to a flat Bott $A$-connection on $\BdiffI$:
	\begin{align}\label{eq:BottDBEI}
		\nabla_{a}(\delta,b)=
		\big([\barnab_a,\delta]+\barnab_{\eth_{_b} a}, D_a b\big),
		\qquad \forall\, a\in \Gamma(A),\,\,(\delta,b)\in\Gamma\big(\BdiffI\big),
	\end{align}
	where $\eth_{_b} a$ is defined in~\eqref{eq:eth}  and $D_a b$ is the flat Bott $A$-connection~\eqref{eq:Bott} on $B$.
\end{prop}

\begin{thm}\label{thm:DBEI}
	After the transfer of Bott $A$-representations, we have that
        \begin{enumerate}[a)]
        \item the identification specified in~\eqref{eq:BdiffI=BdiffE} $$\BdiffI \cong \BdiffE$$
        automatically becomes an $A$-representation isomorphism,
        \item the upper short exact sequence appeared in~\eqref{Dia:BdiffE-pullback}:
        \begin{equation}\label{seq:AtiyahDBEI}
        	\xymatrix@C=2pc@R=2pc{
        		0\ar[r]&
        		\mathend(E)\ar[r]^{\tilde{\iota}_{_B}}& 
        		\BdiffI \ar[r]^{\quad\widetilde{\sigma}_{_B}}&
        		B\ar[r]&
        		0,
        	}
        \end{equation}
        respects the transferred $A$-representation structures and fits in the commutative diagram of $A$-representation sequences:
	\begin{equation*}
		\xymatrix@C=2pc@R=2pc{
			0\ar[r]&
			\mathend(E)\ar[r]^{\tilde{\iota}_{_B}}\ar@{=}[d]& 
			\BdiffI \ar[r]^{\quad\widetilde{\sigma}_{_B}}\ar[d]^{\cong}&
			B\ar[r]\ar@{=}[d]&
			0\,\\
			0 \ar[r]
			&\End(E)\ar[r]^{\iota_{_{B}}}&
			\BdiffE\ar[r]^{\qquad\sigma_{_B}}
			&B\ar[r]&0.
		}
	\end{equation*}
        
        \end{enumerate}
\end{thm}

\begin{defn}
	We call Sequence~\eqref{seq:AtiyahDBEI} the \textbf{$B$-Atiyah sequence of $E$ with respect to the embedding $\iB\colon B\to L$}.
	Because of the inclusion ${f_{_{\iB}}:\BdiffI\longhookrightarrow \LdiffE}$,
	we call $\BdiffI$ the \textbf{embedded $B$-Atiyah extension/bundle associated with $(L,A,E)$}.
\end{defn}

\begin{rmk}\label{rmk:indepiB}
	Another choice of embedding $i'_{_B}:B\to L$
	induces an $A$-representation isomorphism:
	$\Gamma\big(\BdiffI\big)\overset{\cong}{\longto}\Gamma\big(\fD_{i'_{_B}}(E)\big):
	\,(\delta,b)\mapsto (\delta+\barnab_{I(b)},b)$
	where $I=i'_{_B}-\iB$ is the difference map.
\end{rmk}

\vskip 10pt

\subsection{Complementary $B$-connections and extending $L$-connections}

We will see that the essential data for extending $A$-connections to $L$-connections is a relative/complementary type of connections furnished by the vector-bundle right splittings of $B$-Atiyah sequence~\eqref{seq:AtiyahDBEI}.

\vskip 5pt
\subsubsection{$B$-connection and its equivalents}

\begin{defn}\label{def:B-Conn}
	A \textbf{(complementary) $B$-connection $(\nabla^\circ,\iB)$ on $E$} consists of
	\begin{enumerate}
		\item a vector-bundle right splitting $\iB\colon B\to L$, and
		\item a bilinear map 
		$\nabla^\circ\colon\Gamma(B)\times\Gamma(E)\longto\Gamma(E),\quad (b,e)\longmapsto \nabla^\circ_{_b} e,$
	\end{enumerate}
	satisfying $\nabla^\circ_{_{fb}}e = f\nabla^\circ_{_b} e$ and $\nabla^\circ_{_b}(fe) = f\nabla^\circ_{_b}e+\big(\rhoB(b)f\big)e$,
	where $\rho_{_B}(b)=\rho_{_L}\circ\iB(b)\in\Gamma(TM)$, for all $b\in\bsection{B}$, $e\in\bsection{E}$, and $f\in C^{\infty} (M)$.
\end{defn}

\begin{rmk}\label{rmk:B-Conn}
	The defining property of $\nabla^\circ$ means that,
	the covariant derivative operator $\nabla^\circ_{_b}$ has a $TM$-symbol $\rho_{_B}(b)\in\Gamma(TM)$,
	or equivalently has an $L$-symbol $i_{_B}(b)\in\Gamma(L)$.
	Hence, we see that $\big(\nabla^\circ_{_b},\iB(b)\big)\in\Gamma\big(\LdiffE\big)$ by~\eqref{eq:ΓDLE}
	and $\big(\nabla^\circ_{_b},b\big)\in\Gamma\big(\BdiffI\big)$ by~\eqref{eq:ΓDBEI}.
\end{rmk}

\begin{rmk}
	The justification of the adjective ``complementary'' will be clear in the next Subsubsection.
\end{rmk}

Similar to Theorem~\ref{thm:conn=splitting}, a $B$-connection amounts to a vector-bundle splitting of the $B$-Atiyah sequence.
\begin{thm}\label{thm:Bconn=splitting}
	Fix an embedding $\iB\colon B\to L$, the following are equivalent:
	\begin{enumerate}[(1)]
		\item a $B$-connection $\nabla^\circ$ on $E$ with respect to $\iB$;
		\item a vector-bundle right splitting
		$\tilde{s}_{_B}\colon B\to \BdiffI$ of $B$-Atiyah sequence~\eqref{seq:AtiyahDBEI}, 
		or equivalently, a vector-bundle left splitting $\tilde{\theta}_{_B}\colon\BdiffI \longto \End(E)$ 
		with the mutual determination $\tilde{\iota}_{_B}\circ\tilde{\theta}_{_B}+\tilde{s}_{_B}\circ\widetilde{\sigma}_{_B}=id_{\BdiffI}$:
		\begin{equation}\label{seq:BdiffE}
			\xymatrix@C=2pc@R=2pc{
				0\ar[r]&  \End(E) \ar@<0.5ex>[r]^{\,\tilde{\iota}_{_B}}&
				{\BdiffI} \ar@<0.5ex>[r]^(0.65){\widetilde{\sigma}_{_B}} \ar@<0.5ex>@{-->}[l]^{\,\tilde{\theta}_{_B}}&
				B\ar[r]\ar@<0.5ex>@{-->}[l]^(0.35){\tilde{s}_{_B}}& 0;
			}
		\end{equation}
		\item \label{item:BconnDecomp} vector-bundle decomposition $\BdiffI\cong \mathend(E)\oplus B$ under which we have $\tilde{\iota}_{_{B}}\big(\mathend(E)\big)=\mathend(E)\oplus 0$ and $\widetilde{\sigma}_{_B}(0\oplus L)=L$;
		\item a vector-bundle map $\bar{s}_{_B} \colon B\to \diffE$ that makes the following lower triangle of vector bundles commute:
		\begin{equation*}
			\begin{split}
				\xymatrix@C=2.5pc@R=2pc{
					\BdiffI
					\ar@<0.5ex>[r]^{\widetilde{\sigma}_{_B}}
					\ar[d]_{f_{_B}}
					\ar@{}[dr]|(0.25){\text{\pigpenfont J}}|(0.7){\scalebox{1.2}{$\circlearrowleft$}}
					& B\, \ar[d]^{\rhoB}
					\ar@{-->}@<0.5ex>[l]^{\tilde{s}_{_B}} \ar@{-->}[ld]|(0.5){\scalebox{1}{$\bar{s}_{_B}$}}
					\\
					\diffE 
					\ar[r]^{\sigma}
					&TM.}
			\end{split}
		\end{equation*}
		(The correspondence between $\tilde{s}_{_B}$ and $\bar{s}_{_B}$ is that $\bar{s}_{_B}=f_{_B}\circ \tilde{s}_{_B}$, where both the pullback diagram and the bundle map $f_{_B}=f_{_L}\circ f_{_{\iB}}\colon \BdiffI\to\LdiffE\to\diffE$ are composed from~\eqref{Dia:LdiffE-pullback}~\eqref{Dia:BdiffE-pullback}.)
	\end{enumerate}
\end{thm}

Corollary~\ref{cor:conn=splitting} for $L$-connections also has an analogue for $B$-connections.
Its part~(a) is immediate as below due to the above statement (2), while its part~(b) will be deferred to the end of this Subsection.

\begin{cor} \label{cor:Bconn=splitting}
	$B$-connections $\nabla^\circ$ on $E$ with respect to a fixed embedding $\iB$ always exist,
	and they form an affine space modelled on $\Gamma\big(\Hom\big(B,\End(E)\big)\big)$.
\end{cor}

\vskip 5pt
\subsubsection{Relating $B$-connections and $L$-connections}
In reality, we can link every complementary $B$-connection to an extending $L$-connection.

Combining morphisms~\eqref{Dia:DALE-pullback},~\eqref{Dia:BdiffE-pullback} among Atiyah sequences with respect to $A,\,L,\,B$ induced from the embeddings $\iA,\,\iB$ and the projection $\prB$, we have the diagram:
\begin{equation}\label{Dia:DABLE-pullback}
	\begin{split}
		\xymatrix@C=2.5pc@R=2pc{
			0 \ar[r]
			&\mathend(E) \ar[r]^{{\iota}}\ar@{=}[d]
			&\mathfrak{D}_A(E) \ar@<0.5ex>[r]^{\quad\sigma_{_A}}\ar[d]_{f_{A,L}}
			&A\ar[r]\ar@<0.5ex>[d]^{\iA} \ar@<0.5ex>[l]^{\quad s_{_A}}
			&0\,\,\\
			0 \ar[r]
			&\mathend(E) \ar[r]^{{\iota}}
			&\LdiffE \ar@<0.5ex>[r]^{\quad\sigma_{_L}} \ar@<-0.5ex>[d]_{\widetilde{P}_{_B}}
			&L \ar[r] \ar@<0.5ex>@{-->}[u]^{{\rm pr}_{_A}} 
			\ar@<-0.5ex>[d]_{{\rm pr}_{_B}}
			\ar@{-->}@<0.5ex>[l]^{\quad s_{_L}}
			&0\,\,\\
			0 \ar[r]
			&\mathend(E) \ar[r]^{\tilde{\iota}_{_B}}\ar@{=}[u]
			& \BdiffI \ar@<0.5ex>[r]^{\quad\widetilde{\sigma}_{_B}}
			\ar@<-0.5ex>[u]_{f_{_{\iB}}}
			&B \ar[r]\ar@<-0.5ex>[u]_{\iB} \ar@{-->}@<0.5ex>[l]^{\quad \tilde{s}_{_B}}
			&0,
		}
	\end{split}
\end{equation}
where $f_{_{A,L}},\,f_{_{\iB}},\,\widetilde{P}_{_B}$ are given by~\eqref{eq:DAEtoDLA},\,\eqref{eq:symbol-DBE} and Lemma~\eqref{lem:DBEI}.

\begin{lem}\label{lem:Bconn=Lconn}
	Fix an embedding $\iB\colon B\to L$, 
	there is a one-to-one correspondence between the $B$-connections $\nabla^\circ$ with respect to $\iB$
	and the $L$-connections ${\nabla}$ extending $\barnab$.
\end{lem}
\begin{proof}
	In Diagram~\eqref{Dia:DABLE-pullback},
	the upper sequence already has a splitting related to $\barnab$.
	Hence, a splitting of the middle sequence is equivalent to a splitting of the lower sequence.
	Explicitly, given an $L$-connection ${\nabla}$ extending $\barnab$, one considers the $B$-connection 
	$\nabla^\circ_{_b}e:=\nabla_{_{\iB(b)}}e$.
	On the other hand, given a $B$-connection $\nabla^\circ$ with respect to $\iB$ (and also $\prA$),
	one considers the $L$-connection ${\nabla}_{_l}e:=\barnab_{_{\prA(l)}}e+{}^B\nabla_{_{\prB(l)}}e$.
\end{proof}

From above Lemma~\ref{lem:Bconn=Lconn} and Corollary~\ref {cor:Bconn=splitting}, we see that

\begin{cor}
	$L$-connections ${\nabla}$ extending the $A$-connection $\barnab$ on $E$ form an affine space modelled on $\Gamma\big(\Hom\big(B,\End(E)\big)\big)$.
\end{cor}

We can also unfix the embeddings $\iB$.

\begin{prop}\label{prop:Bconn=Lconn}
	The following data on a Lie triad $(L,A,E)$ are equivalent:
	\begin{enumerate}
		\item a $B$-connection $(\nabla^\circ,\iB)$ on $E$;
		\item an embedding $\iB\colon B\to L$ and an $L$-connection ${\nabla}$ extending the $A$-connection $\barnab$ on $E$;
		\item a vector-bundle right splitting $S_{_B}\colon B\to \LdiffE$ of quotient sequence~\eqref{seq:LiePairDLAE} for $\big(\LdiffE,\AdiffE\big)$:
		\begin{equation*}
			\xymatrix@C=2.5pc@R=2pc{
				0 \ar[r]&\fD_A(E)\ar[r]^{ f_{_{A,L}}}&\mathfrak{D}_L(E)\ar@<0.5ex>[r]^{\quad{\rm pr}_{_B}\circ \sigma_{_L}}
				&B\ar[r] \ar@{-->}@<0.5ex>[l]^{\quad S_{_B}}&0.
			}
		\end{equation*}
	\end{enumerate}
	Such a splitting $S_{_B}$   always exists and provides vector-bundle decompositions 
	$$\LdiffE\cong \AdiffE\oplus B \cong \End(E)\oplus A\oplus B.$$	 
\end{prop}

\begin{proof}
	The equivalence between (1) and~(2) is just Lemma~\ref{lem:Bconn=Lconn}.
	The equivalence between (1) and~(3) can be traced within Diagram~\eqref{Dia:LA-Atiyah3}:
	a splitting of its middle sequence involving $\LdiffE$ is equivalent to the compatible splittings of its upper and lower ones involving $L$ and $\BdiffE$, which correspond to $\iB$ and $\nabla^\circ$.
	
	For the last vector-bundle decomposition of $\LdiffE$, one uses the semi-direct product $\AdiffE\cong \End(E)\rtimes A$ in Proposition~\ref{prop:semiDirect}.
\end{proof}


\begin{cor}
The set of	$B$-connections $(\nabla^\circ,\iB)$ on $E$ is an affine space modelled on $\Gamma\big(\Hom\big(B,\AdiffE\big)\big)$
$\cong \Gamma\big(\Hom\big(B,\End(E)\big)\big)\oplus \Gamma\big(\Hom(B,A)\big)$.
\end{cor}

\vskip 10pt

\subsection{Third construction of Atiyah extensions as split bundles twisted by Atiyah cocycles}\label{subsec:splitAtiyahB}
Further equipping the embedded Atiyah extension $\BdiffI$ with an $L$-connection $\nabla$ on $E$ that extends $\barnab$,
we will construct a finer Atiyah extension $\AtExtE$ which is twisted by the Atiyah cocycle and turns out to be independent of the choice of $\iB$.

In Theorem~\ref{thm:Bconn=splitting}.\eqref{item:BconnDecomp}, we get a decomposition ${\BdiffI\cong \mathend(E)\oplus B}$ 
resulted from a $B$-connection $(\nabla^\circ,\iB)$ or equivalently from an extending $L$-connection $\nabla$.
We can compose that identification with the one ${\BdiffE \cong \BdiffI}$ in~\eqref{eq:BdiffI=BdiffE},
and henceforth obtain the identification $\BdiffE \cong \mathend(E)\oplus B$ evaluated on sections (cf.~\eqref{eq:LdiffE=EndE+L}):
\begin{equation}\label{eq:BdiffE=EndE+B}
	\begin{aligned}
		\Gamma\big(\BdiffE\big) &\overset{\cong}{\longto} \Gamma\big(\End(E)\oplus B\big):
		&\llparen\delta,l\rrparen &\longmapsto \big(\delta-\nabla_l\big)\oplus \prB(l),\\ 
		\Gamma\big(\End(E)\oplus B\big) &\overset{\cong}{\longto} \Gamma\big(\BdiffE\big):
		&\varphi \oplus b &\longmapsto \big(\!\big(\varphi+\nabla^\circ_b,\iB(b)\big)\!\big).  
	\end{aligned}
\end{equation}

\begin{prop}\label{prop:BdiffI=AtExtE}
	Fix a $B$-connection $(\nabla^\circ,\iB)$ on $E$, hence also an $L$-connection $\nabla$ extending $\barnab$ on $E$. 
	The above identifications transfer Bott $A$-representations~\eqref{eq:BottDBE},~\eqref{eq:BottDBEI}
	on $\BdiffE$ and $\BdiffI$ to a Bott $A$-representation on $\End(E)\oplus B$:
	\begin{align}\label{eq:BottEndE+B}
		D_{a}(\varphi\oplus b)=
		\Big(\barnab^{\End(E)}_a \varphi + R^{\nabla}_{A\otimes B}(a\otimes b)\Big) \oplus D_a b,
		\qquad \forall\, a\in \Gamma(A),\,\,\varphi\oplus b\in\Gamma\big(\End(E)\oplus B\big),
	\end{align}
	where $\barnab^{\End(E)}_a\varphi=[\barnab_a,\varphi]$ is $A$-representation~\eqref{eq:BottEndE} on $\End(E)$,
	and $D_a b=\prB\big([a,\iB(b)]\big)$ is Bott $A$-representation~\eqref{eq:Bott} on $B$,
	and $R^{\nabla}_{A\otimes B}(a\otimes b)=\nabla_a\nabla_{\iB(b)}-\nabla_{\iB(b)}\nabla_a-\nabla_{[a,\iB(b)]}$ is Atiyah cocycle~\eqref{eq:AtiyahCocycle} relative to $(L,A,\nabla)$.
\end{prop}

\begin{proof}
	To get a transferred connection on $\End(E)\oplus B$,
	we note that, for a section $\varphi\oplus b$ of $\End(E)\oplus B$,
	Bott $A$-connection~\eqref{eq:BottDBE} of $a\in\Gamma(A)$
	on the section $\big(\!\big(\varphi+\nabla^\circ_b,\iB(b)\big)\!\big)$ of $\BdiffE$ identified by~\eqref{eq:BdiffE=EndE+B} gives
	\begin{align*}
		D^{^{\BdiffE}}_{a}\big(\!\big(\varphi+\nabla^\circ_b,\iB(b)\big)\!\big)
		&=\big(\!\big([\barnab_a,\varphi]+[\barnab_a,\nabla^\circ_b],[a,\iB(b)]\big)\!\big).
	\end{align*}
	The latter can be identified via~\eqref{eq:BdiffE=EndE+B} back to be a section of $\End(E)\oplus B$:
	\begin{align*}
		\big([\barnab_a,\varphi]+[\barnab_a,\nabla^\circ_b]-\nabla_{[a,\iB(b)] }\big) \oplus \prB\big([a,\iB(b)]\big)
		=\Big(\barnab^{\End(E)}_a\varphi+R^{\nabla}_{A \otimes B}(a,b) \Big) \oplus D_a b.
	\end{align*} 
\end{proof}

\begin{thm}\label{thm:BdiffI=AtExtE}
	After the transfer of Bott $A$-representations, we denote $\End(E)\oplus B$ by $\AtExtE$, then
    \begin{enumerate}[(a)]
    \item the identification specified in~\eqref{eq:BdiffE=EndE+B}:
    $$\BdiffE \cong \AtExtE$$
    automatically becomes an $A$-representation isomorphism;
    \item $B$-Atiyah sequence~\eqref{seq:AtiyahDBE} is replaced by a new sequence of $A$-representations:
    \begin{equation}\label{seq:AtiyahEndE+B}
    	\xymatrix@C=2pc@R=2pc{
    		0\ar[r]&
    		\mathend(E)\ar[r]^{\hat{\iota}_{_B}\quad}& 
    		\End(E) \nabplus B \ar[r]^{\qquad\hat{\sigma}_{_B}}&
    		B\ar[r]&0,
    	}
    \end{equation}
     that fits in a commutative diagram of $A$-representation sequences:
	\begin{equation*}
		\xymatrix@C=2pc@R=2pc{
			0\ar[r]&
			\mathend(E)\ar[r]^{\hat{\iota}_{_B}\quad}& 
			\End(E) \nabplus B \ar[r]^{\qquad\hat{\sigma}_{_B}}&
			B\ar[r]&0\,\\
			0 \ar[r]&
			\End(E)\ar[r]^{\iota_{_{B}}}\ar@{=}[u]&
			\BdiffE\ar[r]^{\qquad\sigma_{_B}}\ar[u]^{\cong}&
			B\ar[r]\ar@{=}[u]&0,
		}
	\end{equation*}
    where $\hat{\iota}_{_B}(\varphi)=\varphi\oplus 0$ and $\hat{\sigma}_{_B}(\varphi\oplus b)=b$ for $\varphi\in\Gamma\big(\End(E)\big),\, b\in \Gamma(B)$.
    \end{enumerate}
\end{thm}

\begin{defn}\label{def:AtExtE}
	We call Sequence~\eqref{seq:AtiyahEndE+B} the \textbf{virtually split $B$-Atiyah sequence of $E$}
	and call $\End(E)\nabplus B$ the \textbf{virtually split $B$-Atiyah extension/bundle associated with $(L,A,E)$ and $\nabla$}.
\end{defn}

\begin{rmk}\label{rmk:BdiffI=AtExtE1}
	Each identification $\BdiffE\cong \End(E)\oplus B$ essentially comes from an identification $\LdiffE\cong \End(E)\oplus L$.
	In fact, the Lie bracket (Theorem~\ref{thm:EndE+L}.(b)) on $\Gamma\big(\End(E)\oplus L\big)$ equips the quotient bundle $\big(\End(E)\oplus L\big)/A=\End(E)\oplus B$ with a Bott $A$-connection which is exactly~\eqref{eq:BottEndE+B}.
\end{rmk}

\begin{rmk}\label{rmk:BdiffI=AtExtE2}
	Among the three components of~\eqref{eq:BottEndE+B},
	the two flat connections $\barnab^{^{\End(E)}}_a\varphi$ and $D_a b$ in their definitions~\eqref{eq:BottEndE},~\eqref{eq:Bott} are independent of both the embedding $\iB$ and the $L$-connection $\nabla$,
	meanwhile the Atiyah cocycle $R^{\nabla}_{A \otimes B}$
	is independent of $\iB$ but is dependent of $\nabla$, see~\eqref{eq:AtiyahCocycle}.
	This is why the notation $\AtExtE$ is related to $\nabla$ but free of $\iB$.
\end{rmk}

\begin{rmk}[cf. Remark~\ref{rmk:indepNabla}]\label{rmk:indepNabla1}
	The difference $(\nabla'-\nabla)|_B:B\to \End(E)$ of two $L$-connections $\nabla,\,\nabla'$ extending $\barnab$,
	results in a difference of Atiyah cocycles by a coboundary (see Theorem~\ref{thm:CSX}.\eqref{item:AtiyahIndepOfConn}).
	Such a difference will affect flat $A$-connection~\eqref{eq:BottEndE+B} on $\End(E)\oplus B$,
	but it can be reconciled by an $A$-representation isomorphism:
	\begin{align*}\label{eq:TwoEndE+B}
		I_{\nabla,\nabla'}: \Gamma\big(\mathend(E)\nabplus B\big) \longrightarrow \Gamma\big(\mathend(E)\scalebox{0.8}{$\stackrel{\nabla'}{\oplus}$} B\big),\quad \varphi\oplus b\longmapsto\big(\varphi+\nabla_{\iB(b)}-\nabla'_{\iB(b)}\big)\oplus b.
	\end{align*}
	Indeed, we have $A$-representation isomorphisms using~\eqref{eq:BdiffE=EndE+B}:
	$$\AtExtE \cong \BdiffE\cong \End(E)\scalebox{0.8}{$\stackrel{\nabla'}{\oplus}$} B.$$
\end{rmk}

\vskip 10pt
\subsection{Compatible connections and Atiyah extensions}\label{subsec:compatibleConn}

In view of the correspondence (Lemma~\ref{lem:Bconn=Lconn}) between extending $L$-connections and complementary $B$-connections,
the $A$-compatibility of the former converts into that of the latter,
which we furthermore interpret in terms of certain compatible splittings of the $B$-Atiyah sequence.

\begin{Def}\label{def:compatibleBConn}
	We say that a $B$-connection $(\nabla^\circ,\iB)$ is \textbf{$A$-compatible}  if its corresponding $L$-connection $\nabla$ satisfies
	\[\nabla_a\nabla_{\iB(b)}-\nabla_{\iB(b)}\nabla_a=\nabla_{[a,\iB(b)]},\quad \forall a\in\bsection{A}, b\in\bsection{B}, \]
	or equivalently, Atiyah cocycle~\eqref{eq:AtiyahCocycle} $R^{\nabla}_{A\otimes B}(a\otimes b)=\nabla_a\nabla_{\iB(b)}-\nabla_{\iB(b)}\nabla_a-\nabla_{[a,\iB(b)]}$  is zero.
\end{Def}

There is an $A$-compatible version of Theorem~\ref{thm:Bconn=splitting} on the correspondence between connections and splittings of Atiyah sequence.

\begin{prop}\label{prop:AtiyahCocyle=0}
	Given a $B$-connection $(\nabla^\circ,\iB)$,
	by Theorem~\ref{thm:Bconn=splitting} and Proposition~\ref{prop:Bconn=Lconn},
	let $\tilde{s}_{_B}\colon B\to \BdiffI$ be the corresponding vector-bundle right splitting in~\eqref{seq:BdiffE}
	and $\nabla$ the corresponding $L$-connection $\nabla$ extending $\barnab$.
	The following are equivalent:
	\begin{enumerate}
		\item the $L$-connection $\nabla$ is $A$-compatible, i.e., $R^{\nabla}_{A \otimes L}=0$ according to Definition~\ref{def:compatibleConn};
		\item the $B$-connection $(\nabla^\circ,\iB)$ is $A$-compatible, i.e., Atiyah cocycle $R^{\nabla}_{A \otimes B}=0$ according to Definition~\ref{def:compatibleBConn};
		\item the vector-bundle embedding $\hat{s}_{_B}\colon B\to \AtExtE\colon b\mapsto 0\oplus b$ that canonically splits Atiyah sequence~\eqref{seq:AtiyahEndE+B} is an $A$-representation morphism, or equivalently, the canonical vector-bundle projection $\hat{\iota}_{_B}\colon \AtExtE \to \End(E) \colon \varphi\oplus b\mapsto \varphi$ is an $A$-representation morphism;
		\item the vector-bundle right splitting $\tilde{s}_{_B}\colon B\to \BdiffI$ in Atiyah sequence~\eqref{seq:BdiffE} is an $A$-representation morphism, or equivalently, the associated left splitting $\tilde{\theta}_{_B}\colon\BdiffI \longto \End(E)$ is an $A$-representation morphism.
	\end{enumerate}
\end{prop}

\begin{proof}
	The equivalence $(1) \Leftrightarrow (2)$ is because of the flatness of the $A$-connection $\barnab$. For the equivalence {$(2) \Leftrightarrow (3)$}, we use~\eqref{eq:BottEndE+B} to see that $\nabla_{a}(0\oplus b)=R^{\nabla}_{A\otimes B}(a\otimes b) \oplus \nabla_a b$.
	So, $\hat{s}_{_B}$ respects the flat $A$-connections if and only if $R^{\nabla}_{A\otimes B}=0$. 
	The equivalence $(3) \Leftrightarrow (4)$ holds due to Theorem~\ref{thm:BdiffI=AtExtE}.
\end{proof}

The theorem below is a more general version of Proposition~\ref{prop:AtiyahCocyle=0} without fixing a $B$-connection $(\nabla^\circ,\iB)$ before hand.

\begin{thm}\label{thm:AtClass}
	Given a Lie triad $(L,A,E)$, then
	\begin{enumerate}[(a)]
		\item the Atiyah class $\alpha(L,A,E)\in H^1_{\CE}\big(A,B^\vee\otimes\mathend(E)\big)\cong \Ext_A\big(B,\End(E)\big)$ is represented by any one of $B$-Atiyah sequences~\eqref{seq:AtiyahDBE},\,~\eqref{seq:AtiyahDBEI}, and~\eqref{seq:AtiyahEndE+B} of $A$-representation extensions;
		\item the following are equivalent:
		\begin{enumerate}[(i)]
			\item the Atiyah class $\alpha(L,A,E)$ vanishes,
			\item there is an $A$-compatible $L$-connection on $E$ extending $\barnab$,
			\item there is an $A$-compatible complementary $B$-connection on $E$,
			\item there exists an $A$-representation right splitting $s_{_B}$ to $B$-Atiyah sequence~\eqref{seq:AtiyahDBE} for $\BdiffE$:
			\begin{equation*}\label{seq:AtiyahDBESplit}
				\xymatrix@C=2pc@R=2pc{
					0\ar[r]&\End(E)\ar[r]^{\iota_{_{B}}}
					&\BdiffE\ar@<0.5ex>[r]^(0.65){\sigma_{_B}}
					&B\ar[r]\ar@<0.5ex>@{-->}[l]^(0.35){s_{_B}}
					&0.
				}
			\end{equation*}
		\end{enumerate}
	\end{enumerate}

\end{thm}

\begin{proof}
	By standard homological algebra, the Atiyah cocycle $R^{\nabla}_{A \otimes B}$, appearing as a twisted term in $A$-representation(module) structure~\eqref{eq:BottEndE+B} of $\End(E)\nabplus B$,
	represents extension sequence~\eqref{seq:AtiyahEndE+B} for $\End(E)\nabplus B$.
	Hence, it also represents sequences~\eqref{seq:AtiyahDBE},~\eqref{seq:AtiyahDBEI} due to the transfers we have described.
	
	Moreover, we have the following equivalences:
	\begin{equation*}
		\begin{aligned}
			&\alpha(L,A,E)=0, &&\\
			\Leftrightarrow\quad &\exists\mbox{ an $L$-connection $\nabla$ extending $\barnab$ in an $A$-compatible way with } R^{\nabla}_{A\otimes B}=0, && \quad \mbox{(by Theorem~\ref{thm:CSX}.\eqref{item:AtiyahVanish}),}\\
			\Leftrightarrow\quad &\exists\mbox{ a $B$-connection $(\nabla^\circ,\iB)$ such that } R^{\nabla}_{A\otimes B}=0, && \quad \mbox{(by Proposition~\ref{prop:Bconn=Lconn}),}\\
			\Leftrightarrow\quad &
			\exists
			\left\{
			\begin{split}
				&\mbox{ a vector-bundle	 right splitting $\iB\colon B\longto L$,}\\ 
				&\mbox{ and an $A$-representation right splitting $\tilde{s}_{_B}\colon B\longto \BdiffI$},
			\end{split}\right.
			&& \quad \mbox{(by Proposition~\ref{prop:AtiyahCocyle=0}),}\\
			\Leftrightarrow\quad &\exists\mbox{ an $A$-representation right splitting ${s}_{_B}\colon B\longto \BdiffE$,}  && \quad \mbox{(by Theorem~\ref{thm:DBEI}).}\\
		\end{aligned}
	\end{equation*}
\end{proof}

In a nutshell, the Atiyah class obstructs the existence of $A$-compatible $L$-connections, and equivalently, also obstructs the $A$-representation splittings of $B$-Atiyah sequences.

\begin{cor}\label{cor:AtClass}
	If the Atiyah class $\alpha(L,A,E)$ vanishes, then the collection of $A$-compatible $L$-connections on~$E$ forms an affine space modelled on $\Gamma\big(\Hom\big(B,\End(E)\big)\big)^{\Gamma(A)}=H^0_{\CE}\big(A,B^\vee\otimes\End(E)\big)$.
\end{cor}

\vskip 25pt
\section{Matched connections, holomorphic and equivariant Atiyah algebroids}\label{Sec:matchedpairHolomAtiyahalgebroid}
\vskip 10pt

Our constructions of Atiyah extensions and their applications to compatible connections turn out to be particularly interesting
when $L$ is of the type $A\bowtie B$ matched between two Lie algebroids $A$ and $B$.
In this case, two important classes of examples stand out: holomorphic Lie algebroids and equivariant Lie algebroids.
We will use the theory developed in previous sections to describe holomorphic Lie algebroid connections, holomorphic Atiyah algebroids, and their equivariant counterparts.

\vskip 10pt
\subsection{Matched pairs and matched connections} \label{subsec:matchedPair}

\begin{Def}[\cite{Lu97}]
	Two Lie algebroids $A$ and $B$ over $M$ form a \textbf{matched pair}, 
	if there is a Lie algebroid structure on $A\oplus B$,
	such that $A\cong A\oplus 0$ and $B\cong 0 \oplus B$
	are Lie sub-algebroids of $A\oplus B$.
	The direct sum $A\oplus B$ equipped with such a Lie algebroid structure is called a \textbf{matched sum}, denoted by $A\bowtie B$.
\end{Def}

Given a matched pair $(A,B)$, there are two Lie pairs $(A\bowtie B, A)$ and $(A\bowtie B, B)$
which provide flat Bott connections $D$ of $A$ and $B$ mutually on each other.
The Jacobi, Leibniz and anchor conditions of $A\bowtie B$ are then equivalent to certain compatibility conditions for these Bott connections.

\begin{lem}[\cite{Lu97,Mokri97}]\label{lem:matchedPair}
	Two Lie algebroids $A$ and $B$ over $M$ form a matched pair
	if and only if there are a flat $A$-connection on $B$ and a flat $B$-connection on $A$, both denoted by~$D$, satisfying the relations for all $a,\,a'\in \Gamma(A),\,b,\,b'\in \Gamma(B)$ that	
	\begin{enumerate}[(i)]
		\item $D_a[b,b']_{_B}=[D_a b,b']_{_B}+[b,D_a b']_{_B}-D_{D_{b}a}b'+D_{D_{b'}a}b$,
		\item $D_b[a,a']_{_A}=[D_b a,a']_{_A}+[a,D_b a']_{_A}-D_{D_{a}b}a'+D_{D_{a'}b}a$,
		\item $\rhoB(D_{a}b)-\rhoA(D_{b}a)=[\rhoA(a),\rhoB(b)]_{_{TM}}$.
	\end{enumerate}
	If the above requirements are met, the Lie algebroid structure on the matched sum $A\bowtie B$ is specified as
	\begin{equation}\label{eq:matchedPair}
		\begin{gathered}
			\rho_{_{A\bowtie B}}(a+b)=\rhoA(a)+\rhoB(b),\\
			[a,a']_{_{A\bowtie B}}=[a,a']_{_A},  \quad[b,b']_{_{A\bowtie B}}=[b,b']_{_B},  \quad[a,b]_{_{A\bowtie B}}=-[b,a]_{_{A\bowtie B}}=D_{a}b-D_{b}a,
		\end{gathered}
	\end{equation}
	where we have identified $a$ and $b$ respectively with $a\oplus 0$ and $0\oplus b$ in $\Gamma(A\oplus B)$.
\end{lem}

The Lie triad notation $(L=A\bowtie B, A, E)$ means that we have a matched sum $A\bowtie B$ with a mutual flat Bott connection $D$ between $A$ and $B$, and there is a flat $A$-connection $\barnab$ on $E$.
Such a Lie triad will be called a \textbf{matched-type Lie triad}.

Note that we have a canonical inclusion $\iB\colon B\cong 0\oplus B\subseteq A\bowtie B$.
The relations among $B$-Atiyah algebroids, $B$-Atiyah extensions and $B$-Atiyah sequences are summarized in the following two theorems.

\begin{thm}\label{thm:matchedAtExt}
	Given a matched-type Lie triad $(A\bowtie B, A, E)$, we consider the embedded $B$-Atiyah extension $\BdiffI\subseteq \fD_{_{A\bowtie B}}(E)$,
	and the $B$-Atiyah algebroid $\fD_{_B}(E)$.
	We have that
	\begin{enumerate}[(a)]
		\item $\BdiffI=\fD_{_B}(E)$ canonically as a vector bundle, and Atiyah sequence~\eqref{seq:AtiyahDBEI}:
		\begin{equation*}
			\xymatrix@C=2pc@R=2pc{
				0\ar[r]&
				\mathend(E)\ar[r]^{\tilde{\iota}_{_B}}& 
				\BdiffI \ar[r]^{\quad\widetilde{\sigma}_{_B}}&
				B\ar[r]&
				0
			}
		\end{equation*}
		is simultaneously a sequence of $A$-representations and also a sequence of Lie algebroids;
		\item Atiyah algebroids $\BdiffE$ and $\fD_{_{A\bowtie B}}(E)$ fit
		in a commutative diagram of Lie algebroid sequences:
		\begin{equation*}
			\xymatrix@C=2pc@R=2pc{
				0 \ar[r]&\End(E)\ar[r] \ar@{=}[d]
				&\mathfrak{D}_{B}(E)\ar[r]\ar@{^{(}->}[d]^{f_{\iB}}&B\ar[r]\ar@{^{(}->}[d]^{\iB}&0\,\\
				0 \ar[r]&\End(E)\ar[r]
				&\mathfrak{D}_{A\bowtie B}(E)\ar[r]&A\bowtie B\ar[r]&0;\\
			}
		\end{equation*}
		\item the Atiyah algebroid $\fD_{_{A\bowtie B}}(E)$ is a matched sum:
		$$\fD_{_{A\bowtie B}}(E)\cong A\bowtie \fD_{_B}(E),$$
		where for $a\in \Gamma(A),\,\,(\delta,b)\in\Gamma\big(\fD_{_B}(E)\big)$, the flat Bott connection of $A$ on $\fD_{_B}(E)$ is given by~\eqref{eq:BottDBEI}: ${D_{a}(\delta,b)=
		\big([\barnab_a,\delta]+\barnab_{\nabla_b a}, D_a b\big)}$, and the flat Bott connection of $\fD_{_B}(E)$ on $A$ is given by $D_{(\delta,b)} a:=D_b a$.
	\end{enumerate}
\end{thm}

\begin{proof}
	Both the embedded Atiyah extension $\BdiffI$
	and the Atiyah algebroid $\fD_{_B}(E)$ 
	are the pullback of $\diffE$ via $\rhoL\circ \iB=\rhoB$.
	Thus, we see that $\BdiffI=\fD_{_B}(E)$ and the Atiyah sequences are the same for them, which confirms~(a).
	
	Since $\iB$ is a Lie algebroid embedding, by Proposition~\ref{prop:functoriality} of functoriality,
	it induces a morphism between Atiyah sequences as in~(b).
	In particular,we have a
	Lie algebroid morphism $f_{\iB}\colon \BdiffI=\fD_{_B}(E)\longto\fD_{_{A\bowtie B}}(E)$, which is a vector-bundle right splitting to Sequence~\eqref{seq:DLE-DBE2}: ${0\longto A \overset{s_{_{A,L}}}{\longto} \fD_{_{A\bowtie B}}(E) {\longto} \BdiffI\longto 0}$ in Lemma~\ref{lem:DBEI}.
	So, this right splitting respects the Lie algebroid structures and we have a matched sum $\fD_{_{A\bowtie B}}(E)=A\bowtie f_{\iB}\big(\BdiffI\big)$.
	
	To deduce the mutual flat Bott connections of $A$ and $\BdiffE$, we consider the bracket in $\Gamma\big(\LdiffE\big)$:
	\begin{align*}
		\big[s_{_{A,L}}(a),f_{\iB}(\delta,b)\big]=\big[(\barnab_a,a),(\delta,b)\big]
		&=\big([\barnab_a,\delta],[a,b]\big)\\
		&=\big([\barnab_a,\delta],D_a b-D_b a\big)=\big([\barnab_a,\delta]+\barnab_{D_b a},D_a b \big)-\big(\barnab_{D_b a},D_b a\big),
	\end{align*}
	where in the final term, the first component in $f_{\iB}\big(\BdiffI\big)$ gives $D_{a}(\delta,b)$ and the second component in $s_{_{A,L}}(A)$ gives $D_{(\delta,b)} a$.
\end{proof}

The embedded Atiyah extension $\BdiffI$ in the above theorem can be replaced by a virtually split Atiyah extension $\End(E)\nabplus B$.

\begin{thm} \label{thm:matchedAtExt1}
	Given a matched-type Lie triad $(L=A\bowtie B, A, E)$, we have that
	\begin{enumerate}[(a)]
		\item the virtually split $B$-Atiyah extension $\End(E)\nabplus B$ is also the virtually split $B$-Atiyah algebroid with respect to the $B$-connection $\nabla|_B$,
		and Atiyah sequence~\eqref{seq:AtiyahEndE+B}:
	    \begin{equation*}
			\xymatrix@C=2pc@R=2pc{
				0\ar[r]&
				\mathend(E)\ar[r]^{\hat{\iota}_{_B}\quad}& 
				\End(E) \nabplus B \ar[r]^{\qquad\hat{\sigma}_{_B}}&
				B\ar[r]&0,
			}
		\end{equation*}
		is simultaneously a sequence of $A$-representations and also a sequence of Lie algebroids;
		\item Atiyah algebroids $\End(E) \nabplus B$ and $\End(E) \nabplus L$ fit
		in a commutative diagram of Lie algebroid sequences:
		\begin{equation*}
			\xymatrix@C=2pc@R=2pc{
				0 \ar[r]&\End(E)\ar[r] \ar@{=}[d]
				&\End(E) \nabplus B\ar[r]\ar@{^{(}->}[d]^{f_{\iB}}&B\ar[r]\ar@{^{(}->}[d]^{\iB}&0\,\\
				0 \ar[r]&\End(E)\ar[r]
				&\End(E) \nabplus L\ar[r]& L\ar[r]&0;\\
			}
		\end{equation*}
		\item the virtually split $L$-Atiyah algebroid $\End(E)\nabplus L$ is a matched sum:
		$$\End(E)\nabplus L \cong A\bowtie \big(\End(E)\nabplus B\big),$$
		where for $a\in \Gamma(A),\,\,\varphi\oplus b\in\End(E)\nabplus B$, the flat Bott connection of $A$ on $\End(E)\nabplus B$ is given by~\eqref{eq:BottEndE+B}:
		$D_{a}(\varphi\oplus b)=\Big([\barnab_a,\varphi] + R^{\nabla}_{A\otimes B}(a\otimes b)\Big) \oplus D_a b$,
		and the Bott connection of $\End(E)\nabplus B$ on $A$ is given by $D_{(\varphi \oplus b)} a:=D_b a$.
	\end{enumerate} 
\end{thm}


Next, we turn to the $B$-connections.

\begin{defn}\label{Def:matchedconnection}
	Given a matched-type Lie triad $(L=A\bowtie B, A, E)$, we say that a $B$-connection $\nabla^\circ$ on $E$ is \textbf{$A$-matched} if it is $A$-compatible, i.e., it satisfies
	$$\barnab_a\nabla^\circ_b e-\nabla^\circ_b\barnab_a e=\nabla_{[a,b]}e,\qquad \forall\, a\in\Gamma(A),\, b\in\Gamma(B),\, e\in\Gamma(E),$$
	where $\nabla$ is the extending $L$-connection corresponding to $\nabla^\circ$.
\end{defn}

Using the structure equation of $[a,b]$ from~\eqref{eq:matchedPair}, we have $\nabla_{[a,b]}=\nabla^\circ_{D_a b}-\barnab_{D_b a}$.
The Atiyah cocycle $R^{\nabla}_{A\otimes B} \in \Omega_A^1\big({B^\vee}\otimes \mathend(E)\big)$ takes the expression
\[
R^{\nabla}_{A\otimes B}(a,b)e=\barnab_a\nabla^\circ_b e-\nabla^\circ_b\barnab_a e-(\nabla^\circ_{D_a b}e-\barnab_{D_b a}e), \qquad \forall\, a\in\Gamma(A),\, b\in\Gamma(B),\, e\in\Gamma(E),
\]
and represents the Atiyah class $\alpha(A\bowtie B,A,E)\in H^1_{\mathrm{CE}}\big(A,\,B^\vee\otimes\mathend(E)\big)$.

Proposition~\ref{prop:AtiyahCocyle=0} certainly holds in the matched-type case to build a one-to-one correspondence between matched Lie algebroid connections and compatible splittings of the $B$-Atiyah sequences.

In addition, we have a refined correspondence for the flat connections:

\begin{prop}\label{prop:matchedFlat}
	Given a matched-type Lie triad $(L=A\bowtie B,A,E)$ and a $B$-connection $\nabla^\circ$ with respect to the canonical $\iB$, we use the same notations as in Proposition~\ref{prop:AtiyahCocyle=0}. The following are equivalent:
	\begin{enumerate}
		\item the $L$-connection $\nabla$ is flat, i.e., $R^{\nabla}=0$;
		\item the $B$-connection $\nabla^\circ$ is $A$-compatible and flat, i.e., Atiyah cocycle $R^{\nabla}_{A \otimes B}=0$ and $R^{\nabla^\circ}=0$;
		\item the vector-bundle embedding $\hat{s}_{_B}\colon B\to \AtExtE\colon b\mapsto 0\oplus b$ is an $A$-representation morphism and also a Lie algebroid morphism;
		\item the vector-bundle splitting $\tilde{s}_{_B}\colon B\to \BdiffI$ in~\eqref{seq:BdiffE} is an $A$-representation morphism and also a Lie algebroid morphism.
	\end{enumerate}
\end{prop}

\begin{proof}
	Since $L=A\oplus B$ as vector bundles, the total curvature $R^{\nabla}$, which is an $\End(E)$-valued $2$-form of $L$, splits into $R^{\nabla}=R^{\nabla}|_{A\wedge A}+R^{\nabla}_{A \otimes B}+R^{\nabla}|_{B\wedge B}$, where the first term is zero by our assumption that $\nabla|_A=\barnab$ is flat and the third  term is $R^{\nabla^\circ}$.
	Thus, the vanishing of $R^{\nabla}$ is the same as that of $R^{\nabla^\circ}$ and the Atiyah cocycle $R^{\nabla}_{A \otimes B}$ together. This proves the equivalence $(1) \Leftrightarrow (2)$.
	
	The equivalence between (2) and (3)\&(4) follows from Proposition~\ref{prop:AtiyahCocyle=0}  after combining Clause~(b) of Corollary~\ref{cor:conn=splitting} on the relation between flat connections and Lie algebroid morphisms.
\end{proof}

Two important instances of  matched connections are demonstrated in the next two subsections.

\vskip 10pt
\subsection{Holomorphic Atiyah algebroids and holomorphic algebroid connections}\label{SubSec:holomAtiyahalgebroid}
In the following, we use the curly notation $\E$ for a holomorphic vector bundle and $E$ for its underlying smooth $\C$-vector bundle.
The smooth complex valued section space is denoted by $\Gamma_U(E)$ on open subsets $U\subseteq M$,
whereas the holomorphic section space is $\E(U)=\{s\in\Gamma_U(E) \mid \parbar s=0\}$.

\begin{defn}
	A \textbf{holomorphic Lie algebroid} $(\L,[~\cdot,\cdot~]_{_\L},\rho_{_\L})$ over a complex manifold $M$ consists of
	\begin{enumerate}[(i)]
		\item a holomorphic vector bundle $\L\to M$,
		\item a Lie bracket $[~\cdot,\cdot~]_{_\L}$ on each space $\L(U)$ of holomorphic $\L$-sections over an open subset $U\subseteq M$,
		making the functor $\L(\text{-})$ a sheaf of Lie algebras,
		\item a holomorphic vector bundle map $\rho_{_\L}$, called the anchor, from $\L$ to the holomorphic tangent bundle~$\mathscr{T}M$,
	\end{enumerate}
	such that $\rho_{_\L}$ induces a Lie algebra sheaf morphism from $\L(\text{-})$ to $ \mathscr{T}M(\text{-})$ satisfying the Leibniz rule
	$$[l_1,fl_2]_{_\L}=\big(\rho_{_\L}(l_1)f\big)l_2+f[l_1,l_2]_{_\L},$$
	for all holomorphic functions $f \in \O(U)$ and holomorphic sections $l_1,l_2 \in \L(U)$ over every open subset $U\subseteq M$. 
\end{defn}

Similarly, the definitions of Lie algebroid morphisms and Lie subalgebroids carries over to the holomorphic setting. 

In~\cite{LSX08}, Laurent-Gengoux, Sti\'{e}non and Xu showed that the holomorphic Lie algebroid structure on $\L$ can be extended to be a smooth Lie $\C$-algebroid structure on the  $\C$-vector bundle $L$ that underlies $\L$.
For instance, the holomorphic tangent bundle $\mathscr{T}M$ is a holomorphic Lie algebroid which extends to its underlying smooth $\C$-vector bundle $T^{1,0}M$ as a smooth Lie $\C$-algebroid.

Moreover, there is a matched pair between $L$ and the anti-holomorphic tangent Lie algebroid $T^{0,1}M$, viewed as a smooth Lie $\C$-algebroid.

\begin{thm}[{\cite[Thm.~4.8]{LSX08}}]\label{thm:hol-match}
	Let $L$ be a smooth Lie $\C$-algebroid with anchor $\rhoL$ valued in $T^{1,0}M\subseteq T_\C M$. 
	Then, every holomorphic Lie algebroid $\L$ that $L$ underlies
	is equivalent to a matched sum $T^{0,1}M\bowtie L$,
	where the Bott connection of $T^{0,1}M$ on $L$ is given by the Dolbeault $(0,1)$-operator: $\parbar_{X} l$ for $X\in\Gamma(T^{0,1}M),\, l\in\Gamma(L)$
	and the Bott connection of $L$ on $T^{0,1}M$ is given by the $(1,0)$-operator: $\partial_{\rhoL(l)}X$.
\end{thm}

Combining Theorems~\ref{thm:matchedAtExt} and~\ref{thm:hol-match}, we deduce the following consequence. 
\begin{prop}\label{prop:holAtiyah}
	Let $\L$ be a holomorphic Lie algebroid and $\E$ a holomorphic vector bundle on a complex manifold $M$,
	with their smooth $\C$-underlyings denoted by $L$ and $E$.
	We have
	\begin{enumerate}[(a)]
		\item $(T^{0,1}M\bowtie L,\,T^{0,1}M,\, E)$ is a Lie triad;
		\item the smooth Lie $\C$-algebroid $\LdiffE$ satisfies $T^{0,1}M\bowtie \LdiffE\cong\fD_{_{(T^{0,1}M)\bowtie L}}(E)$ and underlies a holomorphic Lie algebroid $\mathfrak{D}_\L(\E)$;
		\item for any open set $U\subseteq M$, the local holomorphic section space
		\[
		\big(\mathfrak{D}_\L(\E)\big)(U)=\big\{(\delta,l) \in \Gamma_U\big(\LdiffE\big) \mid \parbar_X (\delta,l) = 0,\, \forall X\in \Gamma_U(T^{0,1}M)\big\}
		\]
		consists of holomorphic covariant $\L$-differential operators $\delta$ on $\E(U)$ with holomorphic symbols ${l\in\L(U)}$.
	\end{enumerate}
\end{prop}
\begin{proof}
	 As in Example~\ref{Ex:triadOfHoloBundle},
	 the underlying smooth $\C$-vector bundle $E$ is a representation of $T^{0,1}M$ by the Dolbeault operator $\parbar$.
	 Since Theorem~\ref{thm:hol-match} gives us a matched sum $T^{0,1}M\bowtie L$,
	 we get in general a Lie triad $(T^{0,1}M\bowtie L,\, T^{0,1}M ,\, E)$, and in particular a Lie triad $(T_\C M = T^{0,1}M\bowtie T^{1,0}M,\, T^{0,1}M ,\, E)$.
	 
	 Next, applying Theorem~\ref{thm:matchedAtExt},
	 we have a bigger matched sum $T^{0,1}M\bowtie \LdiffE\cong\fD_{_{(T^{0,1}M)\bowtie L}}(E)$ where the anchor of $\LdiffE$ is pulled back from $\rhoL\colon L\to T^{1,0}M$.
	 In fact, due to functoriality~\eqref{diag:functoriality}, the Atiyah algebroid $\LdiffE$ fits in a commutative diagram of Lie algebroids:
	 \begin{equation*}
	 	\xymatrix@C=2.5pc@R=2pc{
	 			0 \ar[r]
	 			&\mathend(E) \ar[r]^{{\iota}}\ar@{=}[d]
	 			&\mathfrak{D}_L(E) \ar[r]^{\quad\sigma_{_L}}\ar[d]_{f_{_{\rhoL}}}
	 			\ar@{}[dr]|(0.25){\text{\pigpenfont J}}
	 			&L \ar[r]\ar[d]^{\rhoL}
	 			&0\,\,
	 			\\
	 			0 \ar[r]
	 			&\mathend(E) \ar[r]^{{\iota}}
	 			& \mathfrak{D}_{T^{1,0}M}(E) \ar[r]^{\quad\sigma_{_{T^{1,0}M}}}
	 			&T^{1,0}M\ar[r] 
	 			&0.
	 	}
	 \end{equation*}
	 All these Lie algebroids are matched with $T^{0,1}M$ and anchored in $T^{1,0}M$, 
	 (the explanation is similar to Example~\ref{ex:holomEnd(E)} on ${\End(E)\rtimes T^{0,1} M}$).
	 Thus, by Theorem~\ref{thm:hol-match}, they all underlie holomorphic Lie algebroids.
	 Moreover, the Lie algebroid morphisms among them also respect $T^{0,1}M$-representation structures by Clause~(a) of Theorem~\ref{thm:matchedAtExt}, 
	 hence are holomorphic Lie algebroid morphisms:
	 \begin{equation*}
	 	\xymatrix@C=2.5pc@R=2pc{
	 		0 \ar[r]
	 		&\mathend(\E) \ar[r]^{{\iota}}\ar@{=}[d]
	 		&\mathfrak{D}_{\L}(\E) \ar[r]^{\quad\sigma_{_\L}}\ar[d]_{f_{_{\rhoL}}}
	 		\ar@{}[dr]|(0.25){\text{\pigpenfont J}}
	 		&\L \ar[r]\ar[d]^{\rho_{_\L}}
	 		&0\,\,
	 		\\
	 		0 \ar[r]
	 		&\mathend(\E) \ar[r]^{{\iota}}
	 		& \mathfrak{D}_{\mathscr{T}M}(\E) \ar[r]^{\quad\sigma_{_{\mathscr{T}M}}}
	 		&\mathscr{T}M\ar[r] 
	 		&0.
	 	}
	 \end{equation*}
	 
	 As for the local holomorphic section space $\mathfrak{D}_\L(\E)(U)$,
	 since $\Gamma_U(T^{0,1}M)$ is $C^\infty(U,\C)$-spanned by anti-holomorphic sections,
	 the condition $\parbar_X (\delta,l) = 0$ suffices to hold for $X\in\Gamma_U(T^{0,1}M)$ satisfying $\partial X = 0$.
	 In this case, by Clause~(c) of Theorem~\ref{thm:matchedAtExt},
	 we have $\parbar_X (\delta,l)=\big([\parbar_X,\delta]+\parbar_{\partial_l X}, \parbar_X l\big)=\big([\parbar_X,\delta], \parbar_X l\big)$
	 whose vanishing means that
	 $$[\parbar_X,\delta]=0,\, \parbar_X l=0 \mbox{ for all } X\in\Gamma_U(T^{0,1}M) \mbox{ satisfying } \partial X = 0.$$
	 In other words, $\delta$ is a holomorphic covariant differential operator with holomorphic $L$-symbol $l$.
\end{proof}

Now, we fix a holomorphic Lie algebroid $\L$ and a holomorphic vector bundle $\E$ on a complex manifold $M$.

\begin{Def}\label{def:holLdiffE}
	We call the bundle of holomorphic covariant $\L$-differential operators, 
	$\mathfrak{D}_\L(\E)$, the \textbf{holomorphic $\L$-Atiyah algebroid of $\E$}, 
	whose structure maps are endowed from~\eqref{eq:LdiffE}
	and are evaluated for local holomorphic sections $(\delta,l),\,(\delta',l')\in \mathfrak{D}_\L(\E)(U)$ over any open subset $U\subseteq M$ to be:
	\begin{equation*}
		\big[(\delta,l),({\delta}',l')\big]_{\mathfrak{D}_\L(\E)(U)}=\big([\delta,{\delta}']_{\mathfrak{D}(\E)(U)},\,[l,l']_{_{\L}}\big),\qquad\quad
			\rho_{_{\mathfrak{D}_\L(\E)(U)}}(\delta,l)=\sigma(\delta)=\rho_{_{\L}}(l).
	\end{equation*}
	And we call sequence
	\begin{equation}\label{seq:holLdiffE}
		\xymatrix@C=2pc@R=2pc{
			0\ar[r]&  \mathend(\E) \ar[r]^{{\iota}}& {\mathfrak{D}_\L(\E)} \ar[r]^{\quad \sigma_{_\L}}&
			\L\ar[r]& 0,
		}
	\end{equation}
	the \textbf{holomorphic $\L$-Atiyah sequence of $\E$}.
\end{Def}

Following Atiyah's original treatment, we relate $\mathfrak{D}_\L(\E)$ with holomorphic connections, but now in  $\L$-symbols.

\begin{Def}
	A \textbf{local holomorphic $\L$-connection}  on $\E$, over an open neighbourhood $U\subseteq M$, is a bilinear map
	\[
	\nabla^{\prime,U} \colon \L(U)\times \E(U)\longrightarrow \E(U),\quad (l,e)\longmapsto \nabla^{\prime,U}_{_l} e, 
	\] 
	satisfying 
	$$\nabla^{\prime,U}_{_{fl}} e = f(\nabla^{\prime,U}_{_l} e), \quad \text{and} \quad \nabla^{\prime,U}_{_l}(fe) = f(\nabla^{\prime,U}_{_l} e) + \big(\rhoL(l)f\big)e,\qquad \forall \,e\in \E(U),\, l\in\L(U),\,f\in \O(U).$$
	The \textbf{curvature} of $\nabla^{\prime,U}$ is an $\O(U)$-linear map $R\colon \Lambda^2 \L(U)\to \big(\mathend(\E)\big)(U)$ defined by
	\[R(l_1,l_2)=\nabla^{\prime,U}_{l_1}\nabla^{\prime,U}_{l_2}-\nabla^{\prime,U}_{l_2}\nabla^{\prime,U}_{l_1}-\nabla^{\prime,U}_{[l_1,l_2]},\qquad \forall \,l_1,l_2\in\L(U). \]
\end{Def}

Locally at $U$, the holomorphic connection $\nabla^{\prime,U}$ associates to each holomorphic $L$-section $l$ a holomorphic covariant differential operator $\nabla^{\prime,U}_l\in \big(\mathfrak{D}_\L(\E)\big)(U)$.
This is the local holomorphic version of Theorem~\ref{thm:conn=splitting}:

\begin{prop}\label{prop:localHoloConn}
	A local holomorphic $\L$-connection $\nabla^{\prime,U}$ on $\E$ over
	an open neighbourhood $U\subseteq M$ is equivalent to a holomorphic vector-bundle right splitting
	$s\colon \L|_U\to \mathfrak{D}_\L(\E)|_U$ of $\L$-Atiyah sequence~\eqref{seq:holLdiffE} restricted on $U$.
	Moreover, $\nabla^{\prime,U}$ is flat if and only if the splitting $s$ is a holomorphic Lie algebroid morphism (making $\L|_U$ a holomorphic Lie subalgebroid of $\mathfrak{D}_\L(\E)|_U$).
\end{prop}

We could directly use a global splitting $s$ for the definition of a global connection.

\begin{Def}[cf. {\cite[p.~188]{Atiyah57}} ]\label{dfn:holo}
	A \textbf{global holomorphic $\L$-connection} on $\E$ is a holomorphic vector-bundle right splitting
	$s\colon \L\to \mathfrak{D}_\L(\E)$ of $\L$-Atiyah sequence~\eqref{seq:holLdiffE}:
	\begin{equation*}
		\xymatrix@C=2pc@R=2pc{
			0\ar[r] &  \mathend(\E) \ar[r]^{{\iota}} & {\mathfrak{D}_\L(\E)} \ar@<0.5ex>[r]^{\quad \sigma_{_\L}} &
			\L\ar[r] \ar@<0.5ex>@{-->}[l]^{\quad s} & 0.
		}
	\end{equation*}
\end{Def}

Equivalently, a global holomorphic $\L$-connection $\nabla^\prime$ can be defined by patching from local ones $\nabla^{\prime,i}$ over $U^i$ for an open covering $M=\cup_{i\in I}U_i$,
provided that compatibility condition $\nabla^{\prime,j}|_{U_i\cap U_j}=\nabla^{\prime,i}|_{U_i\cap U_j}$ holds.
If this is the case, the local curvatures $R^i$ also patch into a \textbf{global curvature} as a holomorphic vector bundle map $R\colon \wedge^2 \L\to \End(\E)$.

Over a contractible open neighbourhood $U_i\subseteq M$, a local holomorphic $\L$-connection on $\E$ always exists 
because one can take a holomorphic basis $e^i_1,\,e^i_2,\ldots,\,e^i_n$ of $\E$ over $U_i$, then declares $\nabla^{\prime,i} e^i_k = 0$ and extends $\nabla^{\prime,i}$ to $\E(U_i)$ by Leibniz rule.

However, its global presence over the whole $M$ is not guaranteed. 
Suppose $\mathscr{U}=\{U_i\}_{i\in I}$ is a good open cover of $M$,
where on each $U_i$ we pick a local holomorphic $\L$-connection $\nabla^{\prime,i}$, then
$$\psi(U_i\cap U_j):= \nabla^{\prime,j}|_{U_i\cap U_j}-\nabla^{\prime,i}|_{U_i\cap U_j} \in \big(\L^\vee \otimes End(\E)\big)(U_i\cap U_j)$$
defines a closed \v{C}ech $1$-cocycle up to exact \v{C}ech $1$-cocycle, and hence a \v{C}ech cohomology class $$\check{\alpha} \in \check{H}^1\big(\mathscr{U},\,\L^\vee \otimes End(\E)\big)$$
which represents holomorphic $\L$-Atiyah sequence~\eqref{seq:holLdiffE} as an extension of $\L$ by $\End(\E)$,
and obstructs the passage from the local holomorphic $\L$-connections $\{\nabla^{\prime,i}\}_{i\in I}$ to a global holomorphic $\L$-connection $\nabla^{\prime}$.

By the \v{C}ech-to-Dolbeault procedure,
we can choose a smooth partition of unity $\{\rho_i\}$ for the open cover $\mathscr{U}$, 
and glue the local holomorphic $\L$-connections $\{\nabla^{\prime,i}\}_{i\in I}$
into a global smooth $L$-connection
$$\nabla:=\sum \rho_i \cdot (\nabla^{\prime,i}+\parbar).$$
With respect to a local holomorphic basis of $\E$ over $U_i$,
we can write $\nabla|_{U_i}=d+\theta_i$ for $\theta_i\in \Gamma_{U_i}(L^\vee\otimes \End(E))$, with the curvature form
\[\Theta=d\theta_i+\frac{1}{2}[\theta_i,\theta_i]=\partial\theta_i+\frac{1}{2}[\theta_i,\theta_i]+\parbar\theta_i.\]

The last term assembles into $\{\parbar\theta_i\}_{i\in I}$
and behaves as a tensor under local holomorphic base changes of~$\E$.
Hence, it defines a global $\parbar$-closed element in $\Gamma\big((T^{0,1}M)^\vee \otimes L^\vee\otimes \End(E)\big)$ and consequently a Dolbeault cohomology class
\[\alpha_{_{\parbar}} \in H^1_{\parbar} \big(M,\,\L^\vee\otimes \End(\E)\big)\]
which obstructs the $(1,0)$-part $\nabla'=\nabla-\parbar$ of the global smooth $L$-connection $\nabla$ from being a global holomorphic $\L$-connection.

In fact, the local $(T^{0,1}M \otimes L)$-type curvature forms $\{\parbar\theta_i\}_{i\in I}$ correspond to the global $(T^{0,1}M \otimes L)$-type curvature operator:
\[R_{T^{0,1}M\otimes L}\colon T^{0,1}M\otimes L \longto \mathend(E), \quad X\otimes l\longmapsto \parbar_X\nabla_{l}-\nabla_{l}\parbar_X-\nabla_{[X,l]},\quad \forall\, X\in\bsection{T^{0,1}M},\, l\in\bsection{L}, \]
where $[X,l]=\parbar_X l - \partial_{\rho_{_L}(l)} X \in \Gamma(T^{0,1}M\bowtie L)$.
This $R_{T^{0,1}M\otimes L}$ is the Atiyah cocycle~\eqref{eq:AtiyahCocycle} and hence represents the Atiyah class
\[\alpha \in H^1_{\CE}\big(T^{0,1}M,\,L^\vee\otimes \End(E)\big).\]
Note that the Chevalley-Eilenberg differential for the Lie algebroid $T^{0,1}M$
is the Dolbeault operator $\parbar$.
Thus, the above Lie algebroid cohomology is exactly the Dolbeault cohomology.

In summary, we have three interpretations of the Atiyah class via different cohomology theories:

\begin{tcolorbox}[colback=white]
	\begin{equation}
		\begin{aligned}
			\check{H}^1\big(\mathscr{U},\,\L^\vee \otimes End(\E)\big) &&\cong&& H^1_{\parbar} \big(M,\,\L^\vee\otimes \End(\E)\big) &&=&& H^1_{\CE}\big(T^{0,1}M,\,L^\vee\otimes \End(E)\big):\\
			\check{\alpha}\hspace{1.5cm} &&\leftrightarrow&& \alpha_{_{\parbar}} \hspace{1.5cm}  &&\leftrightarrow&& \alpha.\hspace{2cm} 
		\end{aligned}
\end{equation}
\end{tcolorbox}

And we obtain the global holomorphic version of Proposition~\ref{prop:localHoloConn} and Theorem~\ref{thm:AtClass}:

\begin{thm}\label{thm:globalHoloConn}
	The existence of a global holomorphic $\L$-connection $\nabla'$ on $\E$, i.e., a holomorphic vector-bundle right splitting
	$s\colon \L\to \mathfrak{D}_\L(\E)$ of $\L$-Atiyah sequence~\eqref{seq:holLdiffE},
	is obstructed by the Atiyah class $\alpha \in H^1_{\CE}\big(T^{0,1}M,\,L^\vee\otimes \End(E)\big)$ associated with the Lie triad $(T^{0,1}M\bowtie L,\,T^{0,1}M,\, E)$.
	Moreover, $\nabla'$ is flat if and only if the splitting~$s$ is a holomorphic Lie algebroid morphism (making $\L$ a holomorphic Lie subalgebroid of $\mathfrak{D}_\L(\E)$).
\end{thm}

\begin{Ex}\label{ex:holoGConn}
	Let ${P}$ be a holomorphic principal bundle over a complex manifold $M$ with the structure group $K$ a connected complex Lie group.
	Suppose $M$ is acted on holomorphically by a connected complex Lie group $G$, 
	then there exists a holomorphic action Lie algebroid $\g\times M$ where $\g$ is the Lie algebra of~$G$.
	Biswas~et~al.~\cite{Bis15} constructed a holomorphic $\g$-Atiyah algebroid $\At_\g({P})$ and a holomorphic $\g$-Atiyah sequence:
	\begin{equation*}
		\xymatrix@C=2pc@R=2pc{
			0\ar[r] &  \Ad(P) \ar[r] & \At_\g({P}) \ar[r] &
			\g\times M \ar[r] & 0.
		}
	\end{equation*}
	If one defines a \textbf{holomorphic $\g$-connection} on $P$
	as a holomorphic vector-bundle map ${s\colon \g\times M \to \At_\g({P})}$ that splits the above sequence,
	then the existence of such a connection is obstructed by the Atiyah class
	$\alpha\in H^1_{\parbar} \big(M,(\g^\vee\times M)\otimes\Ad(P)\big)$ of the Lie triad $\big(T^{0,1}M\bowtie (\g\times M),\,T^{0,1}M,\,P\big)$.
\end{Ex}

\begin{Ex}
	A holomorphic Poisson manifold is a complex manifold $M$ with a holomorphic Poisson bivector $\pi \in \wedge^2 T^{1,0} M$, (i.e., $\pi$ satisfies $\parbar \pi =0$ and $[\pi,\pi]=0$ for the Schouten bracket.)
	
	The holomorphic cotangent bundle $\mathscr{T}^\vee M$ becomes a holomorphic Lie algebroid under the anchor 
	$$\pi^\#\colon (\mathscr{T}^\vee M)(U) \longto (\mathscr{T} M)(U), \quad \omega\longmapsto \pi(\omega,\cdot),$$
	and the Lie bracket 
	$$[\omega_1,\omega_2]_\pi=\L_{\pi^\# \omega_1} \omega_2- \L_{\pi^\# \omega_2} \omega_1- d(\pi(\omega_1,\omega_2)),\quad \forall\, \omega_i\in (\mathscr{T}^\vee M)(U),$$
	in any open neighbourhood $U\subseteq M$.
	Thus, the underling smooth bundle $(T^{1,0} M)^\vee$ of $\mathscr{T}^\vee M$  is matched with $T^{0,1}M$ to form a Lie algebroid $T^{0,1}M\bowtie (T^{1,0} M)^\vee$.
	
	The notion of cotangent connection introduced by Fernandes~\cite{Fer00} in the smooth setting certainly applies to the holomorphic setting.
	Let $\E$ be a holomorphic vector bundle on $M$.
	By a \textbf{holomorphic cotangent connection}, we mean a global holomorphic $\mathscr{T}^\vee M$-connection on $\E$,
	which is a holomorphic vector-bundle right splitting~$s$ to the holomorphic $\mathscr{T}^\vee M$-Atiyah sequence:
	\begin{equation*}
		\xymatrix@C=2.5pc@R=2pc{
			0\ar[r]&  \mathend(\E) \ar[r]& {\mathfrak{D}_{\mathscr{T}^\vee M}(\E)} \ar@<0.5ex>[r]^{\quad \sigma_{_{\mathscr{T}^\vee M}}}&
			\mathscr{T}^\vee M \ar[r] \ar@<0.5ex>@{-->}[l]^{\quad s}& 0,
		}
	\end{equation*}
	Here $\mathfrak{D}_{\mathscr{T}^\vee M}(\E)$ is the holomorphic Atiyah algebroid, and the existence of such an $s$ is obstructed by
	the Atiyah class $$\alpha\in H^1_{\parbar}\big(M,((\mathscr{T} M)^\vee)^\vee\otimes \End(\E)\big)=H^1_{\parbar}\big(M,\mathscr{T} M \otimes \End(\E)\big).$$

	In particular, if $\E$ is a holomorphic line bundle on a holomorphic Poisson manifold $M$ whose ${H^1_{\parbar}(M,\mathscr{T} M)}$ vanishes, then the existence of holomorphic cotangent connection on $\E$ is unobstructed.
	For instance, see~\cite{Hong}, any toric Fano manifold $M$ has the vanishing $H^1_{\parbar}(M,\mathscr{T} M)=0$, and can be equipped with toric holomorphic Poisson structures.
	Hence, holomorphic line bundles on a fano, toric holomorphic Poisson manifold always have holomorphic cotangent connections.
\end{Ex}

\vskip 10pt
\subsection{Equivariant Atiyah algebroids and invariant algebroid connections}\label{Sec:equivariantthings}

Let a Lie algebra $\g$ act on~$M$ via a Lie algebra morphism
$$\X\colon \g\longto \Gamma(TM),\, v \longmapsto \X_v.$$
It gives us an action Lie algebroid $\g\ltimes M $ on which the anchor for $V\in \Gamma\big(\g\ltimes M\big)=\C^\infty(M,\g)$ is defined at any point $m\in M$ as $\rho(V)(m)=\X_{V(m)}(m)$, and the Lie bracket $[V,W]\in \C^\infty(M,\g)$ for $V,\,W\in \C^\infty(M,\g)$ at $m$  is
\[
[V,W](m)=[V(m),W(m)]_\g+\big(\mathscr{L}_{\X_{V(m)}}W\big)(m)-\big(\mathscr{L}_{\X_{W(m)}}V\big)(m),
\]
where the first term takes Lie bracket on the $\g$-fiber over each $m$, and the last two terms are Lie derivatives on functions.

\begin{defn}\label{dfn:g-bundle}
	A vector bundle $E$ over $M$ is a \textbf{$\g$-equivariant vector bundle}
	if we can lift $\X$ to be a Lie algebra morphism $\X^E\colon \g\to \Gamma\big(\fD(E)\big)\colon v \mapsto \X^E_v$ that makes the following diagram of Lie algebras commute:
	\begin{equation*}
		\xymatrix@C=2.5pc@R=2pc{
			\ar@{}[dr]|(0.67){\scalebox{1.2}{$\circlearrowright$}}
			&\Gamma(\fD(E)) \ar[d]^{\sigma} \\
			 \g\ar[ur]^{\X^E}\ar[r]_{\X\quad}
			&\Gamma(TM).
		}
	\end{equation*}
	Such an $\X^E$ is called a \textbf{$\g$-action on $E$ that lifts $\X$}.
\end{defn}

Note that $\X^E$ extends to a $(\g\ltimes M)$-connection $D\colon \Gamma(\g\ltimes M) \times \Gamma(E)\to \Gamma(E)$,
defined for $V\in\Gamma(\g\ltimes M)$, $e\in\Gamma(E)$ by setting $D_{V}(e)$ at each point $m\in M$ that
$$\big(D_{V}(e)\big)(m) \triangleq \big(\X^E_{V(m)}(e)\big)(m).$$
The initial requirement that $\X^E\colon \g\to \Gamma\big(\fD(E)\big)$ preserves Lie brackets means the $(\g\times M)$-connection $D$ is flat, and vice versa.
Thus, we have:

\begin{lem}\label{lem:g-bundle}
	A $\g$-equivariant vector bundle is equivalent to a $(\g\ltimes M)$-representation.
\end{lem}

For a Lie algebroid $(L,[~\cdot,\cdot~]_{_L},\rhoL)$ over $M$,
we consider its infinitesimal automorphisms.

\begin{Def}\label{dfn:DerL}
	A \textbf{Lie algebroid derivation} on $L$ is a differential operator
	$\delta\colon\Gamma(L)\to \Gamma(L)$ equipped with a symbol $X_\delta\in \Gamma(TM)$, which satisfies the following conditions for all $f \in \C^\infty(M)$, $l_1,l_2 \in \Gamma(L)$:
	\begin{enumerate}[(i)]
		\item $\delta(f l_2) = f\delta  ( l_1)+ \big(X_\delta  ( f)\big)  l_1$,
		\item $\rhoL\big(\delta (l_1)\big) = [X_\delta,  \rhoL(l_1)]$,
		\item $\delta  [l_1,l_2]_{_L} = [\delta (l_1),l_2]_{_L} +[l_1,\delta  (l_2)]_{_L}$.
	\end{enumerate}
\end{Def}

The set of all Lie algebroid derivations on $L$ is denoted by $\Der(L)$, and has the symbol map
$$\widehat{\rho}_{_L}\colon \Der(L)\longto \Gamma(TM)\colon \delta \longmapsto X_\delta.$$

\begin{defn}\label{dfn:g-algebroid}
	A Lie algebroid $L$ over $M$ is a \textbf{$\g$-equivariant Lie algebroid}
	if we can lift $\X$ to be a Lie algebra morphism $\X^L\colon \g\to \Der(L)\colon v \mapsto \X^L_v$
	that makes the following diagram of Lie algebras commute:
	\begin{equation*}
		\xymatrix@C=2.5pc@R=2pc{
			\ar@{}[dr]|(0.67){\scalebox{1.2}{$\circlearrowright$}}
			&\Der(L)\ar[d]^{\widehat{\rho}_{_L}} \\
			\g\ar[ur]^{\X^L} \ar[r]_{\X\quad}
			&\Gamma(TM).
		}
	\end{equation*}
	Such an $\X^L$ is called a \textbf{$\g$-action on $L$ that lifts $\X$}.
\end{defn}

\begin{rmk}
	This definition of a $\g$-action on Lie algebroid $L$ is a special case of Mackenzie's derivative representation of a Lie algebroid~\cite{Mac05},
	and a special case of Mehta--Zambon's $L_\infty$ action of a Lie algebra~\cite{MehtaZambon} and Zambon--Zhu's strictly Lie 2-algebra action~\cite{ZambonZhu}.
\end{rmk}

As we have observed in Lemma~\ref{lem:g-bundle}, the $\g$-action $\X^L$ on $L$ extends to a $(\g\ltimes M)$-representation $D$ on~$L$.
On the other hand, there is an $L$-representation on $(\g\ltimes M)$ given for $l\in \Gamma(L),\, V\in \Gamma(\g\ltimes M)$ using the Lie derivative:
$$\mathscr{L}_l V\triangleq\mathscr{L}_{\rhoL(l)} V,$$
which is actually determined by its evaluation on constant sections $ \underline{v}\in \Gamma(\g\ltimes M)$ with $\underline{v}\equiv v\in \g$ being that $\mathscr{L}_l (\underline{v})=0$.
These two mutual representations $D$ and $\mathscr{L}$ build a matched pair between $\g\ltimes M$ and $L$.

\begin{prop}\label{prop:g-algebroid}
	A $\g$-equivariant structure on a Lie algebroid $L$ over $M$ is equivalent to a matched sum {$(\g\ltimes M) \bowtie L$} structure whose Bott representation of $L$ on $(\g\ltimes M)$ is the Lie derivative $\mathscr{L}$.
\end{prop}
\begin{proof}
	Note that conditions (i)-(iii) in Lemma~\ref{lem:matchedPair} for the two mutual representations $D,\,\mathscr{L}$ to form a matched pair are tensorial.
	We only need to evaluate them for sections of $L$ and constant sections of $\g\ltimes M$.
	It amounts to requiring that $D_{\underline{v}}$ belong to $\Der{L}$ for any constant section $\underline{v}$ of $\g\ltimes M$,
	hence providing a $\g$-action on the Lie algebroid $L.$
\end{proof}

\begin{thm}\label{thm:equivConn}
	Let $\g$ act on $M$ via $\X$, and suppose it can be lifted to act on a Lie algebroid $L$ via $\X^L$ and on a vector bundle $E$ via $\X^E$.
	Then, 
	\begin{enumerate}[(a)]
        \item we have a matched-type Lie triad $\big((\g\ltimes M) \bowtie L,\,\g\ltimes M,\,E\big)$;
		\item the Atiyah bundle $\fD_L(E)$ is a $\g$-equivariant Lie algebroid, and the Atiyah sequence:
		\begin{equation}
			\xymatrix@C=2pc@R=2pc{
				0\ar[r]&  \mathend(E) \ar[r]^{{\iota}}& {\LdiffE} \ar[r]^{\quad \sigma_{_L}}&
				L\ar[r]& 0,
			}
		\end{equation}
		is simultaneously a sequence of $\g$-equivariant vector bundles and also a sequence of Lie algebroids.
	\end{enumerate} 
\end{thm}
\begin{proof}
    This theorem results from Proposition~\ref{prop:g-algebroid} and Theorem~\ref{thm:matchedAtExt}.
\end{proof}

By Theorem~\ref{thm:AtClass}, the Atiyah class $\alpha$ obstructs the splitting of the  above sequence of $\g\ltimes M$-representations,
i.e., the existence of $(\g\ltimes M)$-compatible $L$-connections on $E$.
So, we make the following definition.

\begin{defn}
	Under the assumptions of Theorem~\ref{thm:equivConn}, an $L$-connection $\nabla$ is said to be \textbf{$\g$-invariant}, if the Atiyah cocycle vanishes, i.e.,
	for $V\in \Gamma(\g\ltimes M),\,l\in \Gamma(L)$
	\[
	R^\nabla_{(\g\ltimes M)\otimes L}(V,l)=\X^E_V\nabla_l-\nabla_l\X^E_V - \nabla_{[V,l]}=0.
	\]
	Since $R^\nabla_{(\g\ltimes M)\otimes L}$ is tensorial, this condition only needs to hold for constant sections $\underline{v}\equiv v\in \g$ where $[\underline{v},l]=\X^L_{\underline{v}} l-\L_l \underline{v}=\X^L_{\underline{v}} l$:
	\[
	R^\nabla_{(\g\ltimes M)\otimes L}(\underline{v},l)=\X^E_{\underline{v}}\nabla_l-\nabla_l\X^E_{\underline{v}} - \nabla_{\X^L_{\underline{v}} l}=0.
	\]
\end{defn}

If the actions of a Lie algebra $\g$ on $L$ and $E$ are actually deduced from smooth actions of a connected Lie group $G$ preserving the Lie algebroid structure on $L$ and the vector bundle structure on $E$,
then by exponentiation, the $\g$-invariance of an $L$-connection $\nabla$ on $E$ is the same as the $G$-invariance of $\nabla$ in the usual sense that $g^*\nabla = \nabla$ for any $g\in G$. Therefore, we have:
\begin{prop}
	If $G$ is compact, then any $L$-connection $\nabla$ on $E$ can be integrated with respect to the Haar measure of $G$ to produce an $L$-connection that is $G$-invariant and hence also $\g$-invariant.
\end{prop}

However, without assuming $G$ being compact, the existence of a $G$-invariant (hence $\g$-invariant) connection is not always guaranteed.
For instance, Nomizu~\cite{Nom54} and Wang~\cite{Wan58} have studied the necessary and sufficient conditions for the existence of invariant connections on homogeneous spaces and their principal bundles.
Now, we can use the $\g$-equivariant Atiyah sequence to give a concise description of their classic results, also see~\cite{KT71}.

\vskip 10pt
\subsection{Invariant connections on homogeneous spaces}\label{subsec:homogeneous}

	Let $M=G/H$ be a homogeneous space of a connected Lie group $G$ with a closed isotropy subgroup $H$, without assuming $G$ being compact.
	A principal $K$-bundle $\pi\colon P \to G/H$ with the Lie group $K$ acting on $P$ from the right is \textbf{($G$-)homogeneous}
	if the $G$-action on $G/H$ can be lifted on $P$ and commutes with the $K$-action, 
	i.e., 
	$$\pi(gp)=g\pi(p),\quad g(pk)=(gp)k,\qquad\forall\, p\in P,\,g\in G,\,k\in K.$$
	
	The isotropy group $H$ acts on the $K$-fiber $\pi^{-1}(eH)$,
	so that it induces a group homomorphism $\varphi\colon H\to K$, called isotropy homomorphism, and gives the identification
	\[
	P=G\times_H K=(G\times K)/\big\{(g,k)\sim(gh^{-1},\varphi(h)k)\big\}.
	\]
	
	The Atiyah algebroid $(TP)/K$ over $G/H$ is endowed with a natural $G$-action and hence becomes a homogeneous vector bundle over $G/H$.
	To describe it, note the fiber of $(TP)/K$ at $eH$ is
	\[
	\g\times_\h \k=(\g\times \k)/\big\{(X,\,Z)\sim(X-Y,\,\rd\varphi(Y)+Z)\big\},
	\]
	for all $X\in \g,\, Y\in \h,\, Z\in \k$ and $\rd\varphi\colon \h \to \k$ is the Lie algebra morphism by differentiating $\varphi$.
	Denote an element of $\g\times_\h \k$ by $(X,\,Z)_\h$ where $(X,\,Z)\in \g\times \k$,
	the isotropy representation of $H$ on $\g\times_\h \k$ can be verified as
	\[
	h\cdot (X,\,Z)_\h = \big(\Ad_h(X),\,\Ad_{\varphi(h)}(Z)\big)_\h,\quad \forall\,h\in H.
	\]
	Hence, $(TP)/K=G\times_H (\g\times_\h \k)$, and similarly, one checks that each item in the Atiyah sequence of the principal $K$-bundle $P$ takes the format of a homogeneous vector bundle:
	\begin{equation}\label{seq:GhomoAtiyah}
		\begin{split}
			\xymatrix@C=2pc@R=2pc{
				0\ar[r]&  P\times_K \k  \ar[r] \ar@{=}[d] & (TP)/K \ar@{=}[d] \ar[r]&
				T(P/K)\ar[r] \ar@{=}[d]& 0\,\\
				0\ar[r]&  G\times_H \k   \ar[r]
				& G\times_H (\g\times_\h \k) \ar[r]&
				G\times_H (\g/\h) \ar[r]& 0,
			}
		\end{split}
	\end{equation}
	where the isotropy representations of $H$ on $\k$ and $\g/\h$ are respectively given by $\Ad_{\varphi(h)}$ and $\Ad_h$ for $h\in H$.
	Under the corresponding infinitesimal actions of $\g$,
	this is a sequence of $\g$-equivariant Lie algebroids by Theorem~\ref{thm:equivConn}.
	
	Moreover, we clearly have the following series of one-to-one correspondences:
	\begin{equation*}
		\begin{aligned}
			&\{\,\g\mbox{-invariant connections of }P\,\} &&\\
			\overset{1:1}{\longleftrightarrow} \,\,& \{\,\g\mbox{-equivariant vector-bundle left splittings of Sequence~\eqref{seq:GhomoAtiyah}}\,\}
			&& (\,\mbox{by Proposition~\ref{prop:AtiyahCocyle=0}}\,)\\
			\overset{1:1}{\longleftrightarrow} \,\,& \{\,G\mbox{-equivariant vector-bundle left splittings of Sequence~\eqref{seq:GhomoAtiyah}}\,\}
			&& (\,\mbox{$G$ is connected}\,)\\
			\overset{1:1}{\longleftrightarrow} \,\,&
			\left\{
 			\begin{split}
				\,H\mbox{-equivariant } &\mbox{vector-space left splittings of the sequence:}\\ 
				& 0\to\k\to\g\times_\h \k \to \g/\h \to 0.\hfill
			\end{split}\,\right\}
			&& (\,\mbox{Sequence~\eqref{seq:GhomoAtiyah} is $G$-homogeneous}\,)
		\end{aligned}
	\end{equation*}
	The map $\k\to\g\times_\h \k$ is given by $\k\ni Z\mapsto (0,Z)_\h$.
	A left inverse $\g\times_\h \k \to \k$ then requires $(0,Z)_\h \mapsto Z$ together with a map $\phi$ that takes $\{(X,0)_\h\mid X\in \g\}$ to $\k$,
	satisfying the compatibility condition that for $Y\in \h$, the element $(0,0)_\h=(-Y,\rd\varphi(Y))_\h$ should be sent to $0=\phi(-Y,0)_\h+\rd \varphi(Y)$.
	
	Diagrammatically, since $\g\times_\h \k$ is a push-out from $\g$ and $\k$ with respect to $\h$,
	a left inverse $\g\times_\h \k \to \k$ to $\k\to\g\times_\h \k$ is equivalent to a map $\phi\colon \g\to \k$ that makes the following upper triangle commute:
	\begin{equation}
		\xymatrix@C=2.5pc@R=2pc{
			\h
			\ar[r]^{\rd \varphi}
			\ar[d]
			\ar@{}[dr]|(0.75){\text{\pigpenfont R}}|(0.3){\scalebox{1.2}{$\circlearrowleft$}}& \k\, \ar[d]
			\\
			\g \ar@{-->}[ur]|{\scalebox{1}{$\phi$}}
			\ar[r]
			&\g\times_\h \k.}
	\end{equation}
	
	Therefore, we obtain Wang's~\cite{Wan58} classic result that
	\begin{equation*}
		\begin{aligned}
			&\{\,\g\mbox{-invariant connections of }P\,\} &&\\
			\overset{1:1}{\longleftrightarrow} \,\,& \{\,H\mbox{-equivariant vector-space maps $\phi\colon \g\to \k$ such that $\phi|_\h=\rd \varphi$}\,\}.
		\end{aligned}
	\end{equation*}
	In particular, consider $P=G$ as a principal $H$-bundle over $G/H$,
	that is, $K=H$ and the isotropy morphism is $\varphi=\id_\h\colon\h\to\h$. 
	Then, by the above correspondence, we see that
	\begin{equation*}
		\begin{aligned}
			&\{\,\g\mbox{-invariant connections of the principal $H$-bundle }G\to G/H\,\} &&\\
			\overset{1:1}{\longleftrightarrow} \,\,& \{\,H\mbox{-equivariant vector-space surjections $\phi\colon \g\to \h$ such that $\phi|_\h=\id_\h$}\,\}\\
			\overset{1:1}{\longleftrightarrow} \,\,& \{\,H\mbox{-equivariant vector-space decomposition $\g \cong \h\oplus \m$}\,\}.
		\end{aligned}
	\end{equation*}
	The non-emptiness of the last set is exactly the definition of $G/H$ being reductive.
	
	Finally, suppose $G/H$ is reductive by fixing an $H$-equivariant decomposition $\g \cong \h\oplus \m$, or equivalently, an $H$-equivariant surjection $\phi_0\colon \g\to \h$ such that $\phi|_\h=\id_\h$.
	Then, for any homogeneous principal $K$-bundle $\pi\colon P \to G/H$ with an isotropy homomorphism $\varphi\colon H\to K$, the set
	\[
	\{\,H\mbox{-equivariant vector-space maps $\phi\colon \g\to \k$ such that $\phi|_\h=\rd \varphi$}\,\}
	\]
	is nonempty,
	because the composition $\rd \varphi \circ \phi_0\colon \g\to\h\to\k$ provides one element in that set.
	The $\g$-invariant connection of $P$ corresponding to $\rd \varphi \circ \phi_0$ is called the canonical connection on $P$ with respect to the fixed decomposition $\g \cong \h\oplus \m$.
	Moreover, the above set is in one-to-one correspondence to
	\[
	\{\,H\mbox{-equivariant vector-space maps $\tilde{\phi}\colon \m\to \k$} \,\},
	\]
	where the canonical connection is represented by the zero map $\tilde{\phi}=0$.

\vskip 25pt
\appendix
\hypertarget{App:A}{}
\section*{Appendix A. Functorial relations among Atiyah bundles}

\renewcommand{\thesection}{A}

\setcounter{equation}{0}

\vskip 10pt

In this appendix, we supply various relations among the Atiyah bundles with respect to $A,\,L,\,B$.

To start with, there is another way to establish $B$-Atiyah sequence~\eqref{seq:AtiyahDBE} besides the proof of Theorem~\ref{thm:DBE}.
The composition $s_{_{A,L}}=f_{_{A,L}}\circ s_{_A}$ in~\eqref{eq:A-DLE} induces a morphism between Lie pairs $\big(\fD_A(E),A\big)$ and $\big(\fD_L(E),A\big)$, 
and hence by Proposition~\ref{prop:morphismOfBottConn} a morphism between their quotient sequences.
Applying the Snake Lemma, we have a commutative diagram of vector bundles:
\begin{equation}\label{dia:snakeOfLiePair2}
	\begin{split}
		\xymatrix@C=2pc@R=2pc{
			&&0\ar[d]&0\ar[d]&\\
			0 \ar[r]&A\ar[r]^{s_{_A}\quad}\ar@{=}[d]
			&\mathfrak{D}_A(E)\ar[r]^{\theta_{_A}}\ar[d]^{f_{_{A,L}}}
			&\End(E)\ar[r]\ar[d]^{\bar{\iota}_{_{B}}}&0\,\\
			0 \ar[r]&A\ar[r]^{s_{_{A,L}}\quad}
			&\mathfrak{D}_L(E)\ar[r]^{P_{_B}}\ar[d]^{\pr_{_B}\circ \sigma_{_L}}
			&\BdiffE\ar[r]\ar[d]^{\bar{\sigma}_{_{B}}}&0.\\
			&&B\ar@{=}[r]\ar[d]&B\ar[d]&\\
			&&0&0&
		}
	\end{split}
\end{equation}
The second vertical short exact sequence has morphisms $\bar{\iota}_{_{B}}$ 
and $\bar{\sigma}_{_{B}}$ actually the same as ${\iota}_{_{B}}$ and ${\sigma}_{_{B}}$ by diagram chasing,
hence is exactly $B$-Atiyah sequence~\eqref{seq:AtiyahDBE}.

Collecting from Diagrams~\eqref{dia:snakeOfLiePair1},~\eqref{dia:snakeOfLiePair2} the short exact sequences involving $A$,
we get a diagram of morphisms among three quotient sequences of Lie pairs:
\begin{equation}\label{Dia:LA-Atiyah1}
	\begin{split}
		\xymatrix@C=2.5pc@R=2pc{
			0\ar[r]&A\ar[r]^{s_{_A}\quad}\ar@{=}[d]
			&\fD_A(E)\ar[r]^{\theta_{_A}}\ar[d]_{f_{_{A,L}}}&\End(E)\ar[d]^{\iota_{_B}}\ar[r]&0\,\\
			0\ar[r]&A\ar[r]^{s_{_{A,L}}\quad}\ar@{=}[d]
			&\mathfrak{D}_L(E)\ar[r]^{P_{_B}}\ar[d]_{\sigma_{_L}}
			&\BdiffE\ar[r]\ar[d]^{\sigma_{_B}}&0\,\\
			0\ar[r]&A\ar[r]^{i_{_A}\quad}&L\ar[r]^{\pr_{_B}}&B\ar[r]&0,
		}
	\end{split}
\end{equation}
where the last vertical sequence is the Atiyah sequence relative to $(L,A)$.

Collecting from Diagrams~\eqref{dia:snakeOfAtiyah},~\eqref{dia:snakeOfLiePair1} the short exact sequences involving $\End(E)$, we get a diagram of morphisms among three Atiyah sequences, cf. Diagram~\eqref{Dia:DALE-pullback}:
\begin{equation}\label{Dia:LA-Atiyah2}
	\begin{split}
		\xymatrix@C=2.5pc@R=2pc{
			0 \ar[r]
			&\mathend(E) \ar[r]^{{\iota}} \ar@{=}[d]
			&\mathfrak{D}_A(E) \ar@<0.5ex>[r]^{\quad\sigma_{_A}}\ar[d]_{f_{_{A,L}}}
			&A\ar[r]\ar[d]^{i_A} \ar@<0.5ex>[l]^{\quad s_{_A}}
			&0\,\\
			0 \ar[r]
			&\mathend(E) \ar[r]^{{\iota}} \ar@{=}[d]
			&\LdiffE \ar[r]^{\quad\sigma_{_L}} \ar[d]_{P_{_{B}}}
			&L \ar[r]  \ar[d]^{{\rm pr}_{_B}}
			&0\,\\
			0 \ar[r]
			&\mathend(E) \ar[r]^{\iota_{_B}}
			&\BdiffE \ar[r]^{\quad\sigma_{_B}}
			&B \ar[r] &0,
		}
	\end{split}
\end{equation}
where the last vertical sequence is the  quotient sequence of $(L,A)$.

Collecting from Diagrams~\eqref{dia:snakeOfAtiyah},~\eqref{dia:snakeOfLiePair2} the short exact sequences involving $B$,
we get a diagram of morphisms among different types of short exact sequences:
\begin{equation}\label{Dia:LA-Atiyah3}
	\begin{split}
		\xymatrix@C=2.5pc@R=2pc{
			0\ar[r]&A\ar[r]^{i_{_A}\quad}&L\ar[r]^{\pr_{_B}}&B\ar[r]&0
			\\
			0 \ar[r]&\fD_A(E)\ar[r]^{ f_{_{A,L}}}
			\ar[u]^{\sigma_{_A}}\ar[d]_{\theta_{_A}}
			&\mathfrak{D}_L(E)\ar[r]^{\quad{\rm pr}_{_B}\circ \sigma_{_L}}
			\ar[u]^{\sigma_{_L}}\ar[d]_{P_{_B}}
			&B\ar[r]\ar@{=}[u]\ar@{=}[d]&0\\
			0 \ar[r]
			&\mathend(E) \ar[r]^{\iota_{_B}}
			&\BdiffE \ar[r]^{\quad\sigma_{_B}}
			&B \ar[r] &0.
		}
	\end{split}
\end{equation}
where the first two horizontal sequences are quotient sequences of Lie pairs 
while the third horizontal one is the relative Atiyah sequence.
One should beware of the different directions of vertical arrows.

All these bundles and arrows can be assembled to be a centered hexagon:
\begin{equation}\label{Dia:Hexagon}
	\begin{split}
		\xymatrix@C=1pc@R=4.8pc{
			&L\ar[rr]^{\pr_{_B}}&&B&\\
			A\ar[ru]^{i_{_A}}\ar[rr]^{s_{_{A,L}}\quad}\ar@<-0.5ex>[dr]_{s_{_A}}
			&&\fD_{L}(E)\ar[lu]^{\sigma_{_L}}
			\ar[ru]^{\pr_{_B}\circ\sigma_{_L}}\ar[rr]^{P_{_B}}
			&&\BdiffE\ar[lu]_{\sigma_{_B}}\\
			&\fD_A(E)\ar[ru]^{f_{_{A,L}}}\ar@<-0.5ex>[lu]_{\sigma_{_A}}
			\ar@<-0.5ex>[rr]_{\theta_{_A}}
			&&\End(E)\ar[lu]^{{\iota}}\ar[ru]_{\iota_{_B}}
			\ar@<-0.5ex>[ll]_{{\iota}}&
		}
	\end{split}
\end{equation}
where we could look for short exact sequences either starting with $A$ or $\End(E)$, or ending with $B$.
As for the commutativity of the triangles,
we should use $s_{_A}$ instead of $\sigma_{_A}$ in the lower left triangle and ${\iota}$ instead of $\theta_{_A}$ in the bottom triangle.

\end{document}